\documentclass[12pt]{amsart}
\usepackage{amsfonts}
\usepackage{amsmath}
\usepackage{amsxtra}
\usepackage{amssymb,latexsym}
\usepackage[mathcal]{eucal}
\usepackage{amscd}

\makeindex

\newtheorem{theo}{{\bfseries Theorem}}[section]
\newtheorem{prop}[theo]{{\bfseries Proposition}}
\newtheorem{lem}[theo]{{\bfseries Lemma}}
\newtheorem{cor}[theo]{{\bfseries Corollary}}
\newtheorem{df}[theo]{{\bfseries Definition}}

\newtheorem{ex}[theo]{{\bfseries Example}}
\newtheorem{exes}[theo]{{\bfseries Examples}}

\newtheorem{queses}[theo]{{\bfseries Questions}}

\def \ol {\overline}
\def \N {\mathbb N}

\def \Z {\mathbb Z}
\def \R {\mathbb R}

\def \A {\mathcal A}
\def \B {\mathcal B}
\def \CC {\mathcal C}

\def \L {\mathcal L}

\def \NN {\mathcal N}

\def \T {\mathcal T}
\def \G {\mathcal G}

\def \U {\mathcal U}
\def \V {\mathcal V}
\def \W {\mathcal W}
\def \a {\alpha }
\def \b {\beta}
\def \lm {\lambda}
\def \ep {\epsilon}
\def \om {\omega}
\def \d {\delta}

\def \s {\sigma}

\usepackage{amssymb,latexsym}
\usepackage[mathcal]{eucal}
\usepackage{amscd}
\usepackage{amsmath}
\usepackage{graphicx}
\usepackage{amsfonts}
\usepackage{amscd}

\numberwithin{equation}{section}

\begin{document}

\title{\bfseries  Chain Recurrence for General Spaces}
\vspace{1cm}
\author{Ethan Akin and Jim Wiseman}
\address{Mathematics Department \\
    The City College \\ 137 Street and Convent Avenue \\
       New York City, NY 10031, USA     }
\email{ethanakin@earthlink.net}

\address{Department of Mathematics \\
Agnes Scott College \\ 141 East College Avenue \\ Decatur, GA 30030, USA }
\email{jwiseman@agnesscott.edu}

\date{July, 2017}

\vspace{.5cm} \maketitle

\emph{In Memory of John Mather}

\tableofcontents

\newpage

\section{ Introduction}
\vspace{.5cm}

Let $f$ be a continuous map on a compact metric space $(X,d)$.  If $\ep \geq 0$ then a sequence
$\{ x_0, \dots, x_n \}$ with $n \geq 1$ is an $\ep$ chain for $f$ if $\max_{i=1}^n \ d(f(x_{i-1}),x_i) \leq \ep$ and
a strong $\ep$ chain for $f$ if $\Sigma_{i=1}^n \ d(f(x_{i-1}),x_i)$ $ \leq \ep$. Thus, a $0$ chain is just an initial piece
of an orbit sequence.

The Conley chain relation $\CC f$
consists of those pairs $(x,y) \in X \times X$ such that there is an $\ep$ chain with $x_0 = x$ and $x_n = y$ for every
$\ep > 0$. The Easton, or Aubry-Mather, strong chain relation $\A_d f$
consists of those pairs $(x,y) \in X \times X$ such that there is a strong
$\ep$ chain with $x_0 = x$ and $x_n = y$ for every
$\ep > 0$. As the notation indicates, $\CC f$ is independent of the choice of
metric, while $\A_d f$ depends on the metric.
See \cite{C} and \cite{E}.

Fathi and Pageault have studied these matters using what they call \emph{barrier functions}, \cite{P}, \cite{FP} and their work has
been sharpened by Wiseman \cite{W16}, \cite{W17}.  $M_d^f(x,y)$  is the infimum of the $\ep$'s
such that there is an $\ep$ chain from $x$ to $y$ and $L_d^f(x,y)$  is the infimum of the $\ep$'s
such that there is a strong $\ep$ chain from $x$ to $y$. Thus, $(x,y) \in \CC f$ iff $M_d^f(x,y) = 0$ and
 $(x,y) \in \A_d f$ iff $L_d^f(x,y) = 0$.

Our purpose here is to extend these results in two ways.

First, while our interest focuses upon homeomorphisms  or continuous maps, it is convenient, and easy, to extend the
results to relations, following \cite{A93}.

A relation\index{relation} $f : X \to Y$ is just a subset of $X \times Y$ with
$f(x) = \{ y \in Y : (x,y) \in f \}$ for $x \in X$,
and let $f(A) = \bigcup_{x \in A} \ f(x)$ for $A \subset X$.\index{f@$f(A)$}
So $f$ is a mapping when $f(x)$ is a singleton set for every $x \in X$, in which case we will use the notation
$f(x)$ for both the singleton set and the point contained therein. For
example, the identity map on a set $X$ is $1_X = \{ (x,x) : x \in X \}$.
If $X$ and $Y$ are topological spaces then $f$ is a closed relation when
it is a closed subset of $X \times Y$ with the product
topology.

The examples $\CC f$ and $\A_d f$ illustrate how relations arise naturally in dynamics.

For a relation $f : X \to Y$ the \emph{inverse relation}\index{relation!inverse}\index{inverse relation}
$f^{-1} : Y \to X $ is $ \{ (y,x) : (x,y) \in f \}$.
Thus, for $B \subset Y$, $f^{-1}(B) = \{ x : f(x) \cap B \not= \emptyset \}.$
We define $f^*(B)  = \{ x : f(x) \subset  B \}$.\index{f@$f^*(B)$} These are equal when $f$ is a map.

If $f : X \to Y$ and $g : Y \to Z$ are relations then
the \emph{composition}\index{relation!composition}  $g \circ f : X \to Z $ is
$\{ (x,z) : $ there exists $y \in Y$ such that $(x,y) \in f$ and $(y,z) \in  g \}$.  That is, $g \circ f$ is the image
of $(f \times Z) \cap (X \times g)$ under the projection $\pi_{13} : X \times Y \times Z \to X \times Z$.
As with maps, composition of relations is clearly associative.

 The \emph{domain}\index{relation!domain}\index{domain} of a relation $f : X \to Y$ is
\begin{equation}\label{introeq0}
Dom(f) \ = \ \{ x : f(x) \not= \emptyset \} \ = \ f^{-1}(Y).
\end{equation}
We call a relation \emph{surjective} if $Dom(f) = X$ and $Dom(f^{-1}) = Y$, i.e.\ $f(X) = Y$ and $f^{-1}(Y) = X$.
\index{relation!surjective}\index{surjective relation}

If $f_1 : X_1 \to Y_1$ and $f_2 : X_2 \to Y_2$ are relations, then the product relation \index{relation!product}\index{product relation}
$f_1 \times f_2 : X_1 \times X_2 \to Y_1 \times Y_2$ is $\{ ((x_1,x_2),(y_1,y_2)) : (x_1,y_1) \in f_1, (x_2,y_2) \in f_2 \}$.

We call $f$ a \emph{relation on $X$} when $X = Y$.  In that case, we define, for $n \geq 1$
$f^{n+1} = f \circ f^n = f^n \circ f$ with $f^1 = f$.  By definition, $f^0 = 1_X$ and $f^{-n} = (f^{-1})^n$.
If $A \subset X$, then
$A$ is called $f$ \emph{$^+$invariant} if $f(A) \subset A$ and $f$ \emph{invariant} if $f(A) = A$.\index{subset!invariant}\index{subset!$^+$invariant}
In general, for $A \subset X$, the \emph{restriction to $A$} is  $f|A = f \cap (A \times A)$.\index{restriction}\index{f@$f"|A$}
If $u$ is a real-valued function on $X$ we will also write $u|A$ for the restriction of $u$ to $A$, allowing context to determine
which meaning is used.

The \emph{cyclic set}\index{relation!cyclic set} \index{cyclic set $"|f"|$}$|f|$ of a relation $f$ on $X$ is $\{ x \in X : (x,x) \in f \}$.

A relation $f$ on $X$ is \emph{reflexive}\index{relation!reflexive}\index{reflexive relation} if $1_X \subset f$,
\emph{symmetric}\index{relation!symmetric}\index{symmetric relation}
if $f^{-1} = f$ and \emph{transitive}\index{relation!transitive}\index{inverse relation} if
$f \circ f \subset f$.

If $d$ is a pseudo-metric on a set $X$ and $\ep > 0$, then $V^d_{\ep} = \{ (x,y) : d(x,y) < \ep \}$ and
 $\bar V^d_{\ep} = \{ (x,y) : d(x,y) \leq \ep \}$.  Thus, for $x \in X$, $V^d_{\ep}(x)$ (or $\bar V^d_{\ep}(x)$)
is the open (resp.\  closed) ball centered at $x$ with radius $\ep$.\index{v@$V^d_{\ep}$}\index{v@$\bar V^d_{\ep}$}

A pseudo-ultrametric\index{pseudo-ultrametric}\index{ultrametric} $d$ on $X$ is a pseudo-metric
with the triangle inequality strengthened to $d(x,y) \leq \max(d(x,z),d(z,y))$
for all $z \in X$. A pseudo-metric $d$ is a pseudo-ultrametric iff the relations $V^d_{\ep}$ and $\bar V^d_{\ep}$ are equivalence relations for all $\ep > 0$.

 If $(X_1,d_1)$ and $(X_2,d_2)$ are pseudo-metric spaces then the product $(X_1 \times X_2, d_1 \times d_2)$
is defined by
$$d_1 \times d_2((x_1,x_2),(y_1,y_2)) \ = \ \max (d_1(x_1,y_1),d_2(x_2,y_2)).$$
Thus,
$V^{d_1 \times d_2}_{\ep} = V^{d_1}_{\ep} \times V^{d_2}_{\ep}$ and
$\bar V^{d_1 \times d_2}_{\ep} = \bar V^{d_1}_{\ep} \times \bar V^{d_2}_{\ep}$.

Throughout this work, all pseudo-metrics are assumed bounded. For example, on $\R$ we use
$d(a,b) = \min(|a - b|,1)$. Thus, if $A$ is a non-empty subset of
$X$ the \emph{diameter} $diam (A) = \sup \{ d(x,y): x,y \in A \}$\index{diameter}\index{d@$diam(A)$}
is finite.

For metric computations, the following will be useful.

\begin{lem}\label{introalglem} Let $a_1,a_2,b_1,b_2 \in R$. With $a \vee b = \max(a,b)$ and $a \wedge b = \min(a,b)$:
\begin{align}\label{introeq}
\begin{split}
|a_1 \vee b_1 - a_2 \vee b_2|, |a_1 \wedge b_1 - a_2 &\wedge b_2| \ \leq \ |a_1 - a_2| \vee |b_1 - b_2|. \\
(a_1 \vee b_1) \wedge (a_2 \vee b_1)   \wedge (a_1 \vee b_2)  &\wedge
(a_2 \vee b_2)  \ = \  (a_1  \wedge a_2)  \vee (b_1  \wedge b_2).
\end{split}
\end{align}
\end{lem}

{\bfseries Proof:} First, we may assume without loss of generality that $a_1 \vee b_1 \geq a_2 \vee b_2 = a_2$ and
so that $a_2 \geq b_2$.  If $a_1 \vee b_1 = a_1$ then $|a_1 \vee b_1 - a_2 \vee b_2| = a_1 - a_2$.
If $a_1 \vee b_1 = b_1$ then $|a_1 \vee b_1 - a_2 \vee b_2| = b_1 - a_2 \leq b_1 - b_2$.
For the $\wedge$ estimate, observe that $a \wedge b = - (-a)\vee(-b)$.

For the second, factor out $b_1$ and $b_2$ to get $(a_1 \vee b_1) \wedge (a_2 \vee b_1)  = (a_1 \wedge a_2) \vee b_1$,
and $(a_1 \vee b_2)  \wedge (a_2 \vee b_2) = (a_1 \wedge a_2) \vee b_2$. Then factor out $a_1 \wedge a_2$.

$\Box$ \vspace{.5cm}

The other extension is to non-compact spaces. This has been looked at in the past, see \cite{H} and \cite{P}.
However, the natural setting for the theory is that of uniform spaces as described in \cite{K} and \cite{B}, and reviewed in
Appendix B below.

A uniform structure $\U$ on a set $X$ is a collection of relations on $X$ which satisfy various
axioms so as to generalize the notion of metric space. To be precise, a $\U$ is a uniformity when
\begin{itemize}
\item $1_X \subset U$ for all $U \in \U$.
\item $U_1, U_2 \in \U$ implies $U_1 \cap U_2 \in \U$.
\item If $U \in \U$ and $W \supset U$, then $W \in \U$.
\item $U \in \U$ implies $U^{-1} \in \U$.
\item If $U \in \U$, then there exists $W \in \U$ such that $W \circ W \subset U$.
\end{itemize}
The first condition says that the relations are reflexive and the next two imply that they form a filter.\index{filter}

A uniformity $\U$ is equivalently given by its gage $\Gamma(\U)$, the set of  pseudo-metrics $d$ on $X$ (bounded by stipulation) with
the metric uniformity $\U(d)$, generated by $\{ V^d_{\ep} : \ep > 0\} $,  contained in $\U$.  The use of covers in \cite{P} and continuous
real-valued functions in \cite{H} are equivalent to certain choices of uniformity. To a uniformity there is an associated topology and
we say that $\U$ is compatible with a topology on $X$ if the uniform topology agrees with the given topology on $X$.  A topological
space admits a compatible uniformity iff it is completely regular.     A completely regular space $X$ has a maximum
uniformity $\U_M$ compatible with the topology. Any continuous function from a completely regular space $X$ to a uniform space is uniformly continuous
from $(X,\U_M)$.

A completely regular, Hausdorff space is called a \emph{Tychonoff space}.
\index{Tychonoff space}
A compact Hausdorff space $X$ has a unique uniformity consisting of all neighborhoods of the diagonal $1_X$.

In Section 2, we define the barrier functions $m^f_d$ and $\ell^f_d$ of a relation $f$ on a set $X$
with respect to a pseudo-metric $d$ and we describe their elementary properties.
We use a symmetric definition which allows a jump at the beginning as well
as the end of a sequence. In Section 6, we show that the alternative definitions yield equivalent results in cases which include when
$f$ is a continuous map.

In Section 3, we describe the properties of the Conley relation $\CC_d f = \{ (x,y) : m^f_d(x,y) = 0 \}$ and
the Aubrey-Mather relation $\A_d f = \{ (x,y) : \ell^f_d(x,y) = 0 \}$.  Following \cite{A93} we regard
$\CC_d$ and $\A_d$ as operators on the set of relations on $X$. We observe that each of these operators is idempotent.

In Section 4, we consider Lyapunov functions. With the pseudo-metric $d$ fixed, a Lyapunov function $L$ for a relation $f$ on $X$
is a continuous map $L : X \to \R$ such that $(x,y) \in f$ implies $L(x) \leq L(y)$, or, equivalently,
$f \ \subset \ \leq_L$ where $\leq_L \ = \ \{ (x,y) : L(x) \leq L(y) \}$.  Notice that we follow \cite{A93} in using Lyapunov functions
which increase, rather than decrease, on orbits. Following  \cite{P} and \cite{FP} we show that the barrier functions can be
used to define Lyapunov functions.  If $g$ is a relation on $X$ with $f \subset g$ and $z \in X$
then $x \mapsto m^g_d(x,z)$ is a Lyapunov function for  $\CC_d f$ and $x \mapsto \ell^g_d(x,z)$ is a
Lyapunov function for  $\A_d f$. Even when $f$ is a map, it is convenient to use associated relations like
$g = f \cup 1_X$ or $g = f \cup \{(y,y)\}$ for $y$ a point of $X$.

In Section 5, we turn to uniform spaces.  The Conley relation
$\CC_{\U} f$ is the intersection of $\{ \CC_d f : d \in \Gamma(\U) \}$ and $\A_{\U} f$ is
the intersection of $\{ \A_d f : d \in \Gamma(\U) \}$. Thus, $(x,y) \in \CC_{\U} f$ iff $m^f_d(x,y) = 0$ for all
$d \in \Gamma(\U)$ and similarly $(x,y) \in \A_{\U} f$ iff $\ell^f_d(x,y) = 0$ for all
$d \in \Gamma(\U)$. While the gage definition is convenient to use, we
show that each of these relations has an equivalent description which uses the uniformity
directly. Each of these is a closed, transitive relation which contains $f$.  We let $\G f$ denote the smallest
closed, transitive relation which contains $f$, so that $f \subset \G f \subset \A_{\U} f \subset \CC_{\U} f$.

If $L$ is a uniformly continuous Lyapunov function for $f$ then it is automatically a Lyapunov function
for $\A_{\U} f$. If $X$ is Hausdorff and we let $L$ vary over all uniformly continuous Lyapunov functions for $f$ then
$1_X \cup \A_{\U} f = \bigcap_L \ \leq_L$.  That is, if $(x,y) \not\in 1_X \cup \A_{\U} f$, then there exists
a uniformly continuous Lyapunov function $L$ such that $L(x) > L(y)$. If, in addition, $X$ is second countable, then
there exists a uniformly continuous Lyapunov function $L$ such that  $1_X \cup \A_{\U} f  =  \ \leq_L$. If $X$ is Hausdorff and
we let $L$ vary over all Lyapunov functions for $\CC_{\U} f$ then
$1_X \cup \CC_{\U} f = \bigcap_L \ \leq_L$. If, in addition, $X$ is second countable, then
there exists a Lyapunov function $L$ such that  $1_X \cup \CC_{\U} f =  \ \leq_L$. These results use the barrier function Lyapunov functions
developed in the preceding section.

For the Conley relation there are special results. A set $A$ is called
\emph{$\U$ inward}\index{i@$\U$ inward} for a relation $f$ on $(X,\U)$ if for some $U \in \U$
$(U \circ f)(A) \subset A$. A continuous function $L : X \to [0,1]$ is called
an \emph{elementary Lyapunov function}\index{elementary Lyapunov function} if $(x,y) \in f$ and $L(x) > 0$ imply $L(y) = 1$.
For a $\U$ uniformly continuous elementary Lyapunov function $L$ the sets $\{ x : L(x) > \ep \}$ for $\ep \geq 0$ are open $\U$ inward sets.
On the other hand, if $A$ is a $\U$ inward set, then there exists a $\U$ uniformly continuous elementary Lyapunov function $L$ such that
$L = 0$ on $X \setminus A$ and $L = 1$ on $f(A)$.
Each set $\CC_{\U} f(x)$ is an intersection of inward sets. If $A$ is an open $\U$ inward set then it is $\CC_{\U} f$ $^+$invariant and the
maximum $\CC_{\U} f$ invariant subset $A_{\infty}$ is called the associated \emph{attractor}\index{attractor}.

Additional results can be obtained when the relation $f$ satisfies various topological conditions.  In Section 6, we consider \emph{upper semicontinuous}
(= usc) and \emph{compactly upper semicontinuous relations} (= cusc) relations and related topological results. Regarded as a relation, a continuous
map is cusc.  If a Hausdorff space $X$ is locally compact and $\sigma-$compact, or locally compact and paracompact with $f$ cusc, then
$\G f = \A_{\U_M} f$. We exhibit a homeomorphism on a metric space for which the inclusion $\G f \subset \A_{\U_M} f$ is proper.

At the end of the section we consider compactifications and the special results which hold for a compact Hausdorff space.
 In the Hausdorff uniform space context, one proceeds by finding a totally bounded uniformity
$\T \subset \U$ which is compatible with the topology on $X$ and then take the uniform completion.

\begin{theo}\label{theointro} Let $f$ be a closed relation on a Hausdorff uniform
space $(X,\U)$ with $X$ second countable.
There exists  $\T \subset \U$  a totally bounded uniformity, with $(\bar X, \bar \T)$ the
completion of $(X, \T)$, such that the space $\bar X$
is a compact Hausdorff space with its unique uniformity $\bar \T$ metrizable. Let $\bar f$ be the closure of $f$ in $\bar X \times \bar X$.
The uniformity $\T$ can be chosen so that
\begin{equation}\label{eq53ucintro}
\begin{split}
\bar f \cap (X \times X) \ = \ f, \quad
1_X \cup \G \bar f \cap (X \times X) \ = \ 1_X \cup \A_{ \U} f, \\
 \CC \bar f \cap (X \times X) \ = \ \CC_{ \U} f. \hspace{3cm}
\end{split}
\end{equation}
If $f$ is cusc, e.g. a continuous map, then $\G \bar f \cap (X \times X) \ = \   \A_{ \U} f$.

If $f$ is a uniformly continuous map then, in addition, we can choose $\T$ so that $\bar f$ is a continuous map on $\bar X$.
If $f$ is a uniform isomorphism then, in addition, we can choose $\T$ so that
 $\bar f$ is a homeomorphism on $\bar X$.
\end{theo}

If $X$ is a compact Hausdorff space, then every closed, $\CC f$ $^+$invariant set $K$ is an intersection of inward sets.
If a closed set $K$ is $\CC f$ invariant then it is an intersection of attractors and $K$ is determined by $K \cap |\CC f|$
which we call its \emph{trace}. In fact, $K = \CC f(K \cap |\CC f|)$. $K$ is an
attractor iff it is closed and $\CC f$ invariant and, in addition, its trace is a clopen
subset of $|\CC f|$.

In Section 7, we consider totally recurrent and chain transitive relations.
Let $f$ be a relation on a uniform space $(X,\U)$ and let $d \in \Gamma(\U)$.
For $F =  \G f, \A_d f, \A_{\U} f$, $\CC_d f$ or $\CC_{\U} f$ we will say that
$f$ is \emph{totally $F$ recurrent} when $F$ is an equivalence relation. If $f$ is a uniformly continuous
map then $f$ is totally $F$ recurrent iff $1_X \subset F$, i.e.\ $F$ is reflexive.

If $\A_{\U} f$ is an equivalence relation then the quotient space $X/\A_{\U} f$ is \emph{completely Hausdorff}, i.e.\ the
continuous real-valued functions distinguish points. On the other hand, there exist examples such that the quotient is not
regular and so the topology is strictly finer than the weak topology generated by the continuous functions.
The latter is completely regular and the barrier functions $\ell^f_d$, when symmetrized, generate the gage of a compatible uniformity.

Similarly, if  $\CC_{\U} f$ is an equivalence relation then the quotient space $X/\CC_{\U} f$ is \emph{ totally disconnected}, i.e.\
the clopen sets distinguish points. Again there exist examples such that the quotient is not
regular and so the topology is strictly finer than the weak topology generated by the clopen subsets, i.e.\ it is not \emph{zero-dimensional}.
The barrier functions $m^f_d$, when symmetrized, are pseudo-ultrametrics generating the gage of a uniformity compatible with the latter
zero-dimensional topology.

The relation $f$ is called \emph{$\U$ chain transitive} when $\CC_{\U} f = X \times X$.
It is called \emph{$\U$ chain-mixing} if for every pair of points $x,y \in X$ and for every $d \in \Gamma(\U)$ and $\ep > 0$
there exists a positive integer $N$ such that  for every $n \geq N$ there are $\ep, d$ chains of length $n$
connecting $x$ and $y$.  A $\U$ chain-transitive relation $f$ is not $\U$ chain-mixing iff there exists a $\U$ uniformly continuous map taking $f$
to a non-trivial periodic cycle. It follows that $f$ is $\U$ chain-mixing iff the product relation $f \times f$ is $\U$ chain-transitive.
If $f$ is a $\U$ uniformly continuous map, then it is $\U$ chain-mixing iff
for every positive integer $n$ the iterate $f^n$ is $\U$ chain-transitive.

In Section 8 we restrict to compact metrizable spaces. The relation $\G f$ is the intersection
of the $\A_d f$'s as $d$ varies over $\hat \Gamma$, the set of metrics compatible with
the topology. If we take the union, which we denote $\W f$,
it is not obvious that the result is closed or transitive.  We prove it is both by giving a
uniformity characterization. The set $|\W f|$ is referred to as the
Ma\~{n}\'{e} set by Fathi and Pageault. Using the uniformity characterization we give an
alternative proof of their  description, for a homeomorphism $f$,
$|\W f| = |f| \cup |\CC (f|(X \setminus |f|^{\circ}))|$.

\newpage

\section{ Barrier Functions}
\vspace{.5cm}

Let $f$ be a  relation on a pseudo-metric space $(X,d)$. That is, $f$ is a  subset of
$X \times X$ and $d$ is a pseudo-metric on the non-empty set $X$.

Let $f^{\times n}$\index{f@$f^{\times n}$} be the $n-$fold product of copies of
$f$, i.e.\ the space of
sequences in $f$ of length $n \geq 1$, so that an element of
$f^{\times n}$ is a sequence $[a,b] = (a_1,b_1), (a_2,b_2),..,(a_n,b_n)$
 of pairs in
$f$. If $[a,b] \in f^{\times n}, [c,d] \in f^{\times m}$, then the
\emph{concatenation}\index{concatenation} $[a,b]\cdot[c,d] \in f^{\times n+m}$
is the sequence of pairs $(x_i,y_i) = (a_i,b_i)$ for $i = 1, \dots , n$ and $(x_i,y_i) = (c_{i-n},d_{i-n})$ for
$i = n+1, \dots, n+m$.

Define for $(x,y) \in X \times X$ and $[a,b] \in f^{\times n} $ the $xy$ \emph{chain-length of}\index{chain-length} $[a,b]$ (with respect to $d$)
to be the sum
\begin{equation}\label{eq1a}
d(x,a_1) + \Sigma_{i = 1}^{n-1} d(b_i,a_{i+1}) + d(b_n,y)
\end{equation}
 and the $xy$ \emph{chain-bound of}\index{chain-bound} $[a,b]$ (with respect to $d$)
to be
\begin{equation}\label{eq1b}
\max (d(x,a_1), d(b_1,a_2),\dots, d(b_{n-1},a_n), d(b_n,y)).
\end{equation}
 That is, for the vector
$(d(x,a_1), d(b_1,a_2),\dots, d(b_{n-1},a_n), d(b_n,y))$, the chain-length is the $L^1$ norm and the chain-bound
is the $L^{\infty}$ norm. We could proceed as below, using the $L^p$ norm for any $1 \leq p \leq \infty$.

For $(x,y) \in X \times X$, define
\begin{equation}\label{eq1}
\begin{split}
\ell_d^f(x,y) \  = \  \inf  \ \{ d(x,a_1) + \Sigma_{i = 1}^{n-1} d(b_i,a_{i+1}) + d(b_n,y) : \hspace{2cm}\\
[a,b] \in f^{\times n}, n = 1,2,... \}. \hspace{3cm} \\
m_d^f(x,y) \  = \  \inf \ \{ \max (d(x,a_1), d(b_1,a_2),\dots, d(b_{n-1},a_n), d(b_n,y)) : \\
[a,b] \in f^{\times n}, n = 1,2,... \}. \hspace{3cm} \\
\end{split}
\end{equation}

The functions
$\ell_d^f$\index{l@$\ell_d^f$} and $m_d^f$\index{m@$m_d^f$} are the
\emph{barrier functions}\index{barrier functions} for $f$. Clearly, $m_d^f \leq \ell_d^f$.

Using  $n = 1$, we see that for all $(a,b) \in f$
\begin{equation}\label{eq2}
\begin{split}
\ell_d^f(x,y) \  \leq \ d(x,a) + d(b,y), \\
m_d^f(x,y) \  \leq \  \max (d(x,a),d(b,y)).
\end{split}
\end{equation}
and so
\begin{equation}\label{eq2a}
(x,y) \in f \qquad \Longrightarrow \qquad m_d^f(x,y) \ = \ \ell_d^f(x,y) \ = \ 0.
\end{equation}
by using $(a,b) = (x,y)$.

For the special case of $f = \emptyset$ we define
\begin{equation}\label{eq3}
m_d^{\emptyset} \  = \ diam(X), \quad \text{and} \quad \ell_d^{\emptyset} \  = \  2 diam(X),
\end{equation}
the constant functions.

%
%

By using equation (\ref{eq2}) with $(a,b) = (y,y)$ and the triangle inequality in (\ref{eq1}) we see that
\begin{equation}\label{eq4}
 \ell_d^{1_X}(x,y)  \ = \ d(x,y).
\end{equation}

Define for the pseudo-metric $d$
\begin{equation}\label{eq36aa1}
Z_d \ = \ \{ (x,y) : d(x,y) = 0 \}.
\end{equation}
Thus, $Z_d$ is a closed equivalence relation which equals $1_X$ exactly when $d$ is a metric.
$Z_d$ is the closure in $X \times X$ of the diagonal $1_X$.

\begin{lem}\label{domlem} Let $f$ be a relation on $(X,d)$ with $A = Dom(f) = f^{-1}(X)$.  If $f \subset Z_d$, then
\begin{equation}\label{eq4aaa}
\ell_d^{f}(x,y) \  = \  \inf \{ d(x,a) + d(a,y) : a \in A \} \ \geq \ d(x,y)
\end{equation}
with equality if either $x$ or $y$ is an element of $A$.

If $d$ is a pseudo-ultrametric then
\begin{equation}\label{eq4aab}
m_d^{f}(x,y) \  = \  \inf \{\max(d(x,a),d(a,y)) : a \in A \} \ \geq \ d(x,y)
\end{equation}
with equality if either $x$ or $y$ is an element of $A$.
\end{lem}

{\bfseries Proof:} If $(a,b) \in f$ then $d(a,b) = 0$ and so the $xy$ chain-length of  $[(a,b)]$ is $d(x,a) + d(a,y)$.
If $[a,b] \in f^{\times n}$ then $d(a_i,b_i) = 0$ for all $i$ implies that with $a = a_1$ the $xy$ chain-length of
$[a,b]$ is at least $d(x,a) + d(a,y)$ by the triangle inequality.

If $d$ is a pseudo-ultrametric then the $xy$ chain-bound of $[(a,b)]$ is $\max(d(x,a),d(a,y)))$ and if  $[a,b] \in f^{\times n}$,
then with with $a = a_1$ the $xy$ chain-bound of
$[a,b]$ is at least $\max(d(x,a),d(a,y)))$ by the ultrametric version of the triangle inequality.

$\Box$ \vspace{.5cm}

In particular, if $A$ is a nonempty  subset of $X$, then
\begin{equation}\label{eq4a}
\ell_d^{1_A}(x,y) \  = \  \inf \{ d(x,a) + d(a,y) : a \in A \} \ \geq \ d(x,y)
\end{equation}
with equality if either $x$ or $y$ is an element of $A$.

It is clear that $f \subset g$ implies $f^{\times n} \subset g^{\times n}$ and so
\begin{equation}\label{eq5}
f \subset g \quad \Longrightarrow \quad \ell_d^g \ \leq \ \ell_d^f \quad
\text{and} \quad m_d^g \ \leq \ m_d^f \quad \mbox{on} \ X \times X .
\end{equation}
In particular, if $A$ is a  subset of $X$, then
\begin{equation}\label{eq6}
 \ell_d^f \ \leq \ \ell_d^{f|A} \quad \text{and} \quad m_d^f \ \leq \ m_d^{f|A}.
\end{equation}
The relation $f$ is reflexive when $1_{X} \subset f$. We see from (\ref{eq4})
\begin{equation}\label{eq7}
1_X \subset f \qquad \Longrightarrow \qquad \ell_d^f \ \leq \ d \quad \mbox{on} \ X \times X .
\end{equation}

If $[a,b] \in f^{\times n}$, then we let $[a,b]^{-1} \in (f^{-1})^{\times n}$ be $ (b_n,a_n), (b_{n-1},a_{n-1}),...,(b_1,a_1)$.
Using these reverse sequences we see immediately that
\begin{equation}\label{eq8}
 \ell_d^f(x,y)  \ = \ \ell_d^{f^{-1}}(y,x) \quad \text{and} \quad  m_d^f(x,y)  \ = \ m_d^{f^{-1}}(y,x)
\end{equation}
for all $x,y \in X$.

\begin{prop}\label{prop00a} Let $f$ be a  relation on $(X,d)$. Let $x,y,z,w \in X$.

\begin{itemize}
\item[(a)] The directed triangle inequalities hold:
\begin{equation}\label{eq9}
\begin{split}
\ell_d^f(x,y) \ \leq \ \ell_d^f(x,z) + \ell_d^f(z,y), \\
m_d^f(x,y) \ \leq \ m_d^f(x,z) + m_d^f(z,y).
\end{split}
\end{equation}

\item[(b)] Related to the ultrametric inequalties, we have:
\begin{equation}\label{eq9a}
m_d^f(x,y) \ \leq \ \max(m_d^f(x,z) + m_d^f(z,z), m_d^f(z,z) + m_d^f(z,y)).
\end{equation}

\item[(c)] From
\begin{equation}\label{eq10}
\begin{split}
\ell_d^f(x,y) \ \leq \ d(x,w) + \ell_d^f(w,z) + d(z,y) \quad \mbox{for all} \ w, x, y,z \in X, \\
m_d^f(x,y) \ \leq \ d(x,w) + m_d^f(w,z) + d(z,y) \quad \mbox{for all} \ w, x, y,z \in X
\end{split}
\end{equation}
we obtain that the functions $\ell_d^f$ and $m_d^f$ from $X \times X$ to $\R$ are Lipschitz with Lipschitz
constant $\leq 2$.
\end{itemize}
\end{prop}

{\bfseries Proof:} (a) For $x,y,z \in X$ and $[a,b] \in f^{\times n}, [c,d] \in f^{\times m}$, we note that
$d(b_n,c_{1}) \leq d(b_n, z) + d(z,c_{1})$. So the $xz$ chain-length of $[a,b]$ plus the $zy$ chain-length of $[c,d]$
is greater than or equal to the $xy$ chain-length of $[a,b]\cdot[c,d]$. Furthermore, the
$xz$ chain-bound of $[a,b]$ plus the $zy$ chain-bound of $[c,d]$
is greater than or equal to the $xy$ chain-bound of $[a,b]\cdot[c,d]$. The directed triangle inequalities (\ref{eq9})
follow.

(b) Let $[u,v] \in f^{\times p}$. We see that $d(b_n,u_{1}) \leq d(b_n, z) + d(z,u_{1})$ and
$d(v_p,c_{1}) \leq d(v_p, z) + d(z,c_{1})$. Hence, the larger of  the $xz$
chain-bound of $[a,b]$ plus the $zz$ chain-bound of $[u,v]$
and the $zz$ chain-bound of $[u,v]$ plus the $zy$ chain-bound of $[c,d]$ bounds the $xy$ chain-bound of
$[a,b] \cdot [u,v] \cdot [c,d]$. This implies (\ref{eq9a}).

(c) Similarly, $d(x,a_1) \leq d(x,w) + d(w,a_1)$ and $d(b_n,y) \leq d(b_n,z) + d(z,y)$
implies (\ref{eq10}) from which the Lipschitz results are clear.

$\Box$ \vspace{.5cm}

If $h$ is a map from $(X_1,d_1)$ to $(X_2,d_2)$ then $h$ is uniformly continuous if for every $\ep > 0$ there exists
$\d > 0$ such that $d_1(x,y) < \d$ implies $d_2(h(x),h(y)) < \ep$ for
all $x,y \in X_1$. We call $\d$ an $\ep$ modulus of uniform continuity.\index{modulus of uniform continuity}
The
map $h$ is Lipschitz with constant $K$ if
$d_2(h(x),h(y)) \leq K d_1(x,y)$ for all $x,y \in X_1$.

If $f_1$ is a relation on $X_1$ and $f_2$ is a relation on $X_2$ then we say that a function $h : X_1 \to X_2$ maps
$f_1$ to $f_2$ if $(h \times h)(f_1) \subset f_2$, i.e.\ $(x,y) \in f$ implies $(h(x),h(y))$. \index{h@$h$ maps $f_1$ to $f_2$}Since $h$ is a map,
$1_{X_1} \subset h^{-1} \circ h$ and $h \circ h^{-1} \subset 1_{X_2}$. From these it easily follows that
\begin{align}\label{eq10a}
(h \times h)(f_1) \ &= \ h \circ f_1 \circ h^{-1}, \\
(h \times h)(f_1) \ \subset \ f_2 \quad &\Longleftrightarrow \quad  h \circ f_1 \ \subset \ f_2 \circ h.
\end{align}
If $h$ maps $f_1$ to $f_2$ then clearly $h$ maps $f_1^{-1}$ and $f_2^{-1}$ and
\begin{equation}\label{eq10b}
h(|f_1|) \  \subset \  |f_2|.  \hspace{3cm}
\end{equation}

\begin{prop}\label{prop00b} Let $f_1$ and $f_2$ be relations on $(X_1,d_1)$
and $(X_2,d_2)$, respectively. Assume $h : X_1 \to X_2$
maps $f_1$ to $f_2$.

(a) If $h$ is uniformly continuous then for  $\ep > 0$ with $\d > 0$  an $\ep$ modulus of uniform continuity,
$m^{f_1}_{d_1}(x,y) < \d$ implies $m^{f_2}_{d_2}(h(x),h(y)) < \ep$ for all $x,y \in X_1$.

(b) If $h$ is Lipschitz with constant $K$ then $\ell^{f_1}_{d_1}(x,y) \leq K \ell^{f_2}_{d_2}(h(x),h(y))$
for all $x,y \in X_1$.
\end{prop}

{\bfseries Proof:} If $[a,b] \in f_1^{\times n}$ then $(h \times h)^{\times n}([a,b]) \in f_2^{\times n}$.  If $\d$ is
an $\ep$ modulus of uniform continuity then if the $xy$ chain-bound of $[a,b]$ is less then $\d$ then the
$h(x)h(y)$ chain-bound of $(h \times h)^{\times n}([a,b])$ is less than $\ep$.  If $h$ is Lipschitz with constant $K$ then the
$h(x)h(y)$ chain-length is at most $K$ times the $xy$ chain-length.

$\Box$

\vspace{1cm}

\section{ The Conley and Aubry-Mather Chain-Relations}
\vspace{.5cm}

For a  relation $f$ on $(X,d)$, the \emph{Conley chain relation}\index{relation!chain relation}\index{Conley chain relation} $\CC_d f$ is defined by
\begin{equation}\label{eq14a}
\CC_d f \ = \  \{ (x,y) : m^d_f(x,y) = 0 \},
\end{equation}
and the \emph{Aubry-Mather chain relation} \index{Aubry-Mather chain relation}is defined by
\begin{equation}\label{eq15}
\A_d f \ = \ \{ (x,y) : \ell^d_f(x,y) = 0 \}.
\end{equation}
\index{cc@$\CC_d f$}\index{aa@$\AA_d f$}

Because $m^f_d$ and $\ell^f_d$ are continuous,
it follows that $\CC_d f$ and $\A_d f$ are closed in $(X \times X, d \times d)$.
From the directed triangle inequalities (\ref{eq9}),
it follows that $\CC_d f$ and  $\A_d f$ are transitive, i.e.\
\begin{equation}\label{eq15a}
\begin{split}
\CC_d  f \circ \CC_d f \ \subset \  \CC_d  f, \hspace{3cm} \\
\A_d f \circ \A_d f \  \subset \  \A_d f. \hspace{3cm}
\end{split}
\end{equation}

From (\ref{eq2a})
we see that,
\begin{equation}\label{eq16}
f \   \subset \  \A_d f \  \subset \  \CC_d f. \hspace{2cm}
\end{equation}

If $A \subset X$ with $f \subset A \times A$ we can regard $f$ as a relation on $(X,d)$ or as a relation on
$(A,d|A \times A)$ where $d|(A \times A)$ is the restriction of the pseudo-metric $d$ to $A \times A$. It is clear that
if $f \subset A \times A$, then
\begin{equation}\label{eq16aa}
m_d^f| (A \times A) \ = \ m_{d|(A \times A)}^f \quad \text{and} \quad
\ell_d^f| (A \times A) \ = \ \ell_{d|(A \times A)}^f .
\end{equation}
and so
\begin{equation}\label{eq16ab}
 (\CC_d f) \cap (A \times A) \ = \ \CC_{d|(A \times A)} f \quad \text{and} \quad
 (\A_d f) \cap (A \times A) \ = \ \A_{d|(A \times A)} f.
\end{equation}
If $A$ is closed and $x,y \in A$ with either $x \not\in A$ or $y \not\in A$, then $\ell_d^f(x,y) \geq m_d^f(x,y) > 0$ and so
\begin{equation}\label{eq16abaa}
 (\CC_d f)  \ = \ \CC_{d|(A \times A)} f \quad \text{and} \quad
 (\A_d f) \ = \ \A_{d|(A \times A)} f.
\end{equation}

From (\ref{eq5}) we get monotonicity
\begin{equation}\label{eq17}
f \ \subset \ g \qquad \Longrightarrow \qquad \CC_d  f \ \subset \
\CC_d  g \quad \text{and} \quad \A_d f \ \subset \ \A_d g.
\end{equation}
and from (\ref{eq8})
\begin{equation}\label{eq18}
\CC_d (f^{-1}) \ = \ (\CC_d f)^{-1} \quad \text{and} \quad  \A_d (f^{-1}) \ = \ (\A_d f)^{-1},
\end{equation}
and so we can omit the parentheses.

\begin{prop}\label{prop1a} Let $f, g$ be  relations on $X$.
\begin{equation}\label{eq21}
m_d^{\CC_d f} \ = \ m_d^f \quad \text{and} \quad \ell_d^{\A_d f} \ = \ \ell_d^f
\end{equation}

The operators $\CC_d$ and $\A_d$ on  relations are idempotent\index{idempotent operator}.  That is,
\begin{equation}\label{eq20}
\CC_d (\CC_d f) \ = \ \CC_d f \quad \text{and} \quad \A_d (\A_d f) \ = \ \A_d f
\end{equation}

In addition,
\begin{equation}\label{eq20ab}
\CC_d (\CC_d f \cap \CC_d g) \ = \ \CC_d f \cap \CC_d g \  \text{and} \  \A_d (\A_d f \cap \A_d g) \ = \ \A_d f \cap \A_d g,
\end{equation}
\end{prop}

{\bfseries Proof:} Since $f \subset \A_d f \subset \CC_d f$ it follows from (\ref{eq5}) that
$m_d^{\CC_d f} \ \leq \  m_d^f$ and
$\ell_d^{\A_d f} \ \leq \  \ell_d^f$.

For the reverse inequality fix $x,y \in X$ an
let $ t > \ell_d^{\A_d f}(x,y)$  be arbitrary. Choose  $t_1$ with $t > t_1 > \ell_d^{\A_d f}(x,y)$.
Suppose that $[a,b] \in (\A_d f)^{\times n}$
whose $xy$ chain-length is less than $t_1$.
Let $\ep = (t - t_1)/2n $
For $i = 1,...,n$ we can choose an element of some $f^{\times n_i}$ whose $a_i b_i$ chain-length
is less than $\ep$. Concatenating these in order we obtain a sequence in $f^{\times m}$ with
$m = \Sigma_{i = 1}^{n} n_i$ whose $xy$ chain-length is at most $t_1 + 2n \ep \leq t$.
Hence, $\ell_d^f(x,y) \leq t$.
Letting $t$ approach $\ell_d^{\A_d f}(x,y)$ we obtain in the limit that
$\ell_d^f(x,y)  \leq \ell_d^{\A_d f}(x,y)$.

The argument to show $m_d^f(x,y)  \leq m_d^{\CC_d f}(x,y)$ is completely similar.

It is clear that (\ref{eq21}) implies (\ref{eq20}).

Finally, $\CC_d f \cap \CC_d g \subset \CC_d(\CC_d f \cap \CC_d g) \subset \CC_d(\CC_d f) = \CC_d f $ and similarly,
$\CC_d f \cap \CC_d g \subset \CC_d(\CC_d f \cap \CC_d g) \subset  \CC_d g $.  Intersect to get (\ref{eq20ab}) for $\CC_d$ and
the same argument yields the $\A_d$ result.

$\Box$ \vspace{.5cm}

\begin{cor}\label{cor00aa}For a relation $f$ on $(X,d)$ let $\bar f^d$ be the closure of $f$ in $(X \times X, d \times d)$.
\begin{align}\label{eq21a}
\begin{split}
m_d^{\bar f^d} \ = \ m_d^f \quad &\text{and} \quad \ell_d^{\bar f^d} \ = \ \ell_d^f,\\
\CC_d (\bar f^d) \ = \ \CC_d f \quad &\text{and} \quad \A_d (\bar f^d) \ = \ \A_d f
\end{split}
\end{align}
\end{cor}

{\bfseries Proof:} This is clear from (\ref{eq5}) and (\ref{eq21}) because
$f \subset \bar f^d \subset \A_d f \subset \CC_d f$.

$\Box$ \vspace{.5cm}

The \emph{Conley set}\index{Conley set} is the cyclic set $|\CC_d f| = \{ x : (x,x) \in \CC_d f \}$.
Since $|\CC_d f| $ is the pre-image of the closed set $\CC_d f \subset X \times X$ via the continuous map
  $x \mapsto (x,x)$ it follows that $|\CC_d f| \subset X$ is closed. The \emph{Aubry Set}\index{Aubry set} is the cyclic set
$|\A_d f| \subset X$ which is similarly closed.

From (\ref{eq16}) we clearly have
 $|\A_d f| \ \subset \ |\CC_d f|.$

On $|\CC_d f|$ the relation $\CC_d f \cap \CC_d f^{-1}$ is a closed equivalence
relation and on  $|\A_d f|$  $\A_d f \cap \A_d f^{-1}$ is a closed equivalence
relation.

Define the symmetrized functions
\begin{align}\label{eq24}
\begin{split}
sm_d^f(x,y) \ = \ &\max \{ m_d^f(x,y), m_d^f(y,x) \}, \hspace{2cm}\\
s\ell_d^f(x,y) \ = \ &\max \{ \ell_d^f(x,y), \ell_d^f(y,x) \}. \hspace{2cm}
\end{split}
\end{align}
\index{s@$s\ell_d^f, sm_d^f$}

\begin{prop}\label{prop1aa} Let $f$ be a  relation on $X$. Let $x,y,z \in X$

\begin{itemize}
\item[(a)] $sm_d^f(x,y) \leq s\ell_d^f(x,y)$

\item[(b)] The functions $sm_d^f$ and $s\ell_d^f$ are symmetric and  satisfy the triangle inequality.

\item[(c)] The functions $sm_d^f, s\ell_d^f : X \times X \to \R$ are Lipschitz with Lipschitz constant
less than or equal to $2$.

\item[(d)]
\begin{equation}\label{eq24a}
\begin{split}
sm_d^f(x,y) = 0 \qquad \Longleftrightarrow \qquad (x,y), (y,x) \in \CC_d f \quad \text{and so} \ x, y \in |\CC_d f|, \\
s\ell_d^f(x,y) = 0 \qquad \Longleftrightarrow \qquad (x,y), (y,x) \in \A_d f \quad \text{and so} \ x, y \in |\A_d f|.
\end{split}
\end{equation}

\item[(e)]
\begin{equation}\label{eq24b}
\begin{split}
y \in |\CC_d f| \qquad \Longrightarrow \qquad sm_d^f(x,y) \ \leq \ d(x,y), \\
y \in |\A_d f| \qquad \Longrightarrow \qquad s\ell_d^f(x,y) \ \leq \ d(x,y).
\end{split}
\end{equation}

\item[(f)] If $z \in |\CC_d f|$ then $m_d^f(x,y) \leq max(m_d^f(x,z),m_d^f(z,y)).$

\end{itemize}
\end{prop}

{\bfseries Proof:}  (a) is obvious as is symmetry in (b),  i.e.\
$sm_d^f(x,y) = sm_d^f(y,x)$ and $s\ell_d^f(x,y) = s\ell_d^f(y,x)$.
The triangle inequality for $s \ell_d^f$ follows from
\begin{equation}\label{eq25}
\begin{split}
s\ell_d^f(x,z) + s\ell_d^f(z,y) \geq \ell_d^f(x,z) + \ell_d^f(z,y) \geq \ell_d^f(x,y), \\
s\ell_d^f(x,z) + s\ell_d^f(z,y) \geq \ell_d^f(z,x) + \ell_d^f(y,z) \geq \ell_d^f(y,x),
\end{split}
\end{equation}
with a similar argument for for $s m_d^f$.
imply that $s\ell_d^f$ satisfies the triangle inequality.

By Proposition \ref{prop00a}(c) $m_d^f$ and $\ell_d^f$ are Lipschitz.  Then (c)
 follows from Lemma \ref{introalglem}.

The equivalences in (d) are obvious. By transitivity, $(x,y), (y,x) \in \CC_d f$ implies $(x,x), (y,y) \in \CC_d f$.
Similarly, for $\A_d f$.

(e) If $s m_d^f(y,y) = 0$ then $ s m_d^f(x,y) = s m_d^f(x,y) - s m_d^f(y,y) \leq d(x,y)$ by
(c). Similarly, for $s \ell_d^f$.

(f) follows from Proposition \ref{prop00a}(b).

$\Box$ \vspace{.5cm}

We immediately obtain the following.

\begin{cor}\label{cor1ab}  The map $s\ell_d^f$ restricts to define a pseudo-metric on
$|\A_d f|$ and induces a metric on the quotient space of
$\A_d f \cap \A_d f^{-1}$ equivalence classes.  Furthermore, the projection map from $|\A_d f|$ to the space of
equivalence classes has Lipschitz constant at most 2 with respect to this metric.

The map $sm_d^f$ restricts to define a pseudo-ultrametric on
$|\CC_d f|$ and induces an ultrametric on the quotient space of
$\CC_d f \cap \CC_d f^{-1}$ equivalence classes.  Furthermore, the projection map from $|\CC_d f|$ to the space of
equivalence classes has Lipschitz constant at most 2 with respect to this metric.
\end{cor}

$\Box$ \vspace{.5cm}

Let $f_1$ and $f_2$ be relations on $X_1$ and $X_2$, respectively. Recall that $h : X_1 \to X_2$
maps $f_1$ to $f_2$ when $h \circ f_1 \circ h^{-1} = (h \times h)(f_1) \subset f_2$, i.e.\ if $(x,y) \in f_1$ implies
$(h(x),h(y)) \in f_2$.  It then follows that $h$ maps $f_1^{-1}$ to $f_2^{-1}$.

\begin{prop}\label{prop00c} Let $f_1$ and $f_2$ be relations on $(X_1,d_1)$ and
$(X_2,d_2)$, respectively. Assume $h : X_1 \to X_2$
maps $f_1$ to $f_2$.

(a) If $h$ is uniformly continuous, then $h$ maps $\CC_{d_1} f_1$ to
$\CC_{d_2} f_2$ and $\CC_{d_1} f_1 \cap \CC_{d_1} f_1^{-1}$ to
$\CC_{d_2} f_2 \cap \CC_{d_2} f_2^{-1}$. So $h$ maps each $\CC_{d_1} f_1 \cap \CC_{d_1} f_1^{-1}$ equivalence class in
$|\CC_{d_1} f_1|$ into a $\CC_{d_2} f_2 \cap \CC_{d_2} f_2^{-1}$ equivalence class in
$|\CC_{d_2} f_2|$.

(b) If $h$ is Lipschitz,  then $h$ maps $\A_{d_1} f_1$ to $\A_{d_2} f_2$ and $\A_{d_1} f_1 \cap \A_{d_1} f_1^{-1}$ to
$\A_{d_2} f_2 \cap \A_{d_2} f_2^{-1}$. So $h$ maps each $\A_{d_1} f_1 \cap \A_{d_1} f_1^{-1}$ equivalence class in
$|\A_{d_1} f_1|$ into a $\A_{d_2} f_2 \cap \A_{d_2} f_2^{-1}$ equivalence class in
$|\A_{d_2} f_2|$.
\end{prop}

{\bfseries Proof:} This obviously follows from Proposition \ref{prop00b}.

$\Box$ \vspace{.5cm}

We conclude this section with some useful computations.

Recall that
\begin{equation}\label{eq36aa}
Z_d \ = \ \{ (x,y) : d(x,y) = 0 \}.
\end{equation}
\index{z@$Z_d$}

\begin{prop}\label{prop3} Let $f$ be a relation on $X$ and $A$ be a nonempty, closed
subset of $X$

(a) For $x,y \in X$
\begin{align}\label{eq37}
\begin{split}
\ell_d^{1_A \cup f}(x,y) \  &= \  \min(\ell_d^{f}(x,y),\ell_d^{1_A}(x,y)), \\
\ell_d^{1_X \cup f}(x,y) \  &= \  \min(\ell_d^{f}(x,y),d(x,y)),\\
s\ell_d^{1_X \cup f}(x,y) \  &= \  \min(s\ell_d^{f}(x,y),d(x,y)), \\
s\ell_d^{1_X \cup f}(x,y) \  &= \  s\ell_d^{f}(x,y) \qquad \mbox{if} \ x \in |\A_d f|.
\end{split}
\end{align}

(b) \begin{align}\label{eq38}
\begin{split}
\A_d (1_A \cup f) \  &= \  Z_d \cap (A \times A) \cup \A_d f, \\
\A_d (1_X \cup f) \  &= \  Z_d \cup \A_d f.
\end{split}
\end{align}
$s\ell_d^{1_X \cup f}$ is a pseudo-metric on $X$ whose associated metric space
is the quotient space of $X$ by the equivalence relation $Z_d \cup (\A_d f \cap \A_d f^{-1})$.
The quotient map has Lipschitz constant at most 2.


\end{prop}

{\bfseries Proof:} (a) By (\ref{eq5})  $\ell_d^{1_A \cup f} \leq \min(\ell_d^{f},\ell_d^{1_A})$.
 By (\ref{eq4}) $\ell_d^{1_X} = d$.

Let $[a,b] \in (1_A \cup f)^n$. If $(a_i,b_i) \in 1_A$ for all $i$ then omit all but one of the
pairs to obtain an element of $1_A^{\times 1}$. Otherwise, omit the pairs
$(a_i,b_i) \in 1_A$ and renumber.  We then obtain
a sequence in $f^{\times m}$ for some $m$ with $1 \leq m \leq n$. Furthermore, in either
case the $xy$ chain-length has not increased.
For example, if $(a_i,b_i) \in 1_A$ for some $1 < i < n$ then since $a_i = b_i$ the triangle inequality implies
$d(b_{i-1},a_{i+1}) \leq d(b_{i-1},a_i) + d(b_i,a_{i+1})$. It follows that
$\ell_d^{1_X \cup f} \geq min(\ell_d^{f},\ell_d^{1_A})$.

 $s\ell_d^{1_X \cup f}(x,y) = \max[ \min(\ell_d^f(x,y),d(x,y)), \min(\ell_d^f(y,x),d(x,y))]$. This is $d(x,y)$
 except when $d(x,y) > \ell_d^f(x,y)$ and $d(x,y) > \ell_d^f(y,x)$, i.e.\ $d(x,y) > s\ell_d^f(x,y)$ in which case
 it is $s\ell_d^f(x,y)$.

(b) If $x \in|\A_d f|$ then by (\ref{eq24b}) $ \min(s\ell_d^{f}(x,y),d(x,y)) = s\ell_d^{f}(x,y)$.

 It follows $x \in |\A_d (1_A \cup f)|$ iff $\ell_d^f(x,y) = 0$ or $\ell_d^{1_A}(x,y) = 0$. By (\ref{eq4a})
 the latter is true iff $x, y \in A$ with $d(x,y) = 0$ since $A$ is closed.  Thus, (\ref{eq38})
 holds and the rest is obvious.

$\Box$ \vspace{.5cm}

 If $A,B$ are  subsets of $X$ then we can regard $A \times B$ as a  relation on $X$.  For any relation
 $g$ on $X$ we clearly have:
 \begin{equation}\label{eq38b}
 (A \times B) \circ g \circ (A \times B) \  \subset \  A \times B.
 \end{equation}

\begin{lem}\label{lem3a} If $A$ and $B$ are nonempty  subsets of $(X,d)$
and $x,y \in X$,then
\begin{align}\label{eq38d}
\begin{split}
m_d^{A \times B}(x,y) \ &= \ \max(d(x,A),d(y,B)) \quad \text{and} \\
 \ell_d^{A \times B}(x,y) \  &= \  d(x,A) + d(y,B)
\end{split}
\end{align}
where
$ d(x,A)  =  \inf\{ d(x,z) : z \in A \}.$
\end{lem}

{\bfseries Proof:}
If $[a,b] \in (A \times B)^n$ then $(a_1,b_n) \in A \times B$ with $xy$ chain-length $d(x,a_1) + d(y,b_n)$ no
larger than the $xy$ chain-length for $[a,b]$ and with $xy$ chain-bound $\max(d(x,a_1), d(y,b_n))$ no
larger than the $xy$ chain-bound for $[a,b]$. This proves (\ref{eq38d}).

$\Box$ \vspace{.5cm}

From Proposition \ref{prop3}
 we immediately get

\begin{cor}\label{cor3b} If $A$ and $B$ are nonempty  subsets of $X$ and $x,y \in X$ then
\begin{equation}\label{eq38f}
\ell_d^{1_X \cup (A \times B)}(x,y) \  = \  \min[d(x,y), d(x,A) + d(y,B)]. \hspace{2cm}
 \end{equation}
 \end{cor}

 $\Box$ \vspace{.5cm}

 {\bfseries Remark:} If $A = B$ then $s\ell_d^{1_X \cup (A \times A)} = \ell_d^{1_X \cup (A \times A)}$ is
 the pseudo-metric on $X$ induced by the equivalence relation $1_X \cup (A \times A)$ corresponding to smashing
 $A$ to a point.
 \vspace{.5cm}

\begin{lem}\label{lem3b} For $x,y,z \in X$
\begin{equation}\label{eq38g}
\begin{split}
m_d^{f \cup \{ (z,z) \}}(x,y) \ = \hspace{4cm} \\
 \min [\ m_d^f(x,y), \ \max [ \min (m_d^f(x,z), d(x,z)),\  \min ( m_d^f(z,y), d(z,y))  \ ].
\end{split}
\end{equation}
In particular, with $z = y$ or $z = x$
\begin{equation}\label{eq38h}
m_d^{f \cup \{ (y,y) \}}(x,y) \ = \ m_d^{f \cup \{ (x,x) \}}(x,y) \ = \\min [ m_d^f(x,y), d(x,y) ]. \hspace{2cm}
\end{equation}

If $(z,z) \in \CC_d f$, i.e.\ $z \in |\CC_d f|$, then $m_d^{f \cup \{ (z,z) \}} = m_d^f$.
\end{lem}

{\bfseries Proof:} Since $f \subset f \cup \{ (z,z) \}$ we have
$m_d^{f \cup \{ (z,z) \}} \leq \min (m_d^{f}, m_d^{\{ (z,z) \}})$.

Let $[a,b] \in (f \cup \{ (z,z) \})^{\times n}$. If $(z,z)$ occurs more than once in $[a,b]$ we can eliminate the
repeat and all of the terms between them without increasing the $xy$
chain-bound. Thus, we may take the infimum over those
$[a,b]$ in which $(z,z)$ occurs at most once.

The infimum of the $xy$ chain-bounds in $f^{\times n}$ is $m_d^f(x,y)$.
\begin{itemize}
\item The $xy$ chain-bound of
$(z,z)  \in (f \cup \{ (z,z) \})^{\times 1}$ is $\max( d(x,z), d(z,y))$.

\item If $[a,b]$ varies in $ (f \cup \{ (z,z) \})^{\times n}$ with $n > 1$ and $(a_i,b_i) = (z,z)$  only for $i = 1$,
then the infimum of the $xy$ chain-bounds is $\max( d(x,z), m_d^f(z,y))$.

\item If $[a,b]$ varies in $ (f \cup \{ (z,z) \})^{\times n}$ with $n > 1$ and $(a_i,b_i) = (z,z)$  only for $i = n$,
then the infimum of the $xy$ chain-bounds is $\max( m_d^f(x,z), d(z,y))$.

\item If $[a,b]$ varies in $ (f \cup \{ (z,z) \})^{\times n}$ with $n > 2$ and $(a_i,b_i) = (z,z)$  only for some $i$
with $1 < i < n$, then the infimum of the $xy$ chain-bounds is $\max( m_d^f(x,z), m_d^f(z,y))$.

\end{itemize}

Equation (\ref{eq38g}) then follows from Lemma \ref{introalglem}.

If $(z,z) \in \CC_d f$ then $f \subset f \cup \{ (z,z) \} \subset \CC_d f$. So
$m_d^{\CC_d f} \leq m_d^{ f \cup \{ (z,z) \} } \leq m_d^f$ by (\ref{eq5}) and so they are equal by (\ref{eq21}).

 $\Box$

\newpage

 \section{ Lyapunov Functions}
\vspace{.5cm}

 A Lyapunov function\index{Lyapunov function} for a  relation $f$ on a pseudo-metric space $(X,d)$ is a
continuous map $L : X \to \R$ such that
 \begin{equation}\label{eq39}
 (x,y) \in f \qquad \Longrightarrow \qquad L(x) \leq L(y).
 \end{equation}
 We follow \cite{A93} in using functions increasing on orbits rather than decreasing.

 The set of Lyapunov functions contains the constants and is closed under addition, multiplication by positive scalars,
 max, min and post composition with any continuous non-decreasing function on $\R$.
 A continuous function which is a pointwise limit
 of Lyapunov functions is itself a Lyapunov function.

 We define for a  real-valued function $L$ the relation
 \begin{equation}\label{eq40}
 \leq_L \  = \  \{ (x,y) : L(x) \leq L(y) \}.
 \end{equation}
 \index{le@$\leq_L$}
 This is clearly reflexive and transitive.
By continuity of $L$ the relation $\leq_L$ is closed and so contains $Z_d$.

The Lyapunov function condition (\ref{eq39}) can be restated as:
\begin{equation}\label{eq41}
f \ \subset \  \leq_L. \hspace{3cm}
\end{equation}

For a Lyapunov function $L$ and $x \in X$ we have
\begin{equation}\label{eq43}
L(z) \leq L(x) \leq L(w) \qquad \mbox{for \ } z \in f^{-1}(x), w \in f(x)
\end{equation}
 The point $x$ is called an $f$-\emph{regular point}\index{regular point} for $L$
when the inequalities are strict for all $z \in f^{-1}(x), w \in f(x)$.  Otherwise $x$ is called
 an $f$-\emph{critical point}\index{critical point} for $L$. Notice, for example, that if
$f^{-1}(x) = f(x) = \emptyset$ then these conditions hold vacuously and so $x$ is an $f$-regular point.

We denote by $|L|_f$ the set of $f$-critical points for $L$. Clearly,
\begin{equation}\label{eq43a}
|L|_f \ = \ \pi_1(A) \cup \pi_2(A) \qquad \text{where} \quad A \ = \ f \cap (L \times L)^{-1}(1_{\R}),
\end{equation}
and $\pi_1, \pi_2 : X \times X \to X$ are the two coordinate projections.

\begin{df}\label{def6}  Let $F$ be a  transitive relation on $(X,d)$ and let $\L$ be a collection
of Lyapunov functions for $F$.
 We define three conditions on $\L$.
\begin{itemize}
\item[ALG] If $L_1,L_2 \in \L$ and $c \geq 0$ then \\ $L_1 + L_2, \max(L_1,L_2), \min(L_1,L_2), c L_1, c, -c \in \L$.\index{condition!ALG}

\item[CON] For every sequence $\{ L_k \}$ of elements of $\L$ there exists a summable sequence of positive real numbers
$\{ a_k \}$ such that $\Sigma_k a_k L_k$ converges uniformly to an element of $\L$.\index{condition!CON}

\item[POIN] If $(x,y) \not\in Z_d \cup F$ then there exists $L \in \L$ such that $L(y) < L(x)$, i.e.\
$Z_d \cup F \ = \ \bigcap_{L \in \L} \ \leq_L$.\index{condition!POIN}

\end{itemize}
\end{df}
\vspace{.5cm}

\begin{theo}\label{theo8} Assume $(X ,d)$ is separable.
Let $F$ be a closed, transitive relation and $\L$ be a collection of
Lyapunov functions for $F$ which satisfies ALG, CON and POIN. There exists a sequence $\{ L_k \}$ in $\L$
such that
\begin{equation}\label{eq44}
\bigcap_k \ \leq_{L_k} \quad = \quad Z_d \cup F.
\end{equation}
If $\{ a_k \}$ is a positive, summable sequence such that
$L = \Sigma_n \ a_k L_k  \in \L$ then $L$ is a Lyapunov function for $F$ such that
$Z_d \cup F \ = \ \leq_L$ and
\begin{equation}\label{eq45}
x \in F(y) \qquad \Longrightarrow \qquad L(y) < L(x) \quad \mbox{unless} \ y \in F(x)
\end{equation}
In particular,
\begin{equation}\label{eq46}
|L|_F \ = \  |F|
\end{equation}
\end{theo}

{\bfseries Proof:} For each $(x,y) \in (X \times X) \setminus (Z_d \cup F)$ use POIN to choose
$L_{xy} \in \L$ such that $L_{xy}(y) < L_{xy}(x)$ and then neighborhoods $V_{xy} $ of $y$ and $ U_{xy}$ of
$x$ such that $\sup L_{xy}|V_{xy} < \inf L_{xy}|U_{xy}$ and so $\leq_{L_{xy}}$ is disjoint from
$U_{xy} \times V_{xy}$.  Because $(X,d)$ is separable, it is second countable and so
$(X \times X) \setminus (Z_d \cup F)$ is Lindel\"{o}f.
Choose a sequence of pairs $(x_k,y_k)$ so that $\{ U_{x_k y_k} \times V_{x_ky_k} \}$
covers $(X \times X) \setminus (Z_d \cup F)$  and let $L_k = L_{x_k y_k}$. Since $Z_d \cup F \subset \leq_L$ for
any Lyapunov function $L$, (\ref{eq44}) holds.

Now with $L = \Sigma_k \ a_k L_k$, (\ref{eq44}) implies $Z_d \cup F \ = \ \leq_L$.
If $x \in F(y)$ and $d(y,x) = d(x,y) = 0$
then $(y,x) \in F$ implies
$(x,x),(y,y), (x,y) \in F$, because $F$ is closed. Hence, $y \in F(x)$.
Assume $(x,y) \not\in Z_d$. Since $x \in F(y) $, $L_k(y) \leq L_k(x) $ for all $k$. If equality holds for all $k$ then
$(x,y) \in \bigcap_k \leq_{L_k} = Z_d \cup F$.  Since $(x,y) \not\in Z_d$ we
have $y \in F(x)$. If, instead, the inequality is
strict for some $k$ then since $a_k > 0, \ L(y) < L(x)$, proving (\ref{eq45}).

If $x \not\in |F|$ then for $z \in F^{-1}(x)$ and $w \in F(x)$ we have $x \in F(z) $
but not $z \in F(x)$  else by transitivity
$x \in |F|$. Hence, $L(z) < L(x)$.  Similarly, $L(x) < L(w)$. Thus, $x \not\in |L|_F$.

$\Box$ \vspace{.5cm}

\begin{df}\label{df9} For a  relation $f$ on $(X,d)$  and $K > 0$, a function $L : X \to \R$ is called
\emph{$K \ell_d^f$ dominated}\index{kl@$K \ell_d^f$ dominated} if for all $x,y \in X$
\begin{equation}\label{eq47}
L(x) - L(y) \quad \leq \quad K \ell_d^f(x,y),
\end{equation}
\emph{$K m_d^f$ dominated}\index{km@$K m_d^f$ dominated}  if for all $x,y \in X$
\begin{equation}\label{eq47aa}
L(x) - L(y) \quad \leq \quad K m_d^f(x,y).
\end{equation}
\end{df}
\vspace{.5cm}

\begin{theo}\label{theo10} Let $f$ be a  relation on $(X,d)$.

(a) If $L$ is a $K \ell_d^f$ dominated function then it is a Lyapunov function for $\A_d f$ and
so is a  Lyapunov function for $f$.
If $L$ is a $K m_d^f$ dominated function then it is a $K \ell_d^f$ dominated
function and is a Lyapunov function for $\CC_d f$ .

(b) If $L$ is a Lyapunov function for $f$ which is Lipschitz with respect to $d$ with
Lipschitz constant at most $K$ then it is a
$K \ell_d^f$ dominated function and so is a $\A_d f$ Lyapunov function.
\end{theo}

{\bfseries Proof:} (a) If $(x,y) \in \A_d f$ then $\ell_d^f(x,y) = 0$ and so for a $K \ell_d^f$ dominated function
$L(x) - L(y) \leq 0$. Similarly, if $(x,y) \in \CC_d f$ and $L$ is $K m_d^f$ dominated, then $L(x) - L(y) \leq 0$.
Since $m_d^f \leq \ell_d^f$ a  $K m_d^f$ dominated function is a $K \ell_d^f$ dominated function.

(b) Assume $L$ is an $f$ Lyapunov function with Lipschitz constant $K$ and $x,y \in X$. For any $[a,b]\in f^{\times n}$
we note that each $L(a_i) - L(b_i) \leq 0$ since $(a_i,b_i) \in f$ and $L$ is a Lyapunov function for $f$.  Hence,
\begin{equation}\label{eq48}
\begin{split}
L(x) - L(y) \quad = \quad L(x) - L(a_1) + L(a_1) - L(b_1) + L(b_1) - L(a_2) + \\
... + L(a_n) - L(b_n) + L(b_n) - L(y) \quad \leq  \hspace{1.5cm} \\ L(x) - L(a_1) + \Sigma_{i=1}^{n-1}
L(b_i) - L(a_{i+1}) + L(b_n) - L(y) \quad \leq K \ell.
\end{split}
\end{equation}
where $\ell$ is the $xy$ chain-length of $[a,b]$. Taking the infimum over the sequences $[a,b]$  we obtain
(\ref{eq47}).   Hence, $L$ is a $\A_d f$ Lyapunov function by part (a).

$\Box$ \vspace{.5cm}

\begin{prop}\label{prop11}  Let $f \subset g$ be  relations on $(X,d)$.
For any $z \in X$, the function
defined by $x \mapsto \ell_d^g(x,z)$ is a bounded, $1 \ell_d^f$ dominated function, and
the function
defined by $x \mapsto m_d^g(x,z)$ is a bounded, $1 m_d^f$ dominated function. \end{prop}

{\bfseries Proof:} By the directed triangle inequalities for $\ell_d^g$ and $m_d^g$ we have
\begin{equation}\label{eq49}
\ell_d^g(x,z) - \ell_d^g(y,z) \ \leq \ \ell_d^g(x,y) \quad \text{and} \quad m_d^g(x,z) - m_d^g(y,z) \ \leq \ m_d^g(x,y)
\end{equation}
Since $f \subset g$, $\ell_d^g(x,y) \leq \ell_d^f(x,y)$ and $m_d^g(x,y) \leq m_d^f(x,y)$ by (\ref{eq5}).

$\Box$ \vspace{.5cm}

\begin{theo} \label{theo13} For $f$  a  relation on $(X,d)$ let $\L_{\ell}$ be the set of bounded, continuous functions which are
$K \ell_d^f$ dominated for some positive $K$. Each $L \in \L_{\ell}$ is a $\A_d f$ Lyapunov function and so satisfies
\begin{equation}\label{eq50}
\A_d f \ \subset \ \leq_L \qquad \mbox{and} \qquad |\A_d f| \ \subset \ = |L|_{\A_d f}.
\end{equation}
The collection $\L_{\ell}$ satisfies the conditions ALG, CON,  and POIN  with respect to $F = \A_d f$.
\end{theo}

{\bfseries Proof:} Each $L$ in $\L_{\ell}$ is a $\A_d f$ Lyapunov function by Theorem \ref{theo10} and so
the first inclusion of (\ref{eq50}) follows by definition.  Clearly, if $(x,x) \in \A_d f$ then
$x$ is a $\A_d f$ critical point.

For $\L_{\ell}$ ALG is easy to check, see, e.g. Lemma \ref{introalglem}. For CON let $\{ L_k \}$
be a sequence in $\L_{\ell}$ and choose for each
$k$, $M_k \geq 1$ which bounds $|L_k(x)|$ for all $x \in X$ and so that $L_k$ is $M_k \ell^d_f$ dominated.  If
$\{ b_k \}$ is any positive, summable sequence with $\sum b_k = 1$, then
$a_k = b_k/M_k > 0$ is summable and $\Sigma_k \ a_k L_k$
converges uniformly to a function which is $1 \ell^d_f$ dominated. Thus, CON holds as well.

Now assume $(x,y) \not\in Z_d \cup \A_d f$. Let $g = 1_X \cup f$. By Proposition \ref{prop11}
$L(w) = \ell_d^g(w,y)$ defines a $1 \ell_d^f$ dominated function which is a $\A_d f$ Lyapunov function
by Theorem \ref{theo10}(a).

 By Proposition \ref{prop3}
$L(w) = min(\ell_d^f(w,y), d(w,y))$. Hence, $L(y) = 0$. Since $(x,y) \not\in Z_d \cup \A_d f$, $L(x) > 0$.
This proves POIN.

$\Box$ \vspace{.5cm}

\begin{theo} \label{theo15a} For $f$  a  relation on $(X,d)$ let $\L_{m}$ be the set of bounded, continuous functions which are
$K m_d^f$ dominated for some positive $K$. Each $L \in \L_{m}$ is a $\CC_d f$ Lyapunov function and so satisfies
\begin{equation}\label{eq57a}
\CC_d f \ \subset \ \leq_L \qquad \mbox{and} \qquad |\CC_d f| \ \subset \  |L|_{\CC_d f}.
\end{equation}
The collection $\L_{m}$ satisfies the conditions ALG, CON,  POIN with respect to $F = \CC_d f$.
\end{theo}

{\bfseries Proof:} Each $L$ in $\L_{m}$ is a $\CC_d f$ Lyapunov function by Theorem \ref{theo10} and so
the first inclusion of (\ref{eq57a}) follows by definition.  Clearly, if $(x,x) \in \CC_d f$ then
$x$ is a $\CC_d f$ critical point.

For $\L_{m}$ ALG again follows from Lemma \ref{introalglem}. For CON let $\{ L_k \}$
be a sequence in $\L_{m}$ and choose for each
$k$, $M_k \geq 1$ which bounds $|L_k(x)|$ for all $x \in X$ and such that $L_k$ is $M_k m_d^f$ dominated.  If
$\{ b_k \}$ is any positive, summable sequence with $\sum b_k = 1$, then
$a_k = b_k/M_k > 0$ is summable and $\Sigma_k \ a_k L_k$
converges uniformly to a function which is  $1 m_d^f$. Thus, CON holds as well.

Now assume $(x,y) \not\in Z_d \cup \CC_d f$. Let $g = f \cup \{(y,y) \}$. By Proposition \ref{prop11}
$L(w) = m_d^g(w,y)$ defines a $1 m_d^f$ dominated function.
  By Equation (\ref{eq38h})
$L(w) = min(\ell_d^f(w,y), d(w,y))$. Hence, $L(y) = 0$. Since $(x,y) \not\in Z_d \cup \CC_d f$, $L(x) > 0$.
This proves POIN.

$\Box$

\vspace{1cm}

\section{ Conley  and Aubry-Mather Relations for Uniform Spaces}
\vspace{.5cm}

Let $\U$ be a  uniformity on $X$ with gage $\Gamma$, the set of all bounded
pseudo-metrics $d$ on $X$ such that the uniformity
$\U(d) $ is contained in $\U$.

For a relation $f$ on $X$ we define the Conley relation and
Aubry-Mather relation associated with the uniformity.

\begin{equation}\label{eq58a}
\CC_{\U} f \ = \ \bigcap_{d \in \Gamma}  \CC_d f, \qquad \text{and} \qquad \A_{\U} f \ = \ \bigcap_{d \in \Gamma} \A_d f
\end{equation}
with $|\CC_{\U} f|$ the Conley set and $|\A_{\U} f|$ the Aubry set.\index{cc@$\CC_{\U} f$}\index{aa@$\A_{\U} f$}

Thus, $\CC_{\U} f$ and $\A_{\U} f$ are closed, transitive relations on $X$ which contain $f$. We define $\G f$ to be the intersection of all the closed, transitive relations which contain $f$.\index{g@$\G f$}  Thus, $\G f$ is the smallest closed, transitive
relation which contains $f$.  Clearly,
\begin{equation}\label{eq58ab}
f \subset \ \G f \ \subset \ \A_{\U} f \ \subset  \CC_{\U} f.
\end{equation}

Thus, $(x,y) \in \CC_{\U} f$ if for every $d \in \G$ and every $\ep > 0$ there
exists $[a,b] \in f^{\times n}$ with $n \geq 1$
such that the $xy$ chain-bound of $[a,b]$ with respect to $d$ is less than $\ep$.

If  $[a,b] \in f^{\times n}$ with $n \geq 1$ and $U \in \U$ we say that $[a,b]$ is an \emph{ $xy, U$ chain} for $f$ if
$(x,a_1),(b_1,a_2), \dots (b_{n-1},a_n), (b_n,y) \in U$. Clearly, then $[a,b]^{-1}$ is a $yx, U^{-1}$ chain for $f^{-1}$.

Since the $V_d^{\ep}$'s for $d \in \Gamma(\U)$ and $\ep > 0$ generate the uniformity,
it is clear that the pair $(x,y) \in \CC_{\U} f$ iff for every $U \in \U$ there
exists an $xy, U$ chain\index{cxa@$xy,U$ chain} for $f$. This provides a uniformity description
of $\CC_{\U} f$.

Similarly, $(x,y) \in \A_{\U} f$ if for every $d \in \G$ and every $\ep > 0$ there
exists $[a,b] \in f^{\times n}$ with $n \geq 1$
such that the $xy$ chain-length of $[a,b]$ with respect to $d$ is less than $\ep$.

Following \cite{W16} we obtain a uniformity description of  $\A_{\U} f$.

If $\xi = \{ U_k : k \in \N \}$ is a sequence of elements of $\U$
 and $(x,y) \in X \times X, $ we call $[a,b] \in f^{\times n}$
an \emph{$ \xi$ sequence chain from $x$ to $y$}\index{cxa@$\xi$ sequence chain} if there is an injective map $\s : \{ 0, \dots, n \} \to \N $ such that
$(b_i,a_{i+1}) \in U_{\s(i)}$ for $i = 0, \dots, n$ with $b_0 = x, a_{n+1} = y$.


\begin{theo}\label{wisetheo2}  For a relation $f$ on a uniform  space $(X,\U)$,
$ (x,y)  \in \A_{\U} f$ iff  for every sequence $\xi$ in $\U$ there is a $\xi$ sequence chain from $x$ to $y$. \end{theo}

{\bfseries Proof:}
Assume $(x,y) $ satisfies the sequence chain condition. If $d \in \Gamma(\U)$ and $ \ep > 0$
the chain-length with respect to $d$ of any sequence chain with $\xi = \{ V^d_{\ep/2^n} \}$
from $x$ to $y$ is less than $\ep$. Hence, $(x,y) \in \A_d f$.  As $d$ was arbitrary,
$(x,y) \in \bigcap_{d \in \Gamma} \ \A_d f = \A_{\U} f$.

Now let $(x,y) \in \A_{\U} f$ and $\xi = \{ U_k : k \in \N \}$ be a sequence in $\U$. We must show that
there is a $\xi$ sequence chain from $x$ to $y$.

Let $V_0 = X \times X$.  For $k \in \N$, inductively choose
$V_k = V_k^{-1} \in \U$ such that $V_k \circ V_k \circ V_k \subset V_{k-1} \cap U_k$. By the Metrization Lemma
\cite{K} Lemma 6.12, there exists a pseudo-metric $d \leq 1$ such that  $V_k \subset V^d_{1/2^{k-1}} \subset V_{k-1}$
for $k \in \N$.  It follows that $d \in \Gamma$ and since $V^d_{1/2^{k}} \subset U_k$ it follows that
 if $\xi' = \{ V^d_{1/2^{k}} \}$ then a $\xi'$ sequence chain is a $\xi$ sequence chain.
 It suffices to show that  there is a $\xi'$ sequence chain from $x$ to $y$.

 \begin{lem}\label{wiselem3} Let $\phi : \R \to [0,\infty)$ be given
 by $\phi(0) = 0$ and $\phi(t) = e^{-1/t^2}$ for $t \not= 0$. So that $\phi$
 is a $C^{\infty}$ such that
 \begin{itemize}
 \item[(i)] For all $t > 0$, $\phi'(t) > 0$ and for all $\sqrt{2/3} > t > 0$, $\phi''(t) > 0$.
 \item[(ii)] For $\ep = e^{-3/2}/2$, $\bar d(x,y) = \phi^{-1}(\min(d(x,y),\ep))$ defines a pseudo-metric
 on $X$ with $\U(\bar d) = \U(d) \subset \U$ and so $\bar d \in \Gamma$.
 \item[(iii)] If $\{ \a_k \}$ is a finite or infinite, non-increasing sequence of
 non-negative numbers with $\sum_k \a_k < \phi^{-1}(\ep) < 1$ then
 $\bar d(x,y) \leq \a_k$ implies $d(x,y) < 2^{-k}$, for all $k \in \N$.
 \end{itemize}
 \end{lem}

 {\bfseries Proof:} (i) is an easy direct computation.

 (ii) Observe that if $\psi : [0,a] \to \R$ is $C^2$ with $\psi(0) = 0, \psi'(t) > 0$ and $\psi''(t) < 0$ for $0 < t < a$
 then for all $t,s \leq a/2$, $\psi(t) + \psi(s) - \psi(t+s) \geq 0$, because with $t$ fixed it is true for $s=0$ and
 the derivative with respect to $s$ is positive for $a - t > s > 0$. It follows that if $d$ is a pseudo-metric with $d \leq a/2$
 then  $\psi(d)$ is a pseudo-metric. Clearly, $\U(\psi(d)) = \U(d)$. For (ii) we apply this with $\psi = \phi^{-1}$.

 (iii) Observe that for all $k \in \N$, $\phi(1/k) = e^{-k^2} < 2^{-k}$. Each $\a_k < \phi^{-1}(\ep)$ and so
 $\bar d(x,y) \leq \a_k$ iff $d(x,y) \leq \phi(\a_k)$. If $\phi(\a_k) \geq 2^{-k}$ then for $1 \leq j \leq k$
 $\phi(\a_j) \geq \phi(\a_k) \geq 2^{-k} \geq \phi(1/k)$ and so $\a_j \geq 1/k$ for $j = 1, \dots, k$. Hence,
 $\sum_j \a_j \geq k(1/k) = 1 > \phi^{-1}(\ep)$, contradicting the assumption on the sum.

$\Box$ \vspace{.5cm}

Since $(x,y) \in \A_{\U} f$, there exists $[a,b] \in f^{\times n}$ for some $n \geq 1$  such that
with respect to the metric $\bar d$, the $xy$ chain-length of $[a,b]$ is less than    $\phi^{-1}(\ep)$. Let $b_0 = x$ and
$a_{n+1} = y$. Let $k \mapsto i(k)$ be a bijection on $\{ 1, \dots,n+1\}$ so that the sequence
$\a_k = d(b_{i(k)-1},a_{i(k)})$ is non-increasing. From (iii) it follows that $ (b_{i(k)-1},a_{i(k)}) \in V^d_{2^{-k}}$
for $k = 1, \dots, n+1$ and so $[a,b]$ is a $\xi'$ sequence chain from $x$ to $y$ as required.

$\Box$ \vspace{.5cm}

It is clear that $(\G f)^{-1}$ is the smallest closed, transitive relation
which contains $f^{-1}$. So from  (\ref{eq18}) we obtain:

\begin{equation}\label{eq58aaab}
\G (f^{-1}) \ = \ (\G f)^{-1}, \quad \A_{\U} (f^{-1}) \ = \ (\A_{\U} f)^{-1},
\quad  \CC_{\U} (f^{-1}) \ = \ (\CC_{\U} f)^{-1},
\end{equation}
and so again we may omit the parentheses.

\begin{prop}\label{domain}  For a relation $f$ on a uniform  space $(X,\U)$, the image $f(X)$ is dense in $\CC_{\U} f(X)$ and
the domain $f^{-1}(X)$ is dense in $\CC_{\U} f^{-1}(X)$.\end{prop}

 {\bfseries Proof:} Let $A = \ol{f(X)}$ and let $y \in \CC_{\U} f(x)$. If $U \in \U$ and $[a,b] \in f^{\times n}$ is an $xy, U$ chain, then
 $b_i \in A$ for all $i$ and so $y \in U(A)$. Because $A$ is closed, it equals the intersection $\bigcap_{U \in \U} \ U(A)$. Thus,
 $\CC_{\U} f(X) \subset A$. Replacing $f$ by $f^{-1}$ we obtain the domain result.

$\Box$ \vspace{.5cm}

From (\ref{eq17}) we obtain monotonicity: If $f \subset g$ are relations on $(X,\U)$ then
\begin{equation}\label{eq58ac}
\G f \ \subset \ \G g, \quad \A_{\U} f \ \subset \ \A_{\U} g, \quad  \CC_{\U} f \ \subset \ \CC_{\U} g,
\end{equation}

Again the operators are idempotent.

\begin{prop}\label{prop6ua}
\begin{align}\label{eq20aa}
\begin{split}
f \ \subset \ g \subset \CC_{\U} f \quad &\Longrightarrow \quad \CC_{\U} f = \CC_{\U} g,\\
f \ \subset \ g \subset \A_{\U} f \quad &\Longrightarrow \quad \A_{\U} f = \A_{\U} g,\\
f \ \subset \ g \subset \G f \ \  \quad &\Longrightarrow \quad \G f = \G g.
\end{split}
\end{align}
\end{prop}

{\bfseries Proof:} For any $d \in \Gamma$, \ $f \ \subset \ g \subset \CC_{\U} f \subset \CC_d f$ and so by
(\ref{eq20}) and montonicity, $\CC_d f = \CC_d g$. Intersect over $d \in \Gamma$. The proof for $\A_{\U}$ is
similar.

Finally, if $F$ is a closed, transitive relation then $F = \G F$.

$\Box$ \vspace{.5cm}

Proceeding just as with (\ref{eq20ab}) we see that for relations $f$ and $g$ on $(X,\U)$
\begin{align}\label{eq20ac}
\begin{split}
\CC_{\U} f \cap \CC_{\U} g \ &=  \ \CC_{\U}(\CC_{\U} f \cap \CC_{\U} g), \hspace{1cm} \\
\A_{\U} f \cap \A_{\U} g \ &= \ \A_{\U}(\A_{\U} f \cap \A_{\U} g), \hspace{1cm}\\
\G f \cap \G g \ &=  \ \G(\G f \cap \G g). \hspace{2cm}
\end{split}
\end{align}

If $\U_1$ and $\U_2$ are uniformities on $X$ then
\begin{equation}\label{eq59a}
\U_1 \subset \U_2 \qquad \Longrightarrow \qquad \CC_{\U_2} f \ \subset \CC_{\U_1} f \quad
\text{and} \quad \A_{\U_2} f \ \subset \A_{\U_1} f.
\end{equation}

More generally, we have

\begin{prop}\label{prop6ub} If $h : (X_1, \U_1) \to (X_2, \U_2)$ is a continuous map which  maps the relation
$f_1$ on $X_1$ to $f_2$ on $X_2$, then $h$ maps $\G f_1$ to $\G f_2$.
If, in addition, $h$ is uniformly continuous, then $h$ maps  $\CC_{\U_1} f_1$ to $\CC_{\U_2} f_2$, and
maps  $\A_{\U_1} f_1$ to $\A_{\U_2} f_2$.\end{prop}

{\bfseries Proof:} If $h$ is continuous then, $(h \times h)^{-1}(\G f_2)$ is a
closed, transitive relation which contains $f_1$ and so contains $\G f_1$.

Now assume that $h$ is uniformly continuous.
Let  $d_2 \in \Gamma(\U_2)$. By uniform continuity, $d_1 = h^*d_2 \in \Gamma(\U_1)$, where
\begin{equation}\label{eq59aaa}
h^*d_2(x,y) = d_2(h(x),h(y)). \hspace{2cm}
\end{equation}\index{h@$h^*d$}\index{d@$h^*d$}
 Thus, $h : (X_1,d_1) \to (X_2,d_2)$ is Lipschitz. In fact, it is an isometry.
By Proposition  \ref{prop00c}, $h$ maps $\A_{\U_1} f_1 \subset \A_{d_1} f$ into $\A_{d_2} f$ and similarly for
$\CC$. Intersect over all $d_2 \in \Gamma(\U_2)$.

$\Box$ \vspace{.5cm}

For a relation $f$ on $X$ let $f^{[1,k]} = \bigcup_{j=1}^k \ f^j$\index{f@$f^{[1,k]}$} for
any positive integer $k$. Let $f^{[0,k]} = 1_X \cup f^{[1,k]}$.
 If $d$ is a pseudo-metric on $X$
and $f$ is a map on $X$ we let $d^k = \max_{j = 0}^k \ (f^j)^*d$.\index{d@$d^k$} Let $d^0 = d$.

\begin{cor}\label{cor6ubb} Let $k \geq 2$ be an integer and $f$ be a continuous map on a uniform space $(X,\U)$.
\begin{equation}\label{eq59aab}
\G f \ = \  f^{[1,k-1]}  \ \cup \ \G (f^k) \circ f^{[0,k-1]}, \hspace{2cm}
\end{equation}
and $|\G (f^k)| = |\G f|$.

If $f$ is a uniformly continuous map, then
\begin{equation}\label{eq59aac}
\A_{\U} f \ = \  f^{[1,k-1]}  \ \cup \ \A_{\U} (f^k) \circ f^{[0,k-1]},
\quad \CC_{\U} f \ = \  f^{[1,k-1]}  \ \cup \ \CC_{\U} (f^k) \circ f^{[0,k-1]},
\end{equation}
and $|\A_{\U}(f^k)| = |\A_{\U} f|, \quad |\CC_{\U}(f^k)| = |\CC_{\U} f| $.
\end{cor}

{\bfseries Proof:} If $F$ is a closed relation on $X$ and $f$ is a continuous map
on $X$ then $F \circ f$ is a closed relation.
For suppose $\{ (x_i,y_i) \}$ is a net in $F \circ f$ converging to $(x,y)$. Then
$\{ f(x_i) \}$ converges to $f(x)$ by continuity
and $\{ (f(x_i),y_i) \}$ is a net in $F$ converging to $(f(x),y)$.  Since $F$ is
closed, $(f(x),y) \in F$ and $(x,y) \in F \circ f$.

Hence, $ f^{[1,k-1]}  \ \cup \ \G (f^k) \circ f^{[0,k-1]}$ is a closed relation
which contains $f$. Since $f \subset \G f$, transitivity
of $\G f$ implies that $f^k \subset \G f$. Hence, $\G (f^k) \subset \G f$. Transitivity again implies $ f^{[1,k-1]}  \ \cup \ \G (f^k) \circ f^{[0,k-1]} \subset \G f$.

Because $f$ maps $f^k$ to $f^k$ it follows from
Proposition \ref{prop6ub} that it maps $\G (f^k)$ to itself.  Hence,
$f^{[0,k-1]} \circ \G (f^k)  \subset \G (f^k) \circ f^{[0,k-1]}$.
Furthermore, $f^{[0,k-1]} \circ f^{[1,k-1]} \subset  f^{[1,k-1]} \cup f^k \circ f^{[0,k-1]}$. It follows that
$ f^{[1,k-1]}  \ \cup \ \G (f^k) \circ f^{[1,k-1]}$ is transitive and so contains
$\G f$ since it is closed and contains $f$.

It clearly, follows that $|\G (f^k)| \subset |\G f|$.  Assume that $x \in |\G f|$. From (\ref{eq59aab}) it follows that
either $x \in f^j(x)$ for some $j \in [1,k-1]$ or $x \in \G (f^k) \circ f^j(x)$ for some $j \in [0,k-1]$.
If $x = f^j(x)$ then $x = (f^j)^k(x) = (f^k)^j(x)$ and so $x \in \G (f^k)(x)$.  Similarly, since $f^j$ maps $\G(f^k)$ to itself,
$$(\G (f^k) \circ f^j)^k \ \subset \ (\G (f^k))^k \circ (f^j)^k \ \subset \ \G (f^k)$$
and so  $x \in |\G (f^k)|$ if $x \in \G (f^k) \circ f^j(x)$.

Transitivity again implies $f^k \subset \A_{\U} f \subset \CC_{\U} f$, and so monotonicity and transitivity imply
\begin{align}\label{eq59aabb}
\begin{split}
\A_{\U} f \ \supset \  &f^{[1,k-1]}  \ \cup \ \A_{\U} (f^k) \circ f^{[0,k-1]}, \hspace{2cm}\\
 \CC_{\U} f \ \supset \  &f^{[1,k-1]}  \ \cup \ \CC_{\U} (f^k) \circ f^{[0,k-1]}.\hspace{2cm}
\end{split}
\end{align}

Now assume that $f$ is a uniformly continuous map.  Notice that if
$[a,b] \in f^{\times n}$ then $b_i = f(a_i)$ for $i = 1,\dots, n$.
Observe that if $x \in X$ and $j \leq k$
\begin{equation}\label{eq59aad}
\begin{split}
 d(f^j(a_1),a_{j+1}) \ \leq \ \Sigma_{i=1}^{j} \ d(f^{j-i+1}(a_i),f^{j-i}(a_{i+1})) \hspace{2cm} \\
  \leq \  \ \Sigma_{i=1}^{j-1} \ d^k(f(a_i),a_{i+1}), \hspace{4cm} \\
\text{and} \quad d(f^j(x),f^j(a_1)) \ \leq \ d^k(x,a_1). \hspace{3cm}
\end{split}
\end{equation}

Let $(x,y) \in \A_{\U}$. For $\a = (d,\ep) \in \Gamma \times (0,\infty)$ there
exists $[a,b]_{\a} \in f^{n_{\a}}$ with $xy$ chain-length
with respect to $d$ less than $\ep$.

If $n_{\a} < k$ frequently then for some $j \in [1,k-1]$ frequently $n_{\a} = j$ and it follows from
continuity of $f$ that $y = f^j(x)$.

Instead assume that eventually $n_{\a} \geq k$.  If $\ep > 0$ and
 $d_1 \in \Gamma(\U)$, there exists $d \geq d_1$ and
$[a,b] \in f^{\times n}$ with $n \geq k$ so that the $xy$ chain-length
of $[a,b]$ with respect to $d^k$ is less than $\ep$.
Let $n = j + qk$ with $j \in [0,k-1]$ and $q \geq 1$. The sequence

\begin{equation}
\begin{split}
[a,b]^k = (a_{j+1},f^k(a_{j+1})),(a_{j+k+1},f^k(a_{j+k+1}))\dots \\
(a_{j+(q-1)k+1},f^k(a_{j+(q-1)k+1})) \in (f^k)^{\times q},
\end{split}
\end{equation}
and with
$y = a_{n_{\a} + 1}$, (\ref{eq59aad})  implies that
the $f^j(x)y$ chain-length with respect to $d$ and so with respect to $d_1$ is less than $\ep$. Since $d_1$ was
arbitrary it follows that  $y \in \G (f^k) \circ f^{[0,k-1]}(x)$.

For $\CC_{\U} f$ we proceed as before, but use chain-bound less than $\ep/k$.

For $|\A_{\U}(f^k)|$ and $|\CC_{\U}(f^k)|$ we use the same argument as for $|\G (f^k)|$ above.

$\Box$ \vspace{.5cm}

If a real-valued function on $X$ is uniformly continuous with respect to some $d \in \Gamma(\U)$ then it is uniformly
continuous from $(X,\U)$. In particular, for every $d \in \Gamma(\U)$ and $f \subset X \times X$, the functions
$\ell_d^f$ and $m_d^f$ are uniformly continuous from $(X \times X, \U \times \U)$. It follows that the sets
$\CC_{\U} f, \A_{\U} f \subset X \times X$ and $ |\CC_{\U} f|, |\A_{\U} f| \subset X$ are closed.

 As before, a Lyapunov function for a  relation $f$ on a uniform space $(X,\U)$ is a
continuous map $L : X \to \R$ such that
$ (x,y) \in f $ implies  $L(x) \leq L(y)$. Hence, the relation $\ \leq_L \ \subset \ X \times X$ is closed.

As in Definition \ref{def6}

\begin{df}\label{def6u}  Let $F$ be a  closed, transitive relation on a Hausdorff
uniform space $(X,\U)$ and let $\L$ be a collection
of Lyapunov functions for $F$.
 We define three conditions on $\L$.
\begin{itemize}
\item[ALG] If $L_1,L_2 \in \L$ and $c \geq 0$ then \\ $L_1 + L_2, \max(L_1,L_2), \min(L_1,L_2), c L_1, c, -c \in \L$.\index{condition!ALG}

\item[CON] For every sequence $\{ L_k \}$ of elements of $\L$ there exists a summable sequence of positive real numbers
$\{ a_k \}$ such that $\Sigma_k a_k L_k$ converges uniformly to an element of $\L$.\index{condition!CON}

\item[POIN] If $(x,y) \not\in 1_X \cup F$ then there exists $L \in \L$ such that $L(y) < L(x)$, i.e.\
$1_X \cup F \ = \ \bigcap_{L \in \L} \ \leq_L$.\index{condition!POIN}

\end{itemize}
\end{df}
\vspace{.5cm}

\begin{theo} \label{theo13u} Let $f$ be a  relation on a Hausdorff uniform space $(X,\U)$ with gage $\Gamma$.

(a) let $\L_{\ell}$ be the set of  bounded, uniformly
continuous functions which are
$K \ell_d^f$ dominated for some $d \in \Gamma$ and some positive $K$.
Each $L \in \L_{\ell}$ is a $\A_{\U} f$ Lyapunov function and so satisfies
\begin{equation}\label{eq50ua}
\A_{\U} f \ \subset \ \leq_L \qquad \mbox{and} \qquad |\A_{\U} f| \ \subset \ = |L|_{\A_{\U} f}.
\end{equation}
The collection $\L_{\ell}$ satisfies the conditions ALG, CON,  and POIN  with respect to $F = \A_{\U} f$.

(b) let $\L_{m}$ be the set of  bounded, uniformly
continuous functions which are
$K m_d^f$ dominated for some $d \in \Gamma$ and some positive $K$.
Each $L \in \L_{m}$ is a $\CC_{\U} f$ Lyapunov function and so satisfies
\begin{equation}\label{eq50ub}
\CC_{\U} f \ \subset \ \leq_L \qquad \mbox{and} \qquad |\CC_{\U} f| \ \subset \ = |L|_{\CC_{\U} f}.
\end{equation}
The collection $\L_{m}$ satisfies the conditions ALG, CON,  and POIN  with respect to $F = \CC_{\U} f$.
\end{theo}

{\bfseries Proof:} If $\{ d_k \}$ is a sequence in $\Gamma$ and $K_k \geq 1$ so that $d_k \leq K_k$ and $\{ K_k a_k \}$
is a summable sequence of positive reals, then
by Lemma \ref{applem01a} $ d = \Sigma_k \ (a_k) d_k \in \Gamma$. Furthermore,
\begin{equation}\label{eq50aa}
a_k \ell^{d_k}_f \ = \ \ell^{a_k d_k}_f \ \leq \ \ell^d_f \quad \text{and}
\quad  a_k m^{d_k}_f \ = \ m^{a_k d_k}_f \ \leq \ m^d_f.
\end{equation}

So if $L$ is $K \ell^{d_k}_f$ dominated then it is $(K/a_k) \ell^{d}_f$ dominated. Thus, if $\{L_k \}$ is a sequence
in $\L_{\ell}$ we can choose $d \in \Gamma$ such that each $L_k$ is $K_k \ell^d_f$ dominated for some $K_k$.
Then ALG and CON follow for $L_{\ell}$ from Theorem \ref{theo13} for $(X,d)$.

Now assume that $(x,y) \not\in 1_X \cup \A_{\U} f$. Because $X$ is Hausdorff there exists $d_1 \in \Gamma$ such that
$d_1(x,y) > 0$. There exists $d_2 \in \Gamma$ such that $(x,y) \not\in \A_{d_2} f$. Let $d = d_1 + d_2$.
Since $\ell^{d_2}_f \leq \ell^d_f$ it follows that $(x,y) \not\in Z_d \cup \A_{d} f$. From
Theorem \ref{theo13} again there exists a function $L$  which is $d$ uniformly continuous, $K \ell^d_f$
dominated for some $K$ and satisfied $L(x) > L(y)$.  Hence, $L \in \L_{\ell}$ with $L(x) > L(y)$, proving POIN.

The results in (b) for $\L_m$ are proved exactly the same way with Theorem \ref{theo13} replaced by Theorem
\ref{theo15a}.

$\Box$ \vspace{.5cm}

\begin{theo}\label{theo14u}  Let $f$  be a  relation on a uniform space $(X,\U)$.

If $L$ is a Lyapunov function for $f$, then $L$ is a Lyapunov function for $\G f$.

If $L$ is a uniformly continuous Lyapunov function for $f$, then $L$ is a Lyapunov function for $\A_{\U} f$.
 \end{theo}

{\bfseries Proof:} If $L$ is a Lyapunov function for $f$ then, by
continuity of $L$,  $\leq_L$ is a closed, transitive relation which contains $f$ and so
contains $\G f$.

If $L$ is bounded and uniformly continuous, then $d_L(x,y) = |L(x) - L(y)|$ is a pseudo-metric in $\Gamma(\U)$.
Let $(x,y) \in \A_{\U} f$ and $\ep \in (0,1)$. There exists $[a,b] \in f^{\times n}$ such that
the $xy$ chain-length of $[a,b]$ with respect to $d_L$ is less than $\ep$. Since $L$ is a Lyapunov function for $f$,
we have that $L(a_i) \leq L(b_i)$ for $i = 1, \dots, n$.
\begin{equation}\label{50ac}
\begin{split}
L(y) - L(x) \ = \ L(y) - L(b_n) + \Sigma_{i = 1}^n \ L(b_i) - L(a_i) + \hspace{1cm} \\
 \Sigma_{i = 1}^{n-1} L(a_i) - L(b_{i+1}) + L(a_1) - L(x).
\end{split}
\end{equation}
The first sum is non-negative and the rest has absolute value at most the chain-length. Hence,
$L(y) - L(x) \geq - \ep$. Since $\ep$ was arbitrary, $L(y) - L(x) \geq 0$.

If $L$ is unbounded then for each positive $K$, $L_K = \max (\min(L,K),-K)$ is a bounded, uniformly continuous Lyapunov function
and so is an $\A_{\U} f$ Lyapunov function.  If $(x,y) \in \A_{\U} f$ then by choosing $K$ large enough we have
$L_K(x) = L(x)$ and $L_K(y) = L(y)$.  So $L(y) - L(x) = L_K(y) - L_K(x) \geq 0$.

$\Box$ \vspace{.5cm}

\begin{cor}\label{cor15u} Let $f$ be a relation on a Tychonoff space $X$ and let $\U_M$ be the maximum uniformity
compatible with the topology. Let $\L$ be the set of all bounded, Lyapunov functions for $f$. Each $L \in \L$ is a
Lyapunov function for $\A_{\U_M} f$ and
\begin{equation}\label{50ad}
1_X \cup \A_{\U_M} f \ = \ \bigcap_{L \in \L} \ \leq_L
 \end{equation}
 \end{cor}

{\bfseries Proof:} With respect to the maximum uniformity every continuous real-valued function is uniformly continuous.
So every $L \in \L$ is a Lyapunov function for $\A_{\U_M} f$ by Theorem \ref{theo14u}. Hence
$1_X \cup \A_{\U_M} f \ \subset \ \bigcap_{L \in \L} \ \leq_L$.  The reverse inclusion follows from POIN in
Theorem \ref{theo13u} (a).

$\Box$ \vspace{.5cm}

\begin{theo}\label{theo8u}
Let $F$ be a closed, transitive relation on a Hausdorff uniform space $(X,\U)$ whose topology is second countable.
Let $\L$ be a collection of
Lyapunov functions for $F$ which satisfies ALG, CON and POIN. There exists a sequence $\{ L_k \}$ in $\L$
such that
\begin{equation}\label{eq44u}
\bigcap_k \ \leq_{L_k} \quad = \quad 1_X \cup F.
\end{equation}
If $\{ a_k \}$ is a positive, summable sequence such that
$L = \Sigma_n \ a_k L_k  \in \L$ then $L$ is a Lyapunov function for $F$ such that
$1_X \cup F \ = \ \leq_L$ and
\begin{equation}\label{eq45uaa}
x \in F(y) \qquad \Longrightarrow \qquad L(y) < L(x) \quad \mbox{unless} \ y \in F(x)
\end{equation}
In particular,
\begin{equation}\label{eq46uaa}
|L|_F \ = \  |F|
\end{equation}
\end{theo}

{\bfseries Proof:} Proceed just as in the proof of Theorem \ref{theo8}
using  the fact that $(X \times X) \setminus (1_X \cup F)$
is Lindel\"{o}f.

$\Box$ \vspace{.5cm}

For a metrizable space $X$ we let $\Gamma_m(X)$\index{g@$\Gamma_m(X)$} be the set of metrics compatible with the topology on $X$.

\begin{theo}\label{theo9u} Let $f$ a a relation on a Hausdorff uniform space $(X,\U)$ whose topology is
second countable. There exist  bounded, uniformly continuous Lyapunov functions $L_{\ell}, L_{m}$ for $f$ such that

\begin{equation}\label{eq45u}
\begin{split}
1_X \cup \A_{\U} f \ = \ \leq_{L_{\ell}}, \quad 1_X \cup \CC_{\U} f \ = \ \leq_{L_{m}} \qquad  \text{and}, \\
x \in \A_{\U} f (y) \qquad \Longrightarrow \qquad L_{\ell}(y) < L_{\ell}(x) \quad \mbox{unless} \ y \in \A_{\U} f (x),\\
x \in \CC_{\U} f (y) \qquad \Longrightarrow \qquad L_m(y) < L_m(x) \quad \mbox{unless} \ y \in \CC_{\U} f (x)
\end{split}
\end{equation}
In particular,
\begin{equation}\label{eq46u}
|L_{\ell}|_{\A_{\U} f} \ = \  |\A_{\U} f|, \quad \text{and} \quad |L_{m}|_{\CC_{\U} f} \ = \  |\CC_{\U} f|
\end{equation}

Furthermore, there exists a metric $d \in \Gamma_m(X) \cap \Gamma(\U)$
 such that $L_{\ell}$ and $L_{m}$ are Lipschitz functions on
$(X,d)$ and
\begin{equation}\label{eq47u}
\A_{\U} f \ = \ \A_{d} f \quad \text{and} \quad \CC_{\U} f \ = \ \CC_{d} f.
\end{equation}
\end{theo}

{\bfseries Proof:} The  pseudo-metrics chosen below are all assumed bounded by 1. We can always replace $d$ by $\min(d,1)$.

 We apply Theorem \ref{theo8u} to $\L_{\ell}$ and $\A_{\U} f$ and to $\L_m$ and $\CC_{\U} f$
and obtain $L_{\ell} \in \L_{\ell}$ and $L_{m} \in \L_{m}$ which satisfy
(\ref{eq45u}) and  (\ref{eq46u}). We may assume that each maps to $[0,1]$. In particular, there
exist $d_1, d_2 \in \Gamma(\U)$ and positive $K_1, K_2$ so that $L_{\ell}$ is $K_1 \ell^{d_1}_f$ dominated and
$L_m$ is $K_2 m^{d_2}_f$ dominated.

Let $\B$ be a countable base and $D$ be a countable dense subset of $X$.  For each pair $(x,U)$ with $U \in \B$ and
$x \in U \cap D$ there
exists $d  = d_{(x,U)} \in \Gamma(\U)$ and a rational $\ep > 0$ such that the ball $V^d_{\ep}(x) \subset U$

For each $x \not\in |\A_{\U} f|$ there exists $d_{x,1} \in \Gamma(\U)$ such that $\ell^{d_{x,1}}_f (x,x) > 0$  and
for each $x \not\in |\CC_{\U} f|$ there exists $d_{x,2} \in \Gamma(\U)$ such that $m^{d_{x,2}}_f (x,x) > 0$.
These
are open conditions and so we can choose a sequence $\{ d_3, d_4, \dots \}$ in $\G$ and a positive
sequence $\{ a_1, a_2, \dots \}$ with sum $= 1$ so that $d$ defined by
$d(x,y)  = \frac{1}{3}[|L_{\ell}(x) - L_{\ell}(y)| + |L_{m}(x) - L_{m}(y)| + \Sigma_{i = 1}^{\infty} a_i d_i] $  satisfies

\begin{itemize}
\item[(i)] $d \in \Gamma(\U)$.
\item[(ii)] The $\U(d)$ topology is that of $X$, i.e.\ $d \in \Gamma_m(X)$.
\item[(iii)] $x \not\in |\A_{\U} f|$ implies $\ell^d_f(x,x) > 0$, and $x \not\in |\CC_{\U} f|$ implies $m^d_f(x,x) > 0$.
\item[(iv)] There exist positive $K_{\ell}$ and $K_m$ so that $L_{\ell}$ is $K_{\ell}\ell^d_f$ dominated and
$L_{m}$ is $K_{m}m^d_f$ dominated.
\end{itemize}

Condition (i) follows from Lemma \ref{applem01a}. Condition
(ii) implies that $d$ is a metric since $X$ is Hausdorff. From condition (iv)
and (\ref{eq45u})
we obtain
\begin{equation}\label{48u}
\begin{split}
1_X \cup \A_{d} f \ \subset \ \leq_{L_{\ell}} \ = \ 1_X \cup \A_{\U} f, \\
1_X \cup \CC_{d} f \ \subset \ \leq_{L_{m}} \ = \ 1_X \cup \CC_{\U} f.
\end{split}
\end{equation}
On the other hand, $d \in \Gamma(\U)$ implies $\A_{\U} f \subset \A_{d} f $ and $\CC_{\U} f \subset \CC_{d} f $.
Hence, if $(x,y) \in \A_{d} f \setminus \A_{\U} f$ then $(x,y) \in 1_X$ and so $\ell^d_f(x,x) = 0$. By condition
(iii) this implies $x \in |\A_{\U} f|$ and so $(x,y) = (x,x) \in \A_{\U} f$. This contradiction proves the first equation
in (\ref{eq47u}).  The second follows similarly.

Clearly, $L_{\ell}$ and $L_{m}$ are Lipschitz with Lipschitz constant at most $3$.

$\Box$ \vspace{.5cm}

If  $\U_M$  the maximum uniformity
compatible with the topology for a metrizable space $X$, then since such a space is paracompact, $\U_M$ consists of all
neighborhoods of the diagonal. The gage $\Gamma(\U_M)$ consists of all pseudo-metrics which are continuous on $X$. In particular,
$\Gamma_m(X) \subset \Gamma(\U_M)$.

\begin{cor}\label{cor10u}  Let $f$ be a relation on a second countable Tychonoff
space $X$ and let $\U_M$ be the maximum uniformity
compatible with the topology. There exists
a metric $d_0 \in \Gamma_m(X)$ such that
\begin{equation}\label{eq49u}
\A_{\U_M} f \ = \ \A_{d_0} f \quad \text{and} \quad \CC_{\U_M} f \ = \ \CC_{d_0} f.
\end{equation}
Furthermore,
\begin{equation}\label{eq50u}
\A_{\U_M} f \ = \ \bigcap_{d \in \Gamma_m(X)} \ \A_{d} f \quad \text{and} \quad \CC_{\U_M} f \ = \ \bigcap_{d \in \Gamma_m(X)} \ \CC_{d} f.
\end{equation}
\end{cor}

{\bfseries Proof:} 
A second countable
Hausdorff space is metrizable, i.e.\ there exists
a metric $\bar d$ with the $\U(\bar d)$ topology that of $X$. Thus,
$\bar d \in \Gamma_m(X) \subset \Gamma(\U_M)$.  If $d_0\in \Gamma(\U_M)$, then $d = \bar d + d_0$ is a metric in $\Gamma(\U_M)$ and so is continuous.
Since $d \geq \bar d$ it follows that the $\U(d)$  topology is that of $X$ as well, i.e.\ $d \in \Gamma_m(X)$.  Furthermore,
\begin{equation}\label{eq51u}
\begin{split}
\A_{\U_M} f \ \subset \ \A_{d} f \  \subset \ \A_{d_0} f \\
\CC_{\U_M} f \ \subset \ \CC_{d} f \  \subset \ \CC_{d_0} f.
\end{split}
\end{equation}
 Hence, the intersection over $\Gamma_m(X)$ yields the
same result as intersecting over the entire gage, $\Gamma(\U_M)$. Furthermore, if $d_0$ is a metric in $\Gamma(\U)$ satisfying (\ref{eq47u})
then (\ref{eq47u}) together with (\ref{eq51u}) implies (\ref{eq49u}).

$\Box$ \vspace{.5cm}

For $d$  a metric on $X$, $\U(d)$ is the uniformity generated by $V^d_{\ep}$ for all $\ep > 0$.
We say that $d$ generates the uniformity
$\U(d)$ and that $\U$ is metrizable if $\U = \U(d)$ for some metric $d$. The Metrization Theorem, Lemma 6.12 of
\cite{K}, implies that a Hausdorff uniformity is metrizable iff it is countably generated. Two metrics $d_1$ and $d_2$
generate the same uniformity exactly when they are uniformly equivalent. That is, the identity maps between
$(X,d_1)$ and $(X,d_2)$ are uniformly continuous. For a metrizable uniformity $\U$ we
let $\Gamma_m(\U) = \{ d : d $ is a metric
with $\U(d) = \U \}$. \index{g@$\Gamma_m(\U)$}

If $(X,d)$ is a metric space and the set of non-isolated points is not compact,
then the maximum uniformity $\U_M$ is not metrizable even if $X$ is second countable.
Since a metric space is paracompact, $\U_M$
consists of all neighborhoods of the diagonal. By hypothesis there is a sequence
$\{ x_1, x_2, \dots \}$ of distinct non-isolated
points with no convergent subsequence and so we can choose open sets $G_i$ pairwise disjoint and with $x_i \in G_i$.
We can choose $y_i \in G_i \setminus \{ x_i \}$
such that $\ep_i = d(x_i,y_i) \to 0$ as $i \to \infty$ and let $\ep_0 = 1$. Let
$G_0$ be the complement of a closed neighborhood of $\{ x_i \}$ in $\bigcup_{i=1}^{\infty} G_i$. Thus,
$\{ G_i \}$ is a locally finite open cover. Choose $\{ \phi_i \}$ a partition of unity, i.e.\ each $\phi_i$ is a continuous
real-valued function with support in $G_i$ and with $\Sigma_i \phi_i = 1$. Define $\psi(x) = \Sigma_i \ep_i \phi_i(x)/2$.
In particular, $\psi(x_i) = \ep_i/2$ for $i = 1, 2, \dots$. Thus, $\psi$ is a continuous, positive
function with infimum 0. So
$U = \{ (x,y) : d(x,y) < \psi(x) \}$ is a neighborhood of the diagonal disjoint from $\{ (x_i,y_i) : i = 1,2, \dots \}$.
But if $\ep_i < \ep$ then $(x_i,y_i) \in V^d_{\ep}$. It follows that for any metric $d$ compatible with the topology of $X$ there
exists a neighborhood of the diagonal, and so an element of $\U_M$, which is not in $\U(d)$.

\begin{theo}\label{theo11u} Let $(X,\U)$ be a  uniform space with $\U$ metrizable and let $f$ be a relation on $X$.

(a) For every $d \in \Gamma_m(\U)$, $\CC_{\U} f = \CC_d f$.

(b) $\A_{\U} f = \bigcap_{d \in \Gamma_m(\U)} \ \A_d f$.
\end{theo}

{\bfseries Proof:}  If $\bar d \in \Gamma(\U)$ and $d_1 \in \Gamma_m(\U)$ then $d = \bar d + d_1 \in \Gamma_m(\U)$ and
  $\CC_{d} f \subset \CC_{\bar d} f$. Thus,
we need only intersect over $\Gamma_m(\U)$ to get  $\CC_{\U} f$. Similarly, for $\A_{\U} f$.

On the other hand, if $d_1, d_2 \in \Gamma_m(\U)$ then $d_1$ and $d_2$ are uniformly equivalent metrics and so Proposition \ref{prop00c} implies that
$\CC_{d_1} f = \CC_{d_2} f$.  Hence, the intersection $\CC_{\U} f$ is this common set.

$\Box$ \vspace{1cm}

There are special constructions for the Conley relations.

  \begin{df}\label{def16ag} Let $f$ be a relation on a uniform space $(X,\U)$.
  \begin{enumerate}
  \item[(a)] A set $A \subset X$ is called $\U$
\emph{ inward }\index{subset!$\U$ inward}if there exists $U \in \U$ such that $U(f(A)) \subset A$, or, equivalently, if there
 exist $d \in \Gamma(\U)$ and $\ep > 0$ such that $A$ is $(V^d_{\ep} \circ f)$  $^+$invariant.

 \item[(b)] A $\U$ uniformly continuous function $L : X \to [0,1]$ is called a $\U$ \emph{ elementary Lyapunov function}\index{Lyapunov function!elementary}
 \index{elementary Lyapunov function} for $f$ if $(x,y) \in f$ and $L(x) > 0$ imply $L(y) = 1$.
 \end{enumerate}
 \end{df}
 \vspace{.5cm}

 If $\U = \U_M$ for the space $X$, then a $\U$ inward set $A$ for $f$ is just called an inward set for $f$. For a paracompact Hausdorff space
 any neighborhood of a closed set is a $\U_M$ uniform neighborhood and so a set $A$ is inward for a relation $f$ on such a space iff
 $\ol{f(A)} \subset A^{\circ}$.  A continuous function $L : X \to [0,1]$
 is $\U_M$ uniformly continuous and we will call a $\U_M$ elementary Lyapunov function just an elementary Lyapunov function.

 Observe for $L : X \to [0,1]$ that if $L(x) = 0$ or $L(y) = 1$ then $L(y) \geq L(x)$.  So an elementary Lyapunov function is a Lyapunov function.
 In addition, the points of $G_L = \{ x : 1 > L(x) > 0 \}$ are regular points for $L$ and so $|L|_f \subset L^{-1}(0) \cup L^{-1}(1)$ with equality if
 $f$ is a surjective relation.

 If $u : X \to \R$ is a bounded real-valued function we define the pseudometric $d_u$\index{d@$d_u$} on $X$ by $d_u(x,y) = |u(x) -u(y)|$.
 If $u$ is uniformly continuous on $(X,\U)$ then $d_u \in \Gamma(\U)$.

 \begin{theo}\label{theo16ai} Let $f$ be a relation on a uniform space $(X,\U)$.
 \begin{enumerate}
 \item[(a)] If $A$ is a $\U$ inward subset for $f$ then there exist  $d \in \Gamma(\U)$ and $\ep > 0$ such that
 $\ol{V^d_{\ep}(f(A))} \subset A^{\circ}$. In particular, $A_1 = A^{\circ}$ and $A_2 = \ol{V^d_{\ep}(f(A))} $
 are $\U$ inward with $A_1$ open, $A_2$ closed and $\ol{f(A)} \subset A_2 \subset A_1 \subset A$.

 \item[(b)] Let $A$ be an open $\U$ inward subset for $f$. If for  $d \in \Gamma(\U)$ and $\ep > 0$ $V^d_{\ep}(f(A)) \subset A$, then
 $V^d_{\ep}(\CC_{\U} f(A)) \subset A $. In particular,  $A$ is a $\U$ inward subset of $X$ for $\CC_{\U}f$ and is
 and is $(V^d_{\ep} \circ \CC_{\U} f) $ $^+$invariant.

 \item[(c)] If $A$ is  a $\U$ inward subset for $f$, then there exists $B$ a closed $\U$ inward subset for $f^{-1}$ such that
 $A^{\circ} \cup B^{\circ} = X$ and $B \cap f(A) = \emptyset = A \cap f^{-1}(B)$.

  \item[(d)] If $A$ is a $\U$ inward subset of $X$, then there exists a $\U$ uniformly continuous elementary Lyapunov function $L$ for $f$
  such that $L^{-1}(0) \cup A = X$ and $f(A) \subset L^{-1}(1)$.

  \item[(e)] If $L$ is a $\U$  elementary Lyapunov function for $f$ and $1 \geq \ep > 0$, then $A = \{x : L(x) > 1 - \ep \}$ is an open
  set such that
  \begin{equation}\label{eq16ai}
  \begin{split}
  f(A) \subset \CC_{\U} f(A) \subset \CC_{d_L} f(A) \subset L^{-1}(1), \hspace{2cm}\\
  V^{d_L}_{\ep}(f(A)) \subset   V^{d_L}_{\ep}(\CC_{\U} f(A)) \subset   V^{d_L}_{\ep}(\CC_{d_L} f(A)) \subset A. \hspace{1cm}
  \end{split}
  \end{equation}
  In particular, $L$ is a $\U(d_L)$  elementary Lyapunov function
 for  $\CC_{d_L} f $ and hence is a $\U$  elementary Lyapunov function for $\CC_{\U} f $ and for $\ol{f}$.

   \item[(f)]  If $L$ is a $\U$ elementary Lyapunov function for $f$, then $1 - L$ is a $\U$ elementary Lyapunov function for $f^{-1}$.
  \end{enumerate}
  \end{theo}

  {\bfseries Proof:} (a) There exist  $d \in \Gamma$ and $\ep > 0$ such that $V^d_{2\ep}(f(A))$ is contained in $A$ and so is contained in $A^{\circ}$.
  For a subset $B$ of $X$, $x \in \ol{B}$ implies $d(x,B) = 0$ and so $\ol{V^d_{\ep}(f(A))} \subset V^d_{2\ep}(f(A))$ and
  $\ol{f(A)} \subset V^d_{\ep}(f(A))$.

  (b) Assume that $x \in A$ and $z \in V^d_{\ep}(\CC_{\U} f(x))$. So there exist $z_1 \in \CC_{\U} f(x)$ and $\ep_1 > 0$ such that
  $d(z_1,z) < \ep$. There exist $d_1 \in \Gamma$ and $\ep_1 > 0$ such that $V^{d_1}{\ep_1}(x) \subset A$ and $d(z_1,z) + \ep_1 < \ep$.
  Let $\bar d = d + d_1$. There exists $[a,b] \in f^{\times n}$ such that the $xz_1$ chain-bound of $[a,b]$ with respect to $\bar d$
  is less than $\ep_1$. Because $d_1(x,a_1) < \bar d(x,z_1) < \ep_1$, $a_1 \in A$. Since $b_1 \in f(A)$ and $d(b_1,a_2) < \bar d(b_1,a_2) < \ep$,
  $a_2 \in A$.  Inductively, we obtain $a_i \in A$ and $b_i \in f(A)$ for $i = 1, \dots, n$. Finally,
  $d(b_n,z) \leq \bar d(b_n,z_1) + d(z_1,z) < \ep$. So $z \in A$.

  (c) Let $d \in \Gamma$ and $\ep > 0$ be such that $V^d_{2\ep}(f(A))$ is contained in $A$ and so is contained in $A^{\circ}$.
  Let $B = X \setminus V^d_{\ep}(f(A))$ so that $B^{\circ} = X \setminus \ol{V^d_{\ep}(f(A))}$. Thus, $B$ is closed,
  $A^{\circ} \cup B^{circ} = X$ and $B \cap f(A) = \emptyset$. Assume that $(x,y) \in f$ and $z \in V_{\ep}(x)$. If $y \in B$ then $x \not\in A$
  and so $x \not\in V^d_{2\ep}(f(A))$ and $z \not\in V^d_{\ep}(f(A))$. That is, $z \in B$. Thus, $V^d_{\ep}(f^{-1}(B)) \subset B$.
  Finally, if $y \in B$ then $y \not\in f(A)$ and so $x \not\in A$. That is, $f^{-1}(B) \cap A = \emptyset$.

  (d) Assume that $V^d_{\ep}(f(A)) \subset A$. Let $L(x) = \max (\ep - d(x,f(A)),0)/\ep$. If $(x,y) \in f$ and $L(x) > 0$ then
  $d(x,f(A)) < \ep$ and so $x \in A$. Then $y \in f(A)$ implies $L(y) = 1$.

  (e) Clearly, $f(A) \subset L^{-1}(1)$. Let $\ep > \ep_1 > 0$. We show that $ V^{d_L}_{\ep_1}(\CC_{d_L} f(A)) \subset \{ y : L(y) > 1 - \ep_1 \}$.
  Assume $x \in A, y \in V^{d_L}_{\ep_1}(\CC_{d_L} f(x)$. So there exists $z \in \CC_{d_L}f(x)$ with $d_L(z,y) < \ep_1$.
  Choose $\ep_2 > 0$ so that $d_L(z,y) + \ep_2 < \ep_1$ and $L(x) > 1 - \ep + \ep_2 $. Since $L$ is uniformly continuous, $d_L \in \Gamma(\U)$ and so
  there exists $[a,b] \in f^{\times n}$ such that the $xz$ chain-bound of $[a,b]$
  with respect to $d_L$ is less than $\ep_2$. Since $d_L(x,a_1) < \ep_2$, $a_1 \in A$. Hence, $b_1 \in L^{-1}(1)$. Inductively,
  $a_i \in A$ and $b_i \in L^{-1}$ for all $i = 1, \dots, n$. Finally, $d_L(b_n,y) \leq d_L(b_n,z) + d_L(z,y) < \ep_1$. Since $L(b_n) = 1$,
  $L(y) > 1 - \ep_1$. Letting $\ep_1 \to 0$ we obtain $\CC_{d_L} f(A) \subset L^{-1}(1)$. Letting $\ep_1 \to \ep$  we obtain
  $ V^{d_L}_{\ep}(\CC_{d_L} f(A)) \subset \{ y : L(y) > 1 - \ep \} = A$.

  (f) The contrapositive of the definition of an elementary Lyapunov function says that if $(x,y) \in f$ with $L(y) < 1$ then $L(x) = 0$.
  It follows that $1 - L$ is an elementary Lyapunov function for $f^{-1}$.

 $\Box$ \vspace{.5cm}

 \begin{prop}\label{prop16af} Let $f$ be a relation on a uniform space
 $(X,\U)$, $\ep > 0$ and $d \in \Gamma(\U)$. Let $K\subset X$ be closed and
 compact.

 (a) For $x \in X$, the set $\{ y : \ell_d^f(x,y) < \ep \} $ is an open subset of $X$ containing $\A_d f(x) \supset \A_{\U} f(x)$.
 It is $\A_d f$ $^+$invariant and so is $\A_{\U} f$ $^+$invariant.

 \begin{align}\label{eq15ua}
 \begin{split}
 \A_{\U} f(K) = &\bigcap_{d \in \Gamma, \ep > 0} \bigcup_{x \in K} \{ y : \ell_d^f(x,y) < \ep \}, \\
K \cup \A_{\U} f(K) = &\bigcap_{d \in \Gamma, \ep > 0} \bigcup_{x \in K} \{ y : \min(\ell_d^f(x,y), d(x,y)) < \ep \}
\end{split}
 \end{align}

 (b) For $x \in X$, the set $\{ y : m_d^f(x,y) < \ep \} $ is an open subset of $X$ containing
 $V^d_{\ep} \circ  \CC_d f\circ V^d_{\ep}(x) \supset V^d_{\ep} \circ  \CC_{\U}f \circ V^d_{\ep}(x)$.
 It is $V^d_{\ep} \circ \CC_d f $ $^+$invariant and so is $V^d_{\ep} \circ \CC_{\U} f $
  and $V^d_{\ep} \circ f $ $^+$invariant. In particular, $\{ (x,y) : m_d^f(x,y) < \ep \} $ is a $\U$ inward set for $f$.

  \begin{align}\label{eq15ub}
   \begin{split}
 \CC_{\U} f(K) = &\bigcap_{d \in \Gamma, \ep > 0} \bigcup_{x \in K} \{ y : m_d^f(x,y) < \ep \}, \\
 K \cup \CC_{\U} f(K) = &\bigcap_{d \in \Gamma, \ep > 0} \bigcup_{x \in K} \{ y : \min(m_d^f(x,y),d(x,y)) < \ep \}
\end{split}
 \end{align}
 \end{prop}

 {\bfseries Proof:}  The sets are open because $\ell^d_f$ and $m^d_f$ are continuous. The set in (a) clearly contains
 $\A_d f(x) = \{ y : \ell_d^f(x,y) = 0 \}$. If $(y,z) \in \A_d f$ then by Proposition \ref{prop00a}
 $\ell_d^f(x,z) \leq \ell_d^f(x,y) + \ell_d^f(y,z) = \ell_d^f(x,y) < \ep.$

 If $y \in V^d_{\ep} \circ   \CC_d f (z)$ with $m_d^f(x,z) < \ep$ then there exists $z_1 \in \CC_d f(z)$ with
 $d(z_1,y) < \ep$. Let $\ep_1 > 0$ and such that $d(z_1,y) + \ep_1, m_d^f(x,z) + 2\ep_1 < \ep$. There exist
 $[a,b] \in f^{\times n}$ and $[c,d] \in f^{\times m}$ such that with respect to $d$ the $xz$
 chain-bound of $[a,b]$ is less than $m_d^f(x,z) + \ep_1$
 and the $zz_1$ chain-bound of $[c,d]$ is less than  $\ep_1$. Notice that $d(b_n,c_1) \leq d(b_n,z) + d(z,c_1) < \ep$
 and $d(c_m,y) \leq d(c_m,z_1) + d(z_1,y) < \ep$. Hence, the $xy$ chain-bound of the concatenation $[a,b] \cdot [c,d]$
 is less than $\ep$. Thus, $\{ y : m_d^f(x,y) < \ep \} $  is $(V^d_{\ep} \circ   \CC_d f)$ $^+$invariant.

Similarly, if $y \in V_{\ep} \circ   \CC_d f (z)$ with $d(x,z) < \ep$ then there exists $z_1 \in \CC_d f(z)$ with
 $d(z_1,y) < \ep$. Let $\ep_1 > 0$ and such that $d(z_1,y) + 2\ep_1, d(x,z) + 2\ep_1 < \ep$. There exists
  $[c,d] \in f^{\times m}$ such that with respect to $d$ the  $zz_1$ chain-bound of $[c,d]$ is less than  $\ep_1$.
   Notice that $d(x,c_1) \leq d(x,z) + d(z,c_1) < \ep$
 and $d(c_m,y) \leq d(c_m,z_1) + d(z_1,y) < \ep$. Hence, the $xy$ chain-bound of the concatenation $[c,d]$
 is less than $\ep$. Thus, $\{ y : m_d^f(x,y) < \ep \} $ contains $V^d_{\ep} \circ  \CC_d f\circ V^d_{\ep}(x)$.

 If $Q: X \times X \to \R$ is a continuous function with $Q \geq 0$, then we let $Q(K,y) = \inf \{ Q(x,y) : x \in K \}$. Clearly,
 $Q(K,y) \leq \ep$ iff there exists $x \in K$ such that $Q(x,y) < \ep$. Also,
 \begin{equation}\label{eq15uab}
 \{x : Q(K,y) = 0 \} \ =  \bigcap_{\ep > 0} \ \{x : Q(K,y) < \ep \}. \hspace{2cm}
 \end{equation}
 Furthermore, if $K$ is compact then $Q(K,y) = 0$ iff
 there exists $x \in K$ such that $Q(x,y) = 0$.

 Recall from (\ref{eq37}) that $\ell_d^{f \cup 1_X}(x,y) = \min(\ell_d^f(x,y),d(x,y))$ and from (\ref{eq38h}) that
 $m_d^{f \cup \{(x,x) \}}(x,y) = \min(m_d^f(x,y),d(x,y))$.

 Let $Q_d(x,y) = m_d^{f \cup \{ (x,x) \}}(x,y)$ so that $Q_d(K,y) = \min(m_d^f(K,y),d(K,y))$. Observe that if $d_1, d_2 \in \Gamma(\U)$ and
 $\ep_1, \ep_2 \geq 0$ then with $d = d_1 + d_2$ and $\ep = \min(\ep_1,\ep_2)$,
 \begin{equation}\label{eq15abb}
 \{(x,y) : Q_d(x,y) \leq \ep \} \subset \{ (x,y) : Q_{d_1}(x,y) \leq \ep_1 \} \cap \{ (x,y) : Q_{d_2}(x,y) \leq \ep_2 \}.
 \end{equation}

 So if $K$ is compact, and $y \in \bigcap_{d \in \Gamma, \ep > 0} \bigcup_{x \in K} \{ y : $ $\min(m_d^f(x,y),d(x,y)) < \ep \}$
  the collection of closed subsets $\{ \{x \in K :  Q_d(x,y) = 0 \} : d \in \Gamma(\U) \}$ satisfies the finite intersection property and
  so has a nonempty intersection. If $x \in K$ is a point of the intersection, then $y \in K \cup \CC_{\U} f(x)$. This proves the
  second equation in (\ref{eq15ub}). The three remaining equations in  (\ref{eq15ua}) and (\ref{eq15ub}) follow from a similar argument
  with $Q_d$ equal to  $\ell_d^{f}$, $\ell_d^{f \cup 1_X}$ and $m_d^{f}$.

 Notice that as functions of $y$ $\ell_d^{f}(x,y)$ and $m_d^{f}(x,y)$ are $d$ Lipschitz
 with Lipschitz constant at most $1$. Hence, for any $K \subset X$,  as functions of $y$, $\ell_d^f(K,y)$ and
 $m_d^f(K,y)$ are $d$ Lipschitz\index{l@$\ell_d^f(K,y)$}\index{m@$m_d^f(K,y)$}
 with Lipschitz constant at most $1$ as are $\min(\ell_d^f(K,y),d(K,y))$ and $\min(m_d^f(K,y),d(K,y))$.

%

 $\Box$ \vspace{.5cm}

 \begin{theo}\label{theo16ah} Let $f$ be a relation on a uniform space $(X,\U)$.
 \begin{enumerate}
 \item[(a)] If $(x,y) \not\in 1_X \cup \CC_{\U} f$, then there exists a $\U$  elementary Lyapunov function $L$ such that
 $L(y) = 0$ and $L(x) = 1$.

 \item[(b)] If $x \not\in |\CC_{\U} f|$,  then there exists a $\U$  elementary Lyapunov function $L$ such that
 $1 > L(x) > 0$.
 \end{enumerate}
 \end{theo}

 {\bfseries Proof:} (a) With $g = f \cup \{ (x,x) \}$, $m^g_d(y) = \min ( m^f_d(x,y), d(x,y))$ by Lemma \ref{lem3b}.
 By hypothesis, there exist $d \in \Gamma$ and $\ep > 0$ so that $m^g_d(x,y) > \ep$. By Proposition \ref{prop16af} (b), the
 set $A = \{ y : m^g_d(x,y) < \ep \}$ is a $\U$ inward set for $g$. By Proposition \ref{theo16ai}  (d) there is a
 $\U$ uniformly continuous elementary Lyapunov function $L$ for $g$ (and hence for $f$) so that $L^{-1}(0) \cup A = X$ and
 $g(A) \subset L^{-1}(1)$.  Since $x \in A$ and $(x,x) \in g$, $x \in g(A)$ and so $L(x) = 1$. Since $y \not\in A$, $L(y) = 0$.

 (b)  By hypothesis, there exist $d \in \Gamma$ and $1 > \ep > 0$ so that $m^f_d(x,x) > 2 \ep$. Let $A_0 = V^d_{\ep}(x)$ and
 $A_1 = \{ y : m^f_d(x,y) < \ep \}$. Since $m^f(x,x) \leq m^f_d(x,y) + d(y,x)$, it follows that $A_0$ and $A_1$ are disjoint.

 By Proposition \ref{prop16af} (b) $V^d_{\ep}(f(A_0 \cup A_1)) \subset A_1$. Let
 $B = f(A_0 \cup A_1)$.  Define $L(y) = \max( [ \ep - d(y,B) ]/\ep, \ep - d(y,x), 0)$. If $(y_1,y_2) \in f$ and $L(y_1) > 0$ then
 $y_1 \in A_0 \cup A_1$ and so $y_2 \in B$. Thus, $L(y_2) = 1$. Thus, $L$ is a $\U$  elementary Lyapunov function.
 Since $x \in A_0$, $d(x,B) > \ep$. Hence, $L(x) = \ep$.

 $\Box$ \vspace{.5cm}

 \begin{df}\label{def16aff}   Let $f$ be a relation on a uniform space $(X,\U)$. We denote by $\L_e$  the set of
 $\U$  elementary Lyapunov functions for  $f$.  We say that a set $\L \subset \L_e$ satisfies the condition
 POIN-E\index{condition!POIN-E} for $\CC_{\U} f$
 if it satisfies POIN for $\CC_{\U} f$ and, in addition,
 \begin{itemize}
 \item If $x \not\in |\CC_{\U} f|$,  then there exists  $L \in \L_e$ such that
 $1 > L(x) > 0$.
 \end{itemize}\end{df}
 \vspace{.5cm}

 By Proposition \ref{theo16ah}, the set $\L_e$ satisfies POIN-E for $\CC_{\U} f$.

 \begin{theo}\label{theo16ak} For  $f$  a relation on a uniform space $(X,\U)$. If $\L \subset \L_e$ satisfies
POIN-E for $\CC_{\U} f$ then
 \begin{equation}\label{eq15uc}
 \begin{split}
 \CC_{\U} f \ = \ \bigcap_{L \in \L} \ \leq_L, \hspace{4cm}\\
 |\CC_{\U} f| \ = \ \bigcap_{L \in \L} \ [L^{-1}(0) \cup L^{-1}(1)] \ = \ \bigcap_{L \in \L} \ |L|_f.
 \end{split}
 \end{equation}
 \end{theo}

 {\bfseries Proof:} The first equation follows from POIN for $\CC_{\U} f$.

 If $L \in \L_e$ then it is an elementary Lyapunov function for $\CC_{\U} f$ by Proposition \ref{theo16ai} (e) and
 $1 - L$ is an elementary Lyapunov function for $\CC_{\U}f^{-1}$ by  Proposition \ref{theo16ai} (f). So with
 $ G_L = \{ x : 1 >  L(x)   > 0 \}$,
 \begin{equation}\label{eq15ud}
   \CC_{\U} f(G_L) \subset L^{-1}(1) \quad \text{and} \quad \CC_{\U} f^{-1}(G_L) \subset L^{-1}(0).
\end{equation}
 Hence, $G_L \cap |\CC_{\U} f| = \emptyset$, i.e.\ $|\CC_{\U} f| \subset |L|_f $.

 On the other hand, if $x \not\in  |\CC_{\U} f|$ then by POIN-E there exists $L \in \L_e$ such that
 $x \in G_L$.

  $\Box$ \vspace{.5cm}

  If $A$ is a $^+$invariant subset for a relation $f$ we denote by $f^{\infty}(A)$ \index{f@$f^{\infty}$} the (possibly empty) maximum
  invariant subset of $A$,\index{subset!maximum invariant subset}\index{maximum invariant subset} i.e.\ the union of all
  $f$ invariant subsets of $A$.
  We can obtain it by a transfinite construction
  \begin{equation}\label{eq15e}
  A_0 = A, \quad A_{\a + 1} = f(A_{\a}), \quad A_{\a} = \bigcap_{\b < \a} A_{\b} \ \ \text{for} \ \a \ \text{a limit ordinal}.
  \end{equation}
  The process stabilizes at $\a$ when $A_{\a + 1} = A_{\a}$ which then equals  $f^{\infty}(A)$.

  \begin{df}\label{def16al} If $A$ is a $\U$ inward set for a relation $f$ then $(\CC_{\U} f)^{\infty}(A)$ is called the
  \emph{$\U$ attractor associated with $A$}. \index{attractor} A $\U$ attractor for $f^{-1}$ is called a \emph{$\U$ repellor} for $f$.
  \index{repellor} If $A$ is a $\U$ inward set for $f$ and $B$ is a  $\U$ inward set for $f^{-1}$ such that
  $A \cup B = X, f(A) \cap B = \emptyset = f^{-1}(B) \cap A$ then the pair $(A_{\infty}, B_{\infty}) =$
  $ ( (\CC_{\U} f)^{\infty}(A), (\CC_{\U} f^{-1})^{\infty}(B) )$
  is called a $\U$ attractor-repellor pair \index{attractor-repellor pair} with $B_{\infty} = (\CC_{\U} f^{-1})^{\infty}(B)$ the repellor dual to
  $A_{\infty} = (\CC_{\U} f)^{\infty}(A)$ and vice-versa.\index{attractor!dual}\index{repellor!dual} \end{df}
   \vspace{.5cm}

   Again, if $\U = \U_M$ we will drop the label $\U$.

 \begin{prop}\label{prop16ahh} Let $f$ be a relation on a  uniform space $(X,\U)$ and let $x,y \in X$ with $y \not\in \ol{\{ x \}}$. The following are equivalent.
 \begin{itemize}
 \item[(i)] $y \in \CC_{\U} f(x)$.
 \item[(ii)] For every $\U$ elementary Lyapunov function $L$ for $f$, $L(x) > 0$ implies $L(y) = 1$.
 \item[(iii)] For every open $\U$ inward set $A$ for $f$, $x \in A$ implies $y \in A$.
 \end{itemize}
 If $x \in |\CC_{\U} f|$, then these conditions are further equivalent to
 \begin{itemize}
 \item[(iv)] For every $\U$ attractor $A_{\infty}$ for $f$, $x \in A_{\infty}$ implies $y \in A_{\infty}$.
 \end{itemize}
 \end{prop}

 {\bfseries Proof:} (i) $\Rightarrow$ (ii): A $\U$ elementary Lyapunov function for $f$ is a $\U$ elementary Lyapunov function for $\CC_{\U} f$ by
 Theorem \ref{theo16ai}(e).

 (i) $\Rightarrow$ (iii): A $\U$ inward set for $f$ is  $\CC_{\U} f$ $^+$invariant by
 Theorem \ref{theo16ai}(b).

 (ii) $\Rightarrow$ (i): Apply Theorem \ref{theo16ah} (a).

 (iii) $\Rightarrow$ (i): By Proposition \ref{prop16af} (b), with $g = f \cup \{(x,x) \}$,
  $\{ y : m_d^g(x,y) < \ep \} = \{ y : \min(m_d^f(x,y), d(x,y)) < \ep \}$ is a $\U$ inward set for $g$ and hence for $f$. So
  (\ref{eq15ub}) implies that $\ol{\{ x \}} \cup \CC_{\U} f(x)$ is the intersection of $\U$ inward sets.

 If $x \in |\CC_{\U} f|$, then $\CC_{\U} f(x)$ is $\CC_{\U} f$ invariant and so $x$ is contained in an inward set $A$ iff it is contained in the
 associated attractor. Hence (iii) $\Leftrightarrow$ (iv) in this case.

 Notice that if $x \in |\CC_{\U} f|$ then $\ol{\{ x \}}$ is contained in the closed set $\CC_{\U} f(x)$.

 $\Box$ \vspace{.5cm}

   \begin{prop}\label{prop16am} If $A_{\infty}$ is the $\U$ attractor associated with the $\U$ inward set $A$, then
   $A \cap |\CC_{\U}f| \subset A_{\infty}$. Furthermore,
 \begin{equation}\label{eq15udd}
 |\CC_{\U} f| \ = \ \bigcap \{ A_{\infty} \cup B_{\infty} : \ (A,B) \ \text{a} \ \U \ \text{attractor-repellor pair for} \ f \}.
  \end{equation}

  If $x \in |\CC_{\U} f|$ then the $\CC_{\U} f \cap \CC_{\U} f^{-1}$ equivalence class of $x$ in $|\CC_{\U} f|$ is given by
  \begin{equation}\label{eq15ude}
  (\CC_{\U} f \cap \CC_{\U} f^{-1})(x) \ = \
  \bigcap \ \{ B :   B \ \text{ a } \ \U \ \text{attractor or repellor with } \ x \in B \}.
  \end{equation}
 \end{prop}

 {\bfseries Proof:} For any $\CC_{\U} f$ $^+$invariant set $A$, if $x \in |\CC_{\U} f|$ then $\CC_{\U} f(x)$ is a
 $\CC_{\U} f$ invariant subset of $A$ and so is contained $(\CC_{\U} f)^{\infty}(A)$. So if $(A,B)$ is an attractor-repellor pair
 then $|\CC_{\U} f| = |\CC_{\U} f| \cap (A \cup B) \subset |\CC_{\U} f| \cap (A_{\infty} \cup B_{\infty})$.

 In particular, if $L$ is a
 $\U$ elementary Lyapunov function then with $A = \{ x : L(x) > 0 \}$ and $B = \{ x : L(x) < 1 \}$, the associated
 attractor-repellor pair $(A_{\infty},B_{\infty})$ satisfies $A_{\infty} \subset L^{-1}(1)$, $B_{\infty} \subset L^{-1}(0)$ and
 so $|\CC_{\U} f| \cap L^{-1}(1) = |\CC_{\U} f| \cap A_{\infty}$ and $|\CC_{\U} f| \cap L^{-1}(0) = |\CC_{\U} f|
\cap B_{\infty}$. Hence, (\ref{eq15udd}) follows from (\ref{eq15uc}).

Finally, $(\CC_{\U} f \cap \CC_{\U} f^{-1})(x) \ = \ \CC_{\U} f(x) \cap \CC_{\U} f^{-1}(x)$. By Proposition \ref{prop16ahh} $\CC_{\U} f(x)$
is the intersection of the attractors containing $x$ and $\CC_{\U} f^{-1}(x)$ is the intersection of the repellors containing $x$.

  $\Box$

\newpage

\section{ Upper-semicontinuous Relations and Compactifications}\label{cusc}
\vspace{.5cm}

Up to now we have generally imposed no topological conditions on the relation $f$. Consider $f : X \to Y$ a relation with
$X$ and $Y$ Tychonoff spaces, i.e.\ $f \subset X \times Y$. Call $f$ a \emph{closed relation} when it is a closed
subset of $X \times Y$.  Call $f$ \emph{pointwise closed}\index{relation!pointwise closed} when $f(x)$ is closed for every $x \in X$. Call $f$
\emph{pointwise compact}\index{relation!pointwise compact} when $f(x)$ is compact for every $x \in X$. Since $f(x)$ is the pre-image of $f$
by the continuous map $y \mapsto (x,y)$ it follows that a closed relation is pointwise closed. Since $Y$
is Hausdorff a pointwise compact relation is pointwise closed.

If $f : X \to Y$ is a relation and $B \subset Y$, recall that  $f^*(B) = \{ x \in X : f(x) \subset B \}$.\index{f@$f^*(B)$}
For example, $f^*(\emptyset) = \{ x : f(x) = \emptyset \}$ which is the complement of the domain of $f$,
$Dom(f) = f^{-1}(X)$.

We will need the properties of proper maps.  These are reviewed in Appendix C.

\begin{theo}\label{theo16aa} Let $f : X \to Y$ be a relation between Tychonoff spaces.

\begin{enumerate}
\item[(a)] If $f$ is a closed relation and $A \subset X$ is compact, then $f(A) \subset Y$ is closed.

\item[(b)] The following conditions are equivalent. When they hold we call $f$ an \emph{upper semi-continuous relation},\index{relation!usc}\index{usc relation}
written $f$ is usc.

\begin{itemize}

\item[(i)] If $B$ is a closed subset of $Y$, then $f^{-1}(B)$ is a closed subset of $X$.

\item[(ii)] If $B$ is an open subset of $Y$ then $f^*(B)$ is an open subset of $X$.

\item[(iii)] If $\{ x_i : i \in I \}$ is a net in $X$ converging to $x \in X$ and $B$ is an open set
containing $f(x)$ then eventually $f(x_i) \subset B$.
\end{itemize}

\item[(c)] A usc relation is closed iff it is pointwise closed.

\item[(d)] If $f$ and $f^{-1}$ are usc, then $f$ and $f^{-1}$ are closed relations.

\item[(e)] Let $\pi_1 : X \times Y \to X$ be the projection map.
If the restriction $\pi_1|f : f \to X$ is a closed map, then $f$ is usc.

\item[(f)] The following conditions are equivalent. When they hold we call
$f$ a \emph{compactly upper semi-continuous relation},
written $f$ is cusc.\index{relation!cusc}\index{cusc relation}

\begin{itemize}
\item[(i)] With $\pi_1 : X \times Y \to X$ the projection map,
the restriction $\pi_1|f : f \to X$ is a proper map.

\item[(ii)] The relation $f$ is pointwise compact and usc.
\end{itemize}

\item[(g)] If $f$ is  cusc  then $f$ is a closed relation and $A $  a compact subset of $X$,
implies that $f(A)$ is a compact subset of $Y$.

\item[(h)] If $X$ is a k-space, $f$ is a closed relation and for every compact subset
$A$ of $X$, the subset $f(A)$ of $Y$ is compact, then
$f$ is  cusc.

\item[(i)] If $f$ is cusc and $g \subset f$ then $g$ is cusc iff $g$ is closed.
\end{enumerate}
\end{theo}

{\bfseries Proof:}
(a) Since $A$ is compact, the trivial  map of $A$ to a point is proper. Hence, $\pi_2 : A \times Y \to Y$ is a closed map.
If $f$ is closed then $\pi_2((A \times Y)\cap f) = f(A)$ is closed.


(b) (i) $\Leftrightarrow$ (ii): $f^*(B) = X \setminus f^{-1}(Y \setminus B)$.

(ii) $\Leftrightarrow$ (iii): If $f^*(B)$ is open then eventually $x_i \in f^*(B)$.
If $f^*B$ is not open then there is a net $\{ x_i \}$ in the complement which converges to a point $x \in f^*(B)$.
Then $f(x) \subset B$ but never $f(x_i) \subset B$, contradicting (iii).

(c) Assume $f$ is usc and pointwise closed. Suppose $\{ (x_i,y_i) \}$ is a
net in $f$ converging to $(x,y)$ but with $(x,y) \not\in f$ and so $y \not\in f(x)$.
Since $f(x)$ is closed and $Y$ is Tychonoff, there is are disjoint open sets $B, G$ with $f(x) \subset B$ and $y \in G$.
Since $f$ is usc, eventually $f(x_i) \subset B$. In particular, eventually $y_i \in B$ and so eventually $y_i \not\in G$.
This contradicts convergence of $\{ y_i \}$ to $y$.

We saw above that a closed relation is always pointwise closed.

(d) If $f^{-1}$ is usc then $f(x) = (f^{-1})^{-1}(x)$ is closed.  Since $f$ is usc, it is closed by (c). Hence, $f^{-1}$ is
closed as well.

(e) If $B$ is a closed subset of $Y$, then $f \cap (X \times B)$ is a closed subset of $f$.  If $\pi_1|f$ is a closed map
then $f^{-1}(B) = \pi_1(f \cap (X \times B))$ is closed.

(f) (i) $\Rightarrow$ (ii): A proper map is closed and so $f$ is usc by (e).
Since $\pi_1|f$ is proper, $(\pi_1|f)^{-1}(x) = \{ x \} \times f(x)$ is compact by  Proposition
\ref{appprop04}(a). Hence, $f$ is pointwise compact.

(ii) $\Rightarrow$ (i): We verify condition (iv) of  Proposition
\ref{appprop04}(a). Let $\{ (x_i,y_i) \}$ be a net in $f$ such that $\{ x_i \}$ converges to $x \in X$. If $B$ is any
open set containing $f(x)$ then eventually $f(x_i) \subset B$ because $f$ is usc. So eventually $y_i \in B$. Because $f(x)$ is compact,
Lemma \ref{applem02} implies that $f(x)$ contains a cluster point of $ \{ y_i \}$. That is, there is a subnet $\{ y_{i'} \}$ which
converges to a point $y \in f(x)$.  Hence $\{ (x_{i'},y_{i'}) \}$ converges to $(x,y) \in (\pi_1|f)^{-1}(x)$.

(g) A pointwise compact relation is pointwise closed and so a cusc relation is a closed relation by (b).
If $A \subset X$ is compact then $(\pi_1|f)^{-1}(A)$ is compact by Proposition
\ref{appprop04} (c). Hence, $f(A) = \pi_2[(\pi_1|f)^{-1}(A)]$ is compact, where
$\pi_2 : X \times Y \to Y$ is the other projection.

 (h) If $A$ and $f(A)$ are compact and $f$ is closed then $(A \times f(A)) \cap f = (\pi_1|f)^{-1}(A)$ is compact.
 So the result follows from Proposition \ref{appprop05} (a).

 (i) If $g$ is cusc then it is closed by (e) and (c).  If $\pi|f : f \to X$ is a proper map and $g$ is a closed subset
 of $f$ then $\pi|g$ is proper by Proposition \ref{appprop03} (d).

 $\Box$ \vspace{.5cm}

 {\bfseries Remark:} The condition that a pointwise compact relation be usc, and so cusc, is weaker than the demand that
 $x \mapsto f(x)$ is continuous as a function from $X$ to the space of compact subsets with the Hausdorff topology.
 For a comparison in the compact case, see \cite{A93} Chapter 7.
 \vspace{.5cm}

 We call $f$ a \emph{proper relation}\index{relation!proper}\index{proper relation} when both $f$ and $f^{-1}$ are cusc relations, or, equivalently
 when $\pi_1|f : f \to X$ and $\pi_2|f : f \to Y$ are both proper maps.

 \begin{prop}\label{prop16ab} Let $f : X \to Y$ be a map between Tychonoff spaces.
the following are equivalent:
 \begin{itemize}
 \item[(i)] $f$ is a continuous map.
 \item[(ii)] $f$ is a usc relation.
\item[(iii)] $f$ is a cusc relation.
\end{itemize}

If $f$ is continuous then
 $f$ is a closed map iff $f^{-1}$ is a usc relation, and the following are equivalent
 \begin{itemize}
 \item[(iv)] $f$ is a proper map.
 \item[(v)] $f^{-1}$ is a cusc relation.
 \item[(vi)] $f$ is a proper relation.
\item[(vii)] $f$ is a closed map and $f^{-1}(y)$ is compact for every $y \in Y$.
\end{itemize}
\end{prop}

{\bfseries Proof:} (i) $\Leftrightarrow$ (ii): Both say that $f^{-1}(B)$ is closed when $B$ is.

(ii) $\Leftrightarrow$ (iii): because $f$ is pointwise compact.

The relation $f^{-1}$ is usc iff $f(A)$ is closed when $A$ is.

(iv) $\Leftrightarrow$ (vii): by Proposition \ref{appprop04}.

(v) $\Leftrightarrow$ (vi): Since $f$ is a continuous map it is a cusc relation so it is a proper relation iff
$f^{-1}$ is a cusc relation.

(v) $\Leftrightarrow$ (vii): Condition (vii) says that $f^{-1}$ is usc and pointwise compact.

 $\Box$ \vspace{.5cm}

 \begin{theo}\label{theo16aaa} Let $f : X \to Y$ and $g : Y \to Z$ be  relations between Tychonoff spaces.
\begin{enumerate}
\item[(a)] If $f$ and $g$ are usc then $g \circ f$ is usc.
\item[(b)] If $f, g$ and $g^{-1}$ are usc and closed, then $g \circ f$ is usc and closed.
\item[(c)] If $f$ and $g$ are cusc  then $g \circ f$ is cusc.
\item[(d)] If $f$ is cusc and $g$ is closed then $g \circ f$ is closed.
 \end{enumerate}
 \end{theo}

 {\bfseries Proof:} (a) If $C \subset Z$ is closed then $(g \circ f)^{-1}(C) = f^{-1}(g^{-1}(C))$ is closed.

 (b) By (a) $g \circ f$ is usc. For $x \in X$, $g \circ f(x) = g(f(x))$ is closed since $f$ is pointwise closed
 and $g^{-1}$ is usc. Hence, $g \circ f$ is pointwise closed, and so is closed by \ref{theo16aa} (c).

 (c) By Theorem \ref{theo16aa}(f) $g \circ f(x) = g(f(x))$ is compact since $f$ is pointwise compact and $g$ is
 cusc.

 (d) Since $f$ is cusc, $\pi_{13} : f \times Z \to X \times Z$ is a closed map. Since $g$ is a closed relation,
 $(f \times Z) \cap (X \times g)$ is a closed subset and so its image $g \circ f \subset X \times Z$ is closed.

 $\Box$ \vspace{.5cm}

  \begin{prop}\label{prop16ac} Let $f, g : X \to Y$ be  relations between Tychonoff spaces.
\begin{enumerate}
\item[(a)] If $f$ and $g$ are both closed, usc or cusc then $g \cup f$ satisfies the corresponding property.
\item[(b)] If $f$ is cusc and $g$ is closed, then $g \cap f$ is cusc.
\item[(c)] Assume $Y$ is a normal space. If $f$ and $g$ are both closed and usc  then $g \cap f$ is closed and usc.
 \end{enumerate}
 \end{prop}

  {\bfseries Proof:} (a) For $B \subset Y$, $(f \cup g)^{-1}(B) = (f^{-1} \cup g^{-1})(B) = f^{-1}(B) \cup g^{-1}(B)$.
  Since the union of two closed sets is closed it follows that $g \cup f$ is
  closed or usc when each of $g$ and $f$ is closed or
  usc.  Furthermore, $(f \cup g)(x) = f(x) \cup g(x)$ and so $f \cup g$ is pointwise compact when $f$ and $g$ are.

  (b) Apply Theorem \ref{theo16aa}(i).

  (c) If $U$ is an open set containing $(g \cap f)(x) = g(x) \cap f(x)$ then since $f$ and $g$ are closed
  $g(x) \setminus U$ and $f(x) \setminus U$ are disjoint closed sets. Since $Y$ is normal we can choose
  disjoint open sets $V_1 \supset g(x) \setminus U$ and $ V_2 \supset f(x) \setminus U$. Hence,
  $U_1 = V_1 \cup U \supset g(x)$ and $U_2 = V_2 \cup U \supset f(x)$ with $U_1 \cap U_2 = U$. Since $g$ and $f$ are
  usc, $g^*(U_1) \cap f^*(U_2)$ is an open set containing $x$ and contained in $(g \cap f)^*(U)$. Thus,
  $(g \cap f)^*(U)$ is a neighborhood of $x$.  Hence, $g \cap f$ is usc.

 $\Box$ \vspace{.5cm}

 \begin{ex} For $f$ a relation on $X$ with $\pi_1|f$ the first coordinate projection, $f$ can be usc without
 $\pi_1|f$ being closed. Furthermore, with $g \subset f$  closed, $g$ need not be usc. \end{ex}

 {\bfseries Proof:} Let $X = \R$ and $f = f^{-1} = \{ (t,1/t) : t \not= 0  \in \R \} \cup \{ (t,0), (0,t): t \in \R \}$.
 Let $g = g^{-1} = \{ (t,1/t) : t \not= 0  \in \R \} \cup \{ (0,0) \}$.

  $\Box$ \vspace{.5cm}

  Now we illustrate how these conditions on a relation may be applied.

 \begin{lem}\label{lem16ad} Let $F$ be a closed, reflexive, transitive relation on a normal Hausdorff
 space $X$ with $F$ and $F^{-1}$ usc.
 If $A$ is a closed, $F$ invariant set and $U$ is an open set with $A \subset U$ then
  there exists a closed, $F$ invariant set $B$ such that
 $A \subset B^{\circ}$ and $B \subset U$.\end{lem}

 {\bfseries Proof:} Because $F^{-1}$ is usc, $F(A)$ is closed.
 Since $F$ is usc, $F^*(U)$ is open and since $A$ is $F$ invariant,
$A \subset F^*(U)$. Use normality  to choose a closed set $B_1$ so that
 $A \subset B_1^{\circ}$ and $B_1 \subset F^*(U)$. The set $B = F(B_1) \subset U$ is closed because $F^{-1}$ is usc and
$A \subset B_1^{\circ} \subset B^{\circ}$ because $F$ is reflexive.

 $\Box$ \vspace{.5cm}

 The following is a version of \cite{N} Theorem 2, see also \cite{AA10a} and \cite{AA10b}.

 \begin{theo}\label{theo16ae} Let $F$ be a closed, transitive relation on a normal
 Hausdorff space $X$ with $F$ and $F^{-1}$
 usc. Assume that $X_0$ is a closed subset of $X$ and $L_0 : X_0 \to [a,b]$ is a bounded, Lyapunov function for the
 restriction $F_0 = F \cap (X_0 \times X_0)$.  There exists
 $L : X \to [a,b]$ a  Lyapunov function for $F$ such that $L(x) = L_0(x)$ for $x \in X_0$.
 \end{theo}

 {\bfseries Proof:}
 Replacing $F$ by $F \cup 1_X$,
 we can assume that $F$ is reflexive as well as transitive. Without loss of generality we can assume that
 $[a,b] = [0,1]$.

 We mimic the proof of Urysohn's Lemma. Let $\Lambda = \mathbb{Q}
\cap [0,1]$ counted with $\lambda_0 = 0, \lambda_1 = 1$. Let $B_0 = X, B_1 = \emptyset$.
For all $\lm \in \Lambda$ we define the closed set $B_{\lm} \subset X$
 so that:
 \begin{itemize}
 \item[(a)] $F(B_{\lm}) = B_{\lm}$, i.e.\ $B_{\lm}$ is $F$ invariant.
 \item[(b)] $L_0^{-1}((\lm,1]) \subset B_{\lm}^{\circ}$.
 \item[(c)] $L_0^{-1}([0,\lm)) \cap B_{\lm} = \emptyset$.
 \item[(d)] If $\lm' < \lm \in \Lambda$, then $B_{\lm} \subset B_{\lm'}^{\circ}$.
 \end{itemize}

 Observe that if $x$ were a point of $F(L_0^{-1}([\lm,1]) \cap F^{-1}(L_0^{-1}([0,\lm))$, then there would exist
 $z_1,z_2 \in X_0$ with $L_0(z_1) < \lm \leq L_0(z_2)$ and $(z_2,x),(x,z_1) \in F$ and so $(z_2,z_1) \in F_0$ which
 would contradict the assumption that $L_0$ is a Lyapunov function for $F_0$.

We repeatedly apply Lemma \ref{lem16ad}.  We will use the notation $A \subset \subset B$ to mean
$\ol{A} \subset B^{\circ}$. A space is normal exactly when $A \subset \subset B$ implies there
exists $C$ such that $A \subset \subset C \subset \subset B$. Lemma \ref{lem16ad} says that
if $A$ is closed and $F$ invariant and $A \subset \subset B$ then there exists $C$
closed and $F$ invariant such that $A \subset \subset C \subset \subset B$.

Proceed inductively assuming that $B_{\lambda}$ has been defined
for all $\lambda$ in $\Lambda_n = \{ \lambda_i : i = 0,...,n \}$
with $n \geq 1$. Let $\lambda = \lambda_{n+1}$ and let $\lambda' <
\lambda < \lambda''$ the  nearest points in $\Lambda_{n}$
below and above $\lambda$.

Choose a sequence $\{ t_{n}^{-} \}$ with $t_{0}^{-} = \lambda'$,
increasing with limit $\lambda$ and $ \{ t_{n}^{+} \}$ with
$t_{0}^{+} = \lambda''$, decreasing with limit $\lambda$.

Define $Q_{0}^{-} = B_{\lambda'}$ and $Q_{0}^{+} = B_{\lambda''}$.
Inductively, apply Lemma \ref{lem16ad} to choose
$Q_{n}^{+}$ and then $Q_{n}^{-}$ for $n = 1,2,...$ so that
$F(Q_{n}^{\pm}) = Q_{n}^{\pm}$ and

\begin{equation}\label{eqNach16}
\begin{split}
F(L_{0}^{-1}([t_{n}^{+},1]) \cup Q_{n-1}^{+} \quad \subset
\subset \quad Q_{n}^{+} \quad \subset \subset \quad Q_{n-1}^{-}
\setminus
(F)^{-1}(L_{0}^{-1}([0,\lambda]), \\
F(L_{0}^{-1}([\lambda,1]) \cup Q_{n}^{+} \quad \subset \subset
\quad Q_{n}^{-} \quad \subset \subset \quad Q_{n-1}^{-} \setminus
(F)^{-1}(L_{0}^{-1}([0,t_{n}^{-}]).
\end{split}
\end{equation}
Finally, define

\begin{equation}\label{eqNach17}
B_{\lambda} \quad = \quad \bigcap_n Q_{n}^{-}, \hspace{2cm}
\end{equation}
so that

\begin{equation}\label{eqNach18}
 B_{\lambda} \quad \supset  \quad
\bigcup_n Q_{n}^{+}.\hspace{2cm}
\end{equation}

It is easy to check that $B_{\lambda}$ satisfies the required conditions, thus extending the definitions
to $\Lambda_{n+1}$. By induction they can be defined on the entire set $\Lambda$.

Having defined the $B_{\lambda}$'s we proceed as in Urysohn's Lemma to
define $L(x)$ by the Dedekind cut associated with $x$.  That is,

\begin{equation}\label{eqNach19}
L(x) \quad = \quad inf \{ \lambda : x \not\in B_{\lambda} \} \quad
= \quad sup \{ \lambda : x \in B_{\lambda} \}.
\end{equation}
Continuity follows as in  Urysohn's Lemma.
Because each $B_{\lambda}$ is $F  $ invariant, $L$ is a Lyapunov
function.  The additional conditions on these sets imply that if
$x \in X_0$ then $x \in B_{\lambda}$ if $\lambda < L_0(x)$ and
$x \not\in B_{\lambda}$ if $\lambda > L_0(x)$.
Hence, $L$ is an extension of $L_0$.

$\Box$ \vspace{.5cm}

Fathi and Pageault use a slightly different, asymmetric definition of the barrier functions which yields equivalent
results when $f$ is usc.
\begin{equation}\label{eq11}
\begin{split}
L_d^f(x,y) \ = \ \inf  \ \{ d(x,a_1) + \Sigma_{i = 1}^{n-1} d(b_i,a_{i+1}) + d(b_n,y) : \hspace{1cm}\\
[a,b] \in f^{\times n} \ \mbox{with} \ a_1 = x, n = 1,2,... \}, \hspace{2cm} \\
M_d^f(x,y) \ = \ \inf  \ \{ \max(d(x,a_1) , d(b_1,a_2), \dots,  d(b_{n-1},a_n), d(b_n,y)) :\\
[a,b] \in f^{\times n} \ \mbox{with} \ a_1 = x, n = 1,2,... \}, \hspace{2cm}
\end{split}
\end{equation}
\index{l@$L_d^f$}\index{m@$M_d^f$}
So, of course, the first term, $d(x,a_1) = 0$. For the case where $x$ is not in the domain of $f$ we use the convention
\begin{equation}\label{eq12}
f(x) = \emptyset \qquad \Longrightarrow \qquad M_d^f(x,y) \ = \ diam(X),  L_d^f(x,y) \ = \ 2 diam(X).
\end{equation}

We have
\begin{equation}\label{eq13}
\ell^f_d \  \leq \ L^f_d \qquad \text{and} \qquad m^f_d \  \leq \ M^f_d,
\end{equation}
 because
for $L^f_d$ and $M^f_d$ the infimum is taken over a smaller set.

\begin{prop}\label{prop00uu} Let $f$ be a usc relation on a Hausdorff uniform space $(X,\U)$.
For every $x \in X, d \in \Gamma(\U)$ and $\ep > 0$, there exist $d_1 \in \Gamma(\U)$ and $\d > 0$ such that for all $y \in X$
\begin{equation} \label{eq14}
\begin{split}
\ell^f_{d_1}(x,y) \ < \ \d \qquad \Longrightarrow \qquad L^f_d(x,y) \ < \ \ep, \\
m^f_{d_1}(x,y) \ < \ \d \qquad \Longrightarrow \qquad M^f_d(x,y) \ < \ \ep.
\end{split}
\end{equation}
If $\U = \U(d)$ for a metric $d$ then we can choose $d_1 = d$.
\end{prop}

{\bfseries Proof:}
Because $f$ is usc, there exists $d_0 \in \Gamma$ and $\ep/2 > \d > 0$ so that
$f(V^{d_0}_{\d}(x)) \subset V^d_{\ep/2}(f(x))$.
Let $d_1 =  d_0 + d$. If the metric $d$ determines the topology on $X$ then we can use $d_0 = d$ and use $d_1 =  d$.

Now assume $\ell^f_{d_1}(x,y) < \d$.  We need only consider
sequences $[a,b] \in f^{\times n}$ with $xy$ chain-length with respect to
$d_1$ less than
$\d$. With $m^f_{d_1}(x,y) < \d$, consider sequences $[a,b] \in f^{\times n}$ with $xy$ chain-bound less than
$\d$. In either case, $d_1(x,a_1) < \d$ and so $d_0(x,a_1) < \d$. Hence, $f(a_1) \subset V^d_{\ep/2}( f(x))$
 and we can choose $\bar{b}_1 \in f(x)$
such that $d(\bar{b}_1,b_1) < \ep/2$. Replacing the initial pair $(a_1,b_1)$ in $[a,b]$ by $(x,\bar{b}_1)$ we obtain
a sequence with initial point $x$ and whose chain-length is at most $\ep/2 $
 plus the $xy$ chain-length of $[a,b]$ with respect to $d$
because
$d(\bar{b}_1, a_2) \leq d(\bar{b}_1, b_1) + d(b_1,a_2)$, or, if $n = 1$, the
same inequality is used with $y$ replacing $a_2$.
The $xy$ chain-length of $[a,b]$ with respect to $d$ is at most  the  $xy$ chain-length of $[a,b]$ with respect to $d_1$
and so at most $\d < \ep/2$. So the revised sequence which begins with $x$ has $xy$
chain-length with respect to $d$  less than
 $\ep$. Hence, $L_f^d(x,y) < \ep$.

Similarly, the new chain-bound with respect to $d$ is less than $\ep/2 $ plus
 the $xy$ chain-bound of $[a,b]$ with respect to $d_1$.

Notice in passing that if $f(x) = \emptyset$ then the chosen $\d$ implies
$f(a) = \emptyset$ for all $a \in V^{d_0}_{2\d}(x)$ with
$d(x,a) < \d$. Provided that $\d$ has been chosen less than the $d$ diameter of $X$, then from the convention when $f(x) = \emptyset$
it easily follows that then $\ell^f_d(x,y), m^f_d(x,y) \geq \d$ for all $y \in X$ and so the result holds vacuously.

 $\Box$ \vspace{.5cm}

One advantage of the asymmetric definition $M_d^f$ is that, as Pageault points
out in \cite{P}, we can sharpen (\ref{eq9a}) to get
\begin{equation}\label{eq9aa}
M_d^f(x,y) \ \leq \ \max(M_d^f(x,z), M_d^f(z,y)) \quad \text{for all} \ z \in X.
\end{equation}

From (\ref{eq13}),  Proposition \ref{prop00uu} and Theorem \ref{theo11u} the following is obvious.

 \begin{cor}\label{cor00ua} If $f$ is a usc relation on a Hausdorff uniform space $(X,\U)$,
 then $\A_{\U} f = \{ (x,y) : L_f^d(x,y) = 0$ for all $d \in \Gamma(\U) \}$ and
 $\CC_{\U} f = \{ (x,y) : M_f^d(x,y) = 0$ for all $d \in \Gamma(\U)\}$.

 If $d$ is a metric on $X$ with $\U = \U(d)$ then $\A_{d} f = \{ (x,y) : L_f^d(x,y) = 0 \}$
 and $\CC_{\U} f = \CC_{d} f = \{ (x,y) : M_f^d(x,y) = 0 \}$.  \end{cor}

 $\Box$ \vspace{.5cm}

  \begin{prop}\label{prop16aj}  Let $f$ be a relation on a Hausdorff uniform space $(X,\U)$.

 (a) If $f$ is a cusc relation,  then
 \begin{equation}\label{16ua}
 \begin{split}
 \G f \ = \ f \cup (\G f) \circ f, \hspace{2cm} \\
  \A_{\U} f \ = \ f \cup (\A_{\U} f) \circ f, \hspace{2cm}\\
   \CC_{\U} f \ = \ f \cup (\CC_{\U} f) \circ f,\hspace{2cm}
  \end{split}
  \end{equation}
 and if $d \in \Gamma(\U)$ is a metric whose topology is that of $X$ then

  \begin{equation}\label{16ub}
 \A_{d} f \ = \ f \cup (\A_{d} f) \circ f, \quad \text{and} \quad \CC_{d} f \ = \ f \cup (\CC_{d} f) \circ f.
 \end{equation}

 (b) If $f^{-1}$ is a cusc relation,  then
 \begin{equation}\label{16uc}
 \begin{split}
  \G f \ = \  f \cup f \circ (\G f), \hspace{2cm} \\
 \A_{\U} f \ = \ f \cup f \circ (\A_{\U} f) , \hspace{2cm} \\
  \CC_{\U} f \ = \ f \cup f \circ (\CC_{\U} f),\hspace{2cm}
 \end{split}
 \end{equation}
 and if $d \in \Gamma(\U)$ is a metric whose topology is that of $X$ then
  \begin{equation}\label{16ud}
 \A_{d} f \ = \ f \cup f \circ (\A_{d} f), \quad \text{and} \quad \CC_{d} f \ = \ f \cup f \circ (\CC_{d} f).
 \end{equation}
 \end{prop}

 {\bfseries Proof:} In general, if $f \subset F$ and $F$ is transitive, then
 $f \cup F \circ f, f \cup f \circ F \subset F$.
 Furthermore, each of these relations is transitive:
 $$ (f \cup F \circ f) \circ (f \cup F \circ f) \subset F \circ (1_X \cup F) \circ f \subset F \circ f.$$
  Similarly, for $f \cup f \circ F$.  Since each of $F = \G f, \A_{\U} f $ and $ \CC_{\U} f$
  is a transitive relation containing
  $f$, it suffices to prove the reverse inclusions.

 (a) If $f$ is cusc,  then by  Theorem \ref{theo16aaa} (d) $f \cup \G f \circ f$ is a closed,
 transitive relation which contains $f$
 and so contains $\G f$.

Suppose $(x,y) \in \A_{\U} f$. For every $\a = (d, \ep) \in \Gamma \times \R_+$ there is
 $[a,b]_{\a} \in f^{\times n_{\a}} $ whose $xy$ chain-length with respect to $d$ is less than $\ep$.
 Since $d(x,(a_1)_{\a}) < \ep$ and $d(y,(b_{n_{\a}})_{\a}) < \ep$ it follows that $\{ (a_1)_{\a} \} \to x$
 and $\{ (b_{n_{\a}})_{\a}\} \to y$. Since $((a_1)_{\a},(b_1)_{\a}) \in f$ and
 $\pi_1|f : f \to X$ is proper,  Proposition \ref{appprop04}(iv) implies there is a subnet $ \{ (b_1)_{\a'} \}$ converging
 to a point $z$ with $(x,z) \in f$. Now if  $n_{\a'} = 1$ frequently then $z = y$ and $(x,y) \in f$.
 Otherwise we may assume
 all $n_{\a'} > 1$ and define $[a,b]'_{\a} \in f^{\times (n_{\a'} - 1)} $ by omitting the first pair.

 Now given $d \in \Gamma$ and $\ep > 0$ there exists $\a'_1 = (d_1,\ep_1)$ so that $\a'_1  \prec \a'$ implies
 $d((b_1)_{\a'},z) < \ep/2$.
 If $\a'_1 \prec \a' = (\bar d, \bar \ep)$ with $d \leq \bar d$ and $\ep/2 \geq \bar \ep$ then the
 $zy$ chain-length of $[a,b]'_{\a'}$ with respect to $d$ is bounded by $d((b_1)_{\a'},z) $ ($ < \ep/2$) plus the
 $xy$ chain-length of $[a,b]_{\a'}$ with respect to $\bar d$ ($ < \bar d < \ep/2$). It follows that
 $\ell_d^f(z,y) = 0$ for all $d \in \G$.  That is, $(x,z) \in f$ and $(z,y) \in \A_{\U} f$.

 The proof for $\CC_{\U} f$ uses the same argument with chain-bound replacing chain-length throughout.

 If $d$ is a metric in $\Gamma$ with the topology that of $X$, then we keep the metric fixed in the arguments above to
 prove the results for $\A_{d} f$ and $\CC_{d} f$.

 (b) We apply the results of (a) to $f^{-1}$ and invert both sides of the equation using
 $(f \circ g)^{-1} = g^{-1} \circ f^{-1}$
 and  $(\G f)^{-1} = \G (f^{-1}), (\A_{d} f)^{-1} = \A_{d} (f^{-1})$ and the similar equation for $\CC_{d}$.

 $\Box$ \vspace{.5cm}

 \begin{prop}\label{prop16ak}  Let $f \subset F$ be  relations on a set $X$ with $F$ transitive.

 (a) If $A$ an $F$ $^+$invariant subset of $X$, then $A$ is $f$ $^+$invariant. If, in addition, $F = f \cup f \circ F$, then
 $f(A) = F(A)$ for any subset $A$ of $X$. In particular, $A$ is $f$ invariant iff it is $F$ invariant.

 (b) If $F = f \cup F \circ f$ then $Dom(f) = Dom(F)$ (Recall that $Dom(f) = f^{-1}(X)$).

  (c)  Assume that $F = f \cup F \circ f$ and that $F_1$ is also a transitive relation on $X$ with   $F_1 = f \cup F_1 \circ f$.
  If $1_X \cup F_1 = 1_X \cup F$,
  then $F = F_1$.

 (d) Assume that $f \cup F \circ f = F = f \cup f \circ F$. If $L$ is a Lyapunov function for $F$, then
 $x$ is a regular point for $f$ iff it is a regular point for $F$, i.e.\  $|L|_f = |L|_F$.

 (e)  Assume that $F = f \cup F \circ f$.
 If $f$ is a mapping then $f$ maps $F$ to itself, $F^{-1}$ to itself and $F \cap F^{-1}$ to itself. Hence,
 $f(|F|) \subset |F|$. If, in addition, $F = f \cup f \circ F$, then $f(|F|) = |F|$ and if $E$ is any
 $F \cap F^{-1}$ equivalence class in $|F|$, then $f(E) = E$.
 \end{prop}

 {\bfseries Proof:} (a) $A$ is $f$ $^+$invariant because $f \subset F$. Also, $f(A) \subset F(A)$.
 Conversely, if  $y \in F(A)$ then there exists $x \in A$ with
 $y \in F(x)$. Since $F = f \cup f \circ F$, either $y \in f(x)$ or there exists $z \in F(x)$ such that $y \in f(z)$.
 Since $A$ is $F$ $^+$invariant, $z \in A$. Hence, $y \in f(A)$.

 In particular, $f(A) = A$ iff $F(A) = A$.

 (b) Clearly, $X$ is $F^{-1}$ $^+$invariant.  Inverting the assumed equation, we have $F^{-1} = f^{-1} \cup f^{-1} \circ F^{-1}$
 and so by (a), $f^{-1}(X) = F^{-1}(X)$.

 (c )If $(y,x) \in F \subset 1_X \cup F_1$ with $y \not= x$ then $(y,x) \in F_1$.
 If $(x,x) \in F$ then either $(x,x) \in f \subset F_1$ or there exists $y \not= x$ such that $(x,y) \in f \subset F_1$
 and $(y,x) \in F \subset 1_X \cup F_1$. Since $y \not= x$, $(y,x) \in F_1$.  By transitivity, $(x,x) \in F_1$. Thus,
 $F \subset F_1$. Similarly, $F_1 \subset F$.

 (d) In any case, suppose $x$ is a regular point for $F$, i.e.\ $L$ on $F(x)$ is greater than $L(x)$ and $L$ on
 $F^{-1}(x)$ is less than $L(x)$. Since $f \subset F$, $x$ is a regular point for $F$.
 Conversely, suppose $x$ is regular for $f$ and $y \in F(x)$. Since $f \cup F \circ f = F $ either $y \in f(x)$
 and so $L(y) > L(x)$ or there exists $z \in f(x)$ such that $y \in F(z)$. Hence, $L(y) \geq L(z) > L(x)$.
 The argument for $y \in F^{-1}$ is similar, using $f^{-1} \cup F^{-1} \circ f^{-1} = F^{-1} $.

(e) If $f$ is a map, then $f \circ f^{-1} \subset 1_X$. Hence,
\begin{equation}\label{eq14uaaa}
\begin{split}
F \circ f^{-1} \ = \ (f \cup F \circ f) \circ f^{-1} \ \subset \ 1_X \cup F, \\
\text{and} \qquad f \circ F^{-1} \subset 1_X \cup F^{-1}, \hspace{2cm}
\end{split}
\end{equation}
where the second equation follows from the first by inverting. Hence, $(f \times f)(F) = f \circ F \circ f^{-1} \subset F$.
Since $f$ maps $F$ to itself, it maps $F^{-1}$ to itself and $F \cap F^{-1}$ to itself. In particular, each $F \cap F^{-1}$
equivalence class is mapped into some equivalence class. If $x$ is in the $F \cap F^{-1}$ equivalence class $E$, then,
since $F = f \cup F \circ f$,
either $f(x) = x$ or $(f(x),x) \in F$. Because $(x,f(x)) \in f \subset F$ it follows that $f(x) \in E$.  Thus, each
$E$ is mapped into itself by $f$.

Now assume that $F = f \cup f \circ F$ and that $x,y$ are in the $F \cap F^{-1}$ equivalence class $E$.
Since $(x,y) \in F$, either $y = f(z)$ with $z = x$ or there exists $z$ such that $(x,z) \in F$ and $f(z) = y$. Since
$(z,y) \in f \subset F$, $z \in E$. In either case, there $z \in E$ with $f(z) = y$. Thus, $f(E) = E$.

 $\Box$ \vspace{.5cm}

 \begin{prop}\label{prop16aeee} Let $f$ be a relation on a normal
 Hausdorff space $X$, with $\U_M$ the maximum uniformity on
 $X$. If $\G f$ and $\G f^{-1}$ are usc relations, i.e.\ for every closed subset
 $A$ of $X$, both $\G f(A)$ and $\G f^{-1}(A)$
 are closed, then $1_X \cup \G f = 1_X \cup \A_{\U_M} f$. If, in addition, $f$ is
 cusc, then $\G f = \A_{\U_M} f$. \end{prop}

 {\bfseries Proof:} In any case, $\A_{\U_M} f$ is a closed, transitive relation which contains $f$ and so contains $\G f$.

 If $(x,y) \not\in 1_X \cup \G f$ then let $X_0 = \{ x, y \}$. Let $L_0(x) = 1$ and $L_0(y) = 0$.
 Since $(x,y) \not\in 1_X \cup \G f$, $L_0$ is a Lyapunov function on $X_0$. By Theorem \ref{theo16ae}
 there exists
 a Lyapunov function $L$ for $\G f$ with $L(x) = 1$ and $L(y) = 0$. By Corollary \ref{cor15u} $L$ is an  $\A_{\U_M} f$
 Lyapunov function and so $(x,y) \not\in \A_{\U_M} f$.

 If $f$ is cusc then by  Proposition \ref{prop16aj},
 we can apply Proposition \ref{prop16ak} (c) to obtain $\G f = \A_{\U_M} f$.

 $\Box$ \vspace{.5cm}

 We require the following lemma from  \cite{AA10a}.

 \begin{lem}\label{lemcom1} Let $f$ be a proper relation on a paracompact, locally compact, Hausdorff space $X$.
 There exists a clopen equivalence relation $E_f$ on $X$ such that
 $\CC_{\U_M} f \cup \CC_{\U_M} f^{-1} \subset E_f$ and  $E_f(x)$ is a $\s$ compact set for every $x \in X$. \end{lem}

   {\bfseries Proof:}  Since $X$ is paracompact, $\U_M$ consists of all neighborhoods of the diagonal and there
   exists an open cover $\{ U_i \} $ such that $\{ \overline{U_i} \}$ is a locally finite collection of compacta.
   It follows that $W = \bigcup_i \overline{U_i} \times \overline{U_i} $ is a closed, symmetric element of $\U_M$
   with every $W(x)$ compact, i.e.\ $W$ is a pointwise compact relation.  Since $\U_M$ is a uniformity there exists
   $V$ a closed, symmetric element of $\U_M$ such that $V \circ V \subset W$. If $K$ is any compact subset of $X$ then
   there exists $F$ a finite subset of $X$ such that $\{ V(x) : x \in F \}$ is a cover of $K$. Then
   $V(K) \subset \bigcup \{ V \circ V(x) : x \in F \} \subset \bigcup \{ W(x) : x \in F \}$. Since $W$ is pointwise compact,
   the set on the right is compact.  Since $V$ is closed and $K$ is compact, $V(K)$ is closed by   \ref{theo16aa} (a)
   and so is compact. Since a locally compact space is a k-space, it follows from
   \ref{theo16aa} (h)  that $V = V^{-1}$ is cusc and so is proper.

   Since $f$ is proper, i.e.\ $f$ and $f^{-1}$ are cusc, and $1_X$ is proper, it follows from Proposition \ref{prop16ac} (a)
   that $F = f \cup 1_X \cup f^{-1}$ is symmetric and cusc. By Theorem \ref{theo16aaa} (c) the composition
   $V_f = V \circ F \circ V \supset V$ is a cusc, symmetric element of $\U_M$. Hence, $E_f =  \bigcup_{n=1}^{\infty} (V_f)^n$
   is an equivalence relation.   Since $F \circ E_f \subset V_f \circ E_f \subset E_f$, and $E_f$ is a closed, transitive relation,
   it follows that  $\G f \cup \G f^{-1} \subset \G F \subset E_f$.
   Since $E_f(x) \supset V_f(x)$ is a neighborhood of $x$, each $E_f(x)$ is open and since the equivalence
   classes are disjoint, each is clopen. Hence, $E_f = \bigcup_x E_f(x) \times E_f(x)$ is a clopen subset of
   $X \times X$. Beginning with the compact set $\{ x \}$ we see, inductively, that
   $(V_f)^{n+1}(x) = V_f((V_f)^n(x))$ is compact because $V_f$ is proper. Hence, each $E_f(x)$ is $\s$ compact.

   Finally, since $E_f$ is a neighborhood of the diagonal, $E_f \in \U_M$. For $x,y \in X$ let $[a,b] \in f^{\times n}$ be a
   $xy, E_f$ chain. Let $b_0 = x$, and $a_{n+1} = y$. Hence, $(a_i,b_i) \in f \subset E_f$ for $i = 1,\dots,n$ and $(b_i,a_{i+1}) \in E_f$
   for $i = 0,\dots,n$. By transitivity of $E_f$, $(x,y) \in E_f$. Hence, $\CC_{\U_M} f \subset E_f$ and by symmetry  $\CC_{\U_M} f^{-1} \subset E_f$.

 $\Box$ \vspace{.5cm}

  \begin{theo}\label{theocom2} Let $F$ be a closed, transitive relation on a paracompact, locally compact, Hausdorff space $X$
  with $\U_M$  the uniformity of all neighborhoods of the diagonal.
   Assume that $X_0$ is a closed subset of $X$ and $L_0 : X_0 \to [a,b]$ is a bounded, Lyapunov function for the
 restriction $F_0 = F \cap (X_0 \times X_0)$.  If either
 \begin{itemize}
 \item[(a)] $X$ is $\s$-compact, or,
 \item[(b)] there exists a proper relation $f$ on $X$ such that $F \subset \CC_{\U_M} f$,
 \end{itemize}
then there exists
 $L : X \to [a,b]$ a  Lyapunov function for $F$ such that $L(x) = L_0(x)$ for $x \in X_0$.
 \end{theo}

 {\bfseries Proof:} (a) Because $X$ is locally compact and $\s$ compact
 there is an increasing sequence of compacta $\emptyset = K_0, K_1, \dots $ with union $X$ such that $K_n \subset K_{n+1}^{\circ}$.
 Let $K_{n+\frac{1}{2}} = K_n \cup (X_0 \cap K_{n+1})$. Assume we have a Lyapunov function $L_{n}: X_n \to [a,b]$ for $F \cap (K_n \times K_n)$
 with $L_{n} = L_0$ on $X_0 \cap K_n$. Extend to define $L_{n+\frac{1}{2}}: X_{n+\frac{1}{2}} \to [a,b]$ by using $L_0$ on $X_0 \cap K_{n+1}$.
  By Theorem \ref{theo16ae}
 there exists a
Lyapunov function $L_{n+1}: K_{n+1} \to [a,b]$ for $F \cap (K_{n+1} \times K_{n+1}$ such that $L_{n+1}$ extends $L_{n+\frac{1}{2}}$. Completing the
inductive construction we define
$L : X \to [a.b]$  by $L|K_n = L_n$. Since $X = \bigcup_n (K_n)^{\circ}$, $L$ is continuous and
so is the required Lyapunov function.

(b) Let $E_f$ be a clopen equivalence relation on $X$ as given by Lemma \ref{lemcom1}. Each equivalence class $E$ is $\s$-compact and
$F \cup F^{-1}$ $^+$invariant. Use (a) on $E$ to define  $L_E : E \to [a,b]$ a Lyapunov function for $F\cap (E \times E)$ which extends
$L_0|(X_0 \cap E)$. Define $L$ by $L|E = L_E$ for each equivalence class.  $L$ extends $L_0$. As the equivalence classes are clopen, $L$ is continuous.
Finally, $F = \bigcup_E (F \cap (E \times E))$ and so $L$ is a Lyapunov function for $F$.
 $\Box$ \vspace{.5cm}

 \begin{cor}\label{corcom2a} Let $f$ be a relation on a paracompact, locally compact, Hausdorff space $X$ with $\U_M$
 the uniformity of all neighborhoods of the diagonal.
\begin{itemize}
 \item[(a)] If $X$ is $\s$ compact, then $1_X \cup \G f = 1_X \cup \A_{\U_M} f$. If, in addition, $f$ is
 cusc, then $\G f = \A_{\U_M} f$.

 \item[(b)]If $f$ is a proper relation, then $\G f = \A_{\U_M} f$.
 \end{itemize}
 \end{cor}

 {\bfseries Proof:}  As in Proposition \ref{prop16aeee} it suffices to show that if $(x,y) \not\in 1_X \cup \G f$
 there is a Lyapunov function $L$ for $\G f$ with $L(x) = 1$ and $L(y) = 0$.

 If $(x,y) \not\in 1_X \cup \G f$ then let $X_0 = \{ x, y \}$. Let $L_0(x) = 1$ and $L_0(y) = 0$.
 Since $(x,y) \not\in 1_X \cup \G f$, $L_0$ is a Lyapunov function on $X_0$.
  By Theorem \ref{theocom2}
 there exists a $\G f$
Lyapunov function $L: X \to [0,1]$ which extends $L_0$. Hence,
$L : X \to [0.1]$ is uniquely defined by $L|X_n = L_n$. Since $X = \bigcup_n (X_n)^{\circ}$, $L$ is continuous and
so is the required Lyapunov function.

When $f$ is cusc, as in (b), we obtain $\G f = \A_{\U_M} f$ from Proposition \ref{prop16ak}(c).

 $\Box$ \vspace{.5cm}

Now we consider extensions to  completions and compactifications.

If $X$ is a compact, Hausdorff space, then $\U_M$ is the unique
uniformity on $X$ and we write $\CC f$\index{cc@$\CC f$} for $\CC_{\U_M} f$ in the compact case.
If a  compact space $X$ is metrizable, then by Theorem \ref{theo11u},
$\CC f = \CC_d f$ for every continuous metric $d$ on $X$.  Since a compact Hausdorff space is normal and every closed
relation on a compact Hausdorff space is proper, it follows from
Proposition \ref{prop16aeee} that $\A_{\U_M} f = \G f$ when $X$ is compact.

\begin{prop}\label{prop12u} Let $(X,\U)$ be a Hausdorff uniform space with completion $(\bar X, \bar \U)$ so that
$X$ is a dense subset of $\bar X$ with $\U$ the uniformity on $X$ induced from $\bar \U$.  If $f$ is a closed
relation on $X$ and $\bar f$ is the closure of $f$ in $\bar X \times \bar X$ then
\begin{equation}\label{eq52u}
\begin{split}
\bar f \cap (X \times X) \ = \ f, \quad \CC_{\bar \U} \bar f \cap (X \times X) \ = \ \CC_{ \U} f, \\
 \A_{\bar \U} \bar f \cap (X \times X) \ = \ \A_{ \U} f. \hspace{3cm}
\end{split}
\end{equation}
If, moreover, $f$ is a uniformly continuous map on $(X,\U)$ then $\bar f$ is a uniformly continuous map
on  $(\bar X, \bar \U)$
\end{prop}

{\bfseries Proof:} $\bar f \cap (X \times X) = f$ because $f$ is closed in the topology of $X \times X$ which is the
relative topology from $\bar X \times \bar X$.

Since $f \subset X \times X$, and the pseudo-metrics of $\Gamma(\U)$ are the restrictions of
the pseudo-metrics in $\Gamma(\bar \U)$
it follows from (\ref{eq16ab}) that $\CC_{\bar \U}  f \cap (X \times X) \ = \ \CC_{ \U} f$. On the other hand,
 (\ref{eq20aa}) implies that $\CC_{\bar \U} \bar f = \CC_{\bar \U}  f$. Similarly, for $\A_{\bar \U}$.

 If $f$ is a uniformly continuous map and $\bar x \in \bar X$ then there is a net $\{ x_i : i \in I \}$ in $X$ which
 converges to $\bar x$. Since $f$ is uniformly continuous, $\{ f(x_i) \}$ is Cauchy and so converges to a point
 $y$ with $(x,y) \in \bar f$. If $\{ x_j : j \in J \}$ is another net converging to $x$, then let
 $\{ 0, 1 \}$ be directed by the relation $\{ 0, 1 \} \times \{ 0, 1 \}$. On $I \times J \times \{0, 1 \}$
 define the net by $(i,j,0) \mapsto x_i$ and $(i,j,1) \mapsto x_j$. This net converges to $x$ and so
 the limit points of $\{ f(x_i) \}$ and $\{ f(x_j) \}$ agree.  Thus, $\bar f$ is a well-defined map on $\bar X$.
 Since the uniformity $\bar \U$ is generated by the closures of $U \in \U$ it is easy to see that $\bar f$ is
 uniformly continuous.

$\Box$ \vspace{.5cm}

Let $\B$ be a closed subalgebra of the Banach algebra $\B(X,\U)$ of bounded
continuous functions on a Hausdorff uniform space.
For any transitive relation $F$ on $X$, the set of those $L \in \B$ which are Lyapunov functions for $F$
always satisfies ALG and CON.

 If $\B$ distinguishes points
and closed sets then it generates a totally bounded uniformity $\T(\B) \subset \U$ with topology compatible to that of
$(X,\U)$, see Appendix B. Let $(\bar X, \bar \T(\B))$ be the
completion of $(X, \T(\B))$. The space $\bar X$
is a compact, Hausdorff space with $\bar \T(\B)$ its unique uniformity. The inclusion $(X,\T(\B))$ into
$(\bar X,\bar \T(\B))$ is a uniform isomorphism onto its image and so the inclusion from $(X,\U)$ is a uniformly
continuous homeomorphism.

If $h : (X_1,\U_1) \to (X_2,\U_2)$ is uniformly continuous, then \\ $h^* : \B(X_2,\U_2) \to \B(X_1,\U_1)$
with $h^*(u) = u \circ h$ is a map of Banach algebras with norm $1$. If $\B_1 \subset \B(X_1,\U_1)$ and
$\B_2 \subset \B(X_2,\U_2)$ are closed subalgebras such that $h^*(\B_2) \subset \B_1$ then
$h : (X_1,\T(\U_1)) \to (X_2,\T(\U_2))$ is uniformly continuous because for $u \in \B_2$
$h^*d_u$, that is, $d_{u} \circ (h \times h)$, is equal to $d_{h^*u}$.

\begin{lem}\label{lem13ua} Suppose that $r : X \times X \to \R$ is a bounded, uniformly continuous map and let
$D$ be a dense subset of $X$. For  $z \in X$, the function $r^z : X \to \R$ is defined by $r^z(x) = r(x,z)$.
If for every $z \in D$, the function $r^z $ is contained in a closed subalgebra $\B$ of $\B(X,\U)$ then
$r^z \in \B$ for all $z \in X$. \end{lem}

{\bfseries Proof:} By uniform continuity, $z \mapsto r^z$  is a continuous map from $X$ to $\B(X,\U)$.  If the dense set
$D$ is mapped into the closed subset $\B$ then all of $X$ is.

$\Box$ \vspace{.5cm}

\begin{theo}\label{theo13ua} Let $f$ be a closed relation on a Hausdorff uniform space $(X,\U)$.
There exists $\B$ a closed subalgebra of $ \B(X,\U)$ such that
\begin{itemize}
\item $\B$ distinguishes points and closed sets in $X$.
\item The set of $f$ Lyapunov functions in $\B$ satisfies POIN  for $\A_{\U} f$.
\end{itemize}
With $\T(\B) \subset \U$  the totally bounded uniformity generated by $\B$, let $(\bar X, \bar \T(\B))$ be the
completion of $(X, \T(\B))$ and let $\bar f$ be the closure of $f$ in $\bar X \times \bar X$. The space $\bar X$
is a compact, Hausdorff space with $\bar \T(\B)$ its unique uniformity. Furthermore,
\begin{equation}\label{eq53ua}
\bar f \cap (X \times X) \ = \ f, \quad
\quad 1_X \cup \G \bar f \cap (X \times X) \ = \ 1_X \cup \A_{ \U} f.
\end{equation}
If $f$ is cusc then $\G \bar f \cap (X \times X) \ = \  \A_{ \U} f$.

If $f$ is a uniformly continuous map, and so is cusc,  such that $f^*\B \subset \B$,
then $\bar f$ is a continuous map on $\bar X$.
If $f$ is a uniform isomorphism such that $f^*\B = \B$, then $\bar f$ is a homeomorphism on $\bar X$.
\end{theo}

{\bfseries Proof:} Since $X$ is a Tychonoff space, $\B = \B(X,\U)$ distinguishes points and closed sets.
The set of functions which are $K \ell_d^f$ dominated for some positive
$K$ and some $d \in \Gamma(\U)$ is a collection of $\A_{\U} f$ Lyapunov functions which satisfies POIN.

Now assume that $\B$ is a closed subalgebra which satisfies these two conditions.

 To prove (\ref{eq53ua}) it
suffices, by (\ref{eq52u}) to show that on $X$ that
$ 1_X \cup \A_{\T(\B)} f = 1_X \cup \A_{\U} f$ because $\A_{\T(\B)} \bar f = \G \bar f$ for the compact Hausdorff space $\bar X$.

Because $\T(\B) \subset \U$, $1_X \cup \A_{\U} f \subset 1_X \cup \A_{\T(\B)} f$. If $(x,y) \not \in 1_X \cup \A_{\U} f$
then by POIN there exists $L \in \B$ a Lyapunov function for $f$ such that $L(x) > L(y)$. Because $L \in \B$
it is uniformly continuous with respect to $\T(\B)$. By Theorem \ref{theo14u} $L$ is an $\A_{\T(\B)} f$ Lyapunov function.
Since $L(x) > L(y)$, $(x,y) \not\in \A_{\T(\B)} f$.

If $f$ is cusc then by Proposition \ref{prop16aj}   $\A_{\U} f \ = \ f \cup (\A_{\U} f) \circ f$ and
$ \A_{\T(\B)} f \ = \ f \cup (\A_{\T(\B)} f) \circ f$. So by  Proposition \ref{prop16ak}(b) we may remove $1_X$ from the equation.

If $f^*\B \subset \B$ then $f$ is uniformly continuous on $(X,\T(\B))$ and so
extends to a continuous map on the completion.
If $f$ is invertible and $f^*\B = \B$ then the same applies to the inverse of $f$.

$\Box$ \vspace{.5cm}

\begin{theo}\label{theo13ub} Let $f$ be a closed relation on a Hausdorff uniform space $(X,\U)$.
There exists $\B$ a closed subalgebra of $ \B(X,\U)$ such that
\begin{itemize}
\item $\B$ distinguishes points and closed sets in $X$.
\item The set of elementary $\U$ Lyapunov functions for $f$ in $\B$ satisfies POIN-E  for $\CC_{\U} f$.
\end{itemize}
With $\T(\B) \subset \U$  the totally bounded uniformity generated by $\B$, let $(\bar X, \bar \T(\B))$ be the
completion of $(X, \T(\B))$ and let $\bar f$ be the closure of $f$ in $\bar X \times \bar X$. The space $\bar X$
is a compact, Hausdorff space with $\bar \T(\B)$ its unique uniformity.  Furthermore,
\begin{equation}\label{eq53ub}
\bar f \cap (X \times X) \ = \ f, \quad
\CC \bar f \cap (X \times X) \ = \ \CC_{ \U} f.
\end{equation}

If $f$ is a uniformly continuous map, and so is cusc,  such that $f^*\B \subset \B$,
then $\bar f$ is a continuous map on $\bar X$.
If $f$ is a uniform isomorphism such that $f^*\B = \B$, then $\bar f$ is a homeomorphism on $\bar X$.

\end{theo}

{\bfseries Proof:} Again it suffices to use $\B = \B(X,\U)$ and as before
it suffices to prove on $X$ that $\CC_{\T(\B)} f = \CC_{\U} f$.

Because $\T(\U) \subset \U$, $\CC_{\U} f \subset \CC_{\T(\B)} f$.

If $(x,y) \not \in 1_X \cup \CC_{\U} f$
then by POIN-E there exists $L \in \B$ an elementary  Lyapunov function for $f$ such that $L(x) > L(y)$. Because $L \in \B$
it is uniformly continuous with respect to $\T(\B)$ and so is a $\T(\B)$ elementary Lyapunov function for $f$.
By Theorem \ref{theo16ai} $L$ is anelementary Lyapunov function for $\CC_{\T(\B)} f$.
Since $L(x) > L(y)$, $(x,y) \not\in \CC_{\T(\B)} f$.

In this case we can eliminate the $1_X$ term without assuming that $f$ is cusc.

If $(x,x) \not\in \CC_{\U} f$, i.e.\ $x \not\in |\CC_{\U} f|$, then by POIN-E there exists $L \in \B$ an
elementary  Lyapunov function for $f$ such that $1 > L(x) > 0$. As before $L$ is an elementary Lyapunov function for $\CC_{\T(\B)} f$.
Hence, $L = 1$ on $\CC_{\T(\B)} f(x)$ and so $(x,x) \not\in \CC_{\T(\B)} f$.

The map cases are as before.

$\Box$ \vspace{.5cm}

The spaces we obtain from these theorems are quite large. The conditions may
well require $\B = \B(X,\U)$, leading to the
entire uniform version of the Stone-\v{C}ech compactification.  However, in
the second countable case we are able to obtain
a metric compactification.

\begin{theo}\label{theo13uc}  Let $f$ be a closed relation on a Hausdorff uniform space $(X,\U)$ with $X$ second countable.
There exists $\B$ a separable, closed subalgebra of $ \B(X,\U)$ such that
\begin{itemize}
\item $\B$ distinguishes points and closed sets in $X$.
\item The set of $f$ Lyapunov functions in $\B$ satisfies POIN  for $\A_{\U} f$.
\item The set of elementary $\U$ Lyapunov functions for $f$ in $\B$ satisfies POIN-E  for $\CC_{\U} f$.
\end{itemize}
With $\T(\B) \subset \U$  the totally bounded uniformity generated by $\B$, let $(\bar X, \bar \T(\B))$ be the
completion of $(X, \T(\B))$ and let $\bar f$ be the closure of $f$ in $\bar X \times \bar X$. The space $\bar X$
is a compact, metrizable Hausdorff space with its unique uniformity $\bar \T(\B)$ metrizable.
Furthermore,
\begin{equation}\label{eq53uc}
\begin{split}
\bar f \cap (X \times X) \ = \ f, \quad
1_X \cup \G \bar f \cap (X \times X) \ = \ 1_X \cup \A_{ \U} f, \\
\CC \bar f \cap (X \times X) \ = \ \CC_{ \U} f. \hspace{3cm}
\end{split}
\end{equation}
If $f$ is cusc then $\G \bar f \cap (X \times X) \ = \   \A_{ \U} f$.

If $f$ is a uniformly continuous map then, in addition, we can choose $\B$ so
that $f^*\B \subset \B$ and so $\bar f$ is a continuous map on $\bar X$.
If $f$ is a uniform isomorphism then, in addition, we can choose $\B$ so that $f^*\B = \B$ and so
 $\bar f$ is a homeomorphism on $\bar X$.
\end{theo}

{\bfseries Proof:}  Apply Theorem \ref{theo9u}   to obtain a metric $d \in \Gamma(\U)$ with
the topology that of $X$ and such that $ \A_{\U} f = \A_d f$,
$1_X \cup \A_{\U} f = \ \leq_{L_{\ell}}$ and $\CC_{\U} f = \CC_d f$.
 Let $D$ be a countable dense subset of $X$.

Let $d^z(x) = d(x,z)$. Let $\ell^z(x) = \ell^{f \cup 1_X}(x,z) = \min(\ell^f(x,z), d(x,z))$. If $\B$ is a closed
subalgebra of $\B(X,d) \subset \B(X,\U)$ which contains $d^z$ and $\ell^z$ for all $z$ in $D$ then by Lemma \ref{lem13ua}
 $d^z, \ell^z \in \B $ for all $z \in X$. Since $d^z \in \B$ for all $z$, $\B$ distinguishes points and closed sets. Each $\ell^z$ is
 a Lyapunov function for $\A_d f$ by Theorem \ref{theo10} and Proposition \ref{prop11}. If $(x,y) \not\in \A_d f$ then
 $\ell^y(y) = 0$ and $\ell^y(x) > 0$. So the Lyapunov functions in $\B$ satisfy POIN for $\A_{\U} f = \A_d f$.

 Because the subspaces $X \times X \setminus (1_X \cup \CC_{\U} f)$ and $X \setminus |\CC_{\U} f|$ are Lindel\"{o}f, Theorem \ref{theo16ah} implies
 that we can find a sequence $\{ L_i \}$ of $\U$ elementary Laypunov functions for $f$ such that
 \begin{itemize}
 \item For $(x,y) \in X \times X \setminus (1_X \cup \CC_{\U} f)$
 there exists $i$ such that $L_i(x) > L_i(y)$.
 \item For $x \in X \setminus |\CC_{\U} f|$ there exists $i$ such that $1 > L_i(x) > 0$.
 \end{itemize}

If $\B$ contains $\{ L_i \}$ then the elementary Lyapunov functions in $\B$ satisfy POIN-E.

Thus, if $\B$ is the closed subalgebra generated by $\{ d^z : z \in D \} \cup \{ \ell^z : z \in D \} \cup \{ L_i \}$ then
$\B$ is a separable subalgebra of $\B(X,d)$ which satisfies the required properties.

If $f$ is a uniformly continuous map we extend the countable set of generators $\{ u_i \}$ to include
$\{ (f^n)^* u_i \}$ for all positive integers $n$.  If $f$ is a uniform
isomorphism we use $\{ (f^n)^* u_i \}$ with all  integers $n$. In either case,
we still have a countable set of generators and
so obtain a separable algebra $\B$.

Since $\B$ is separable, the compact space $\bar X$ is metrizable.

The results then follow from Theorems \ref{theo13ua} and  \ref{theo13ub}.

$\Box$ \vspace{.5cm}

Let $f$ be a closed relation relation on $X$ and let $\bar f$ be the extension to one of the
compactifications as above $\bar X$. If the domain of $f$, $f^{-1}(X)$ is all of $X$, then it is dense in $\bar X$.
Since the domain $\bar f^{-1}(\bar X)$ is compact and contains $f^{-1}(X)$, it follows that
$\bar f^{-1}(\bar X) = \bar X$. If $\bar X$ is merely
a completion but $f$ is a uniformly continuous map on $X$ then $\bar f$
is a uniformly continuous map on $\bar X$ and so has
domain all of $\bar X$.  If $f$ is merely continuous, the domain of $\bar f$ need not be all of $\bar X$.  For example,
let $f : (0,\infty) \to (0,\infty)$ be the continuous map with $f(t) = 1/t$.
With the usual metric the completion is $[0,\infty)$
and $\bar f = f$.

We conclude the section by considering the special results when $X$ is a compact Hausdorff space, so the $\U_M$ is its unique
uniformity. We need the following result which is Lemma 2.5 from \cite{A93}. Recall that a closed relation $f$ on a compact Hausdorff
space $X$ is proper and so $f(A)$ is closed if $A \subset X$ is closed.

\begin{lem}\label{comlem1} Let $F$ be a closed, transitive relation on a compact Hausdorff space $X$ and let $B$ be a closed subset with
$B \cap |F| = \emptyset$.  There exists a positive integer $N$ such that if $\{a_0, \dots, a_k \}$ is a finite sequence in $B$ with
$(a_{i-1},a_i) \in F$ for $i = 1, \dots, k$, then $k \leq N$.\end{lem}

{\bfseries Proof:} Since $F \cap (B \times B)$ is disjoint from $1_X$, there exists an open, symmetric $U \in \U$ such that
$F \cap (B \times B) \cap (U \circ U) = \emptyset$. Since $B$ is compact, there is a subset $\{x_0, \dots, x_N \}$ of $B$ such that
$\{ U(x_j) : j = 0, \dots, N \}$ covers $B$. If $\{a_0, \dots, a_k \}$ is a sequence as above with $k > N$ then by the Pigeonhole Principle
there exist $0 \leq i_1 < i_2 \leq k$ which lie in the same $U(x_j)$ and so $(a_{i_1},a_{i_2}) \in U \circ U$. By transitivity of $F$,
$(a_{i_1},a_{i_2}) \in F$, contradicting the choice of $U$.

$\Box$ \vspace{.5cm}

\begin{prop}\label{comprop2} Let $f$ be a closed  relation on a compact Hausdorff space $X$ and let $A$ be a nonempty, closed subset of $X$.
\begin{enumerate}
\item[(a)] If $A$ is $f$ $^+$invariant and $A \subset Dom(f)$,
then maximum closed $f$ invariant subset $f^{\infty}(A)$ is closed and nonempty and equals
$\bigcap_{n \in \N}  f^n(A)$.

\item[(b)] If $F$ is a closed, transitive relation on $X$ such that
$F = f \cup f \circ F$ and $A$ is $F$ $^+$invariant, then $f^{\infty}(A) = F^{\infty}(A)$.
\end{enumerate}
\end{prop}

{\bfseries Proof:} (a) Since $Dom(f) \subset A$, $\{ f^n(A) \}$ is a non-increasing sequence of nonempty compacta and so the intersection
is nonempty. If $y \in \bigcap_{n \in \N}  f^n(A)$ then $\{ f^{-1}(y) \cap f^n(A) \}$ is a  non-increasing sequence of nonempty compacta with
nonempty intersection $f^{-1}(y) \cap \bigcap_{n \in \N}  f^n(A)$. So $\bigcap_{n \in \N}  f^n(A)$ is an $f$ invariant subset.

(b)  By Proposition \ref{prop16ak} and induction, $f^n(A) = F^n(A)$ for all $n$. Hence the intersections are equal.

$\Box$ \vspace{.5cm}

\begin{theo}\label{comtheo3}  Let $F$ be a closed, transitive relation on a compact Hausdorff space $X$. If  $A$ is an $F$ $^+$invariant closed subset,
then
\begin{equation}\label{eqcom1}
F^{\infty}(A) \ = \ F(A \cap |F|). \hspace{4cm}
\end{equation}

If $G$ is an open set containing $A$ then there exists a Lyapunov function $L : X \to [0,1]$ for $F$ such that $L = 0$ on $X \setminus G$ and
$L = 1$ on $A$. In particular, the $F$ $^+$invariant open neighborhoods of $A$ form a base for the neighborhood system of $A$.
\end{theo}

{\bfseries Proof:} Since $x \in F(x)$ for $x \in |F|$, $A \cap |F| \subset F^{\infty}(A)$. From invariance of $F^{\infty}(A)$ we obtain
$F(A \cap |F|) \subset F^{\infty}(A)$.

For $x \in A \setminus F(A \cap |F|)$, $B = (\{ x \} \cup F^{-1}(x)) \cap A$ is closed and disjoint from $|F|$. If $y \in B \cap F^{n}(A)$ then
there exists a sequence $a_0, a_1, \dots, a_n \in A$ with $a_n = y$ and with $a_i \in F(a_{i-1})$ for $i = 1, \dots,n $. From transitivity
of $F^{-1}$ it follows that $a_i \in B$ for all $i$.
From Lemma \ref{comlem1} it then follows that
there exists a positive integer $N$ such that $B$ is disjoint from $F^{N+1}(A)$. Hence, $x \not\in F^{\infty}(A)$.

If $G$ is an open set containing $A$ then we let $X_0 = (X \setminus G) \cup A$. Let $L_0 = 0$ on $X \setminus G$ and $= 1$ on $A$.
Since $A$ is $F$ $^+$invariant, $L_0$ is a Lyapunov function on $X_0$ for $F \cap (X_0 \times X_0)$.  By Theorem \ref{theo16ae}
it extends to an $F$  Lyapunov function $L$ on $X$.

For any $c \in (0,1)$ the set $\{ x : L(x) > c \}$ is an $^+$invariant neighborhood of $A$ which is contained in $G$.

$\Box$ \vspace{.5cm}

These results apply directly to $F = \G f = \A_{\U_M} f $ for $f$ any closed relation on $X$, see  Proposition \ref{prop16aj} and Corollary \ref{corcom2a}.
For $F = \CC f = \CC_{\U_M} f$ we obtain special results.

If $K$ is a closed $\CC f$ invariant set, we call $K \cap |\CC f|$ the \emph{trace}
\index{trace}\index{attractor!trace of} of $K$.

\begin{theo}\label{comtheo4}  Let $f$ be a closed relation on a compact Hausdorff space $X$. Let $K$ be a
subset of $X$.
\begin{enumerate}
\item[(a)] Assume $K$ is closed and $\CC f$ $^+$invariant.
If $G$ is an open set which contains $K$, then there exists an open inward set $A$ with $K \subset A \subset G$ and
there exists an elementary
Lyapunov function $L : X \to [0,1]$ for $f$ such that $L = 0$ on $X \setminus G$ and
$L = 1$ on $K$.
 In particular, the open inward sets which
contain a closed, $\CC f$ $^+$invariant set form a neighborhood base of the set.

The intersection $K \cap |\CC f|$ is a closed, $\CC f \cap(|\CC f| \times |\CC f|)$ invariant subset of $|\CC f|$.

\item[(b)]  If $K$ is closed, then following conditions
are equivalent.
\begin{itemize}
\item[(i)] $K$ is $\CC f$ $^+$invariant and is $f$ invariant.

\item[(ii)] $K$ is $\CC f$ invariant.

\item[(iii)] $K = \CC f(K \cap |\CC f|)$.

\item[(iv)]  $K$ is $\CC f$ $^+$invariant and if $A$ is an inward set which
contains $K$ then the associated attractor $A_{\infty}$ contains $K$.

\end{itemize}

\item[(c)] $K$ is an attractor iff  $K$ is closed, $\CC f$ invariant and $K \cap |\CC f|$ is a
clopen subset of $|\CC f|$. Conversely, if $A_0$ is a clopen $\CC f \cap (|\CC f| \times |\CC f|)$
invariant subset of $|\CC f|$
    then $K_0 = \CC f(A_0)$ is an attractor of which $A_0$ is the trace.
\end{enumerate}
\end{theo}

{\bfseries Proof:} (a) We apply the notation of the proof of Proposition \ref{prop16af}(b). Since $K$ is assumed to be closed and
$\CC f$ $^+$invariant, $K$ is compact and equals $K \cup \CC f(K)$. Let $Q_d(K,y) = \min(m_d^f(K,y),d(K,y)$. From the Proposition we see that
$K$ is the intersection of the inward sets $ \{ y : Q_d(K,y) < \ep \}$ as $(d,\ep)$ varies with $d \in \Gamma(\U_M), \ep > 0$.
Recall that with $d = d_1 + d_2$ and $\ep = \min(\ep_1, \ep_2)$, $\{ y : Q_d(K,y) \leq \ep \} $
$\subset \{ y : Q_{d_1}(K,y) \leq \ep_1 \} \cap \{ y : Q_{d_2}(K,y) \leq \ep_2 \}$. It follows from compactness that for some
$d \in \Gamma(\U_M), \ep > 0$, the compact set $\{ y : Q_d(K,y) \leq \ep \}$ is contained in $G$. Hence, $A = \{ y : Q_d(K,y) < \ep \}$ is
an open inward set with $K \subset A \subset U$.

If $A$ is an open inward set containing $K$ then $X \setminus A$ and $K \cup (\CC f)(A)$ are disjoint
closed sets. Since a compact Hausdorff space is normal, there exists a continuous $L: X \to [0,1]$
which $= 0 $ on $X \setminus A$ and $= 1$ on  $K \cup (\CC f)(A)$. Any such is clearly the required
elementary Lyapunov function.

Since $x \in (\CC f)(x)$ for $x \in |\CC f|$ it is clear that
$K \cap |\CC f|$ is a closed, $\CC f \cap(|\CC f| \times |\CC f|)$ invariant subset of $|\CC f|$.

(b) (i) $\Leftrightarrow$ (ii): By Proposition \ref{prop16aj}
$\CC f \ = \ f \cup f \circ (\CC f) \ = \ f \cup  (\CC f) \circ f$ and so $f(K) = \CC f(K)$ by
Proposition \ref{prop16ak} (a).

(ii) $\Leftrightarrow$ (iii): If $K$ is $\CC f$ $^+$invariant then $K$ is $\CC f$ invariant
iff $K = (\CC f)^{\infty}(K)$ and the latter
equals $\CC f(K \cap |\CC f|)$ by (\ref{eqcom1}).

(ii) $\Leftrightarrow$ (iv): If $A$ is $\CC f$ $^+$invariant and contains a $\CC f$ invariant set $K$ then
the maximum  $\CC f$ invariant set $(\CC f)^{\infty}(A)$ contains $K$.  If $A$ is inward
then $A_{\infty} = (\CC f)^{\infty}(A)$ is the associated attractor.

On the other hand, let $K_1 = \CC f(K)$ and assume there exists $x \in K \setminus K_1$.
Let $G = (\CC f)^*(X \setminus \{ x \})$. Since $K_1 \subset X \setminus \{ x \}$, $K \subset G$.
By (a) there exists $A$ an inward set with $K \subset A \subset G$.
So $A_{\infty} \subset \CC f (G) \subset X \setminus \{ x \}$. That is,the associated attractor does not contain $K$.

(c) If $A$ is an inward set and $\hat A$ is a subset such that $\CC f(A) \subset (\hat A)^{\circ}$ and
$\hat A \subset A$, then $\hat A$ is an inward set with $(\CC f)^{\infty}(\hat A) = (\CC f)^{\infty}( A)$,
i.e.\ with the same associated attractor.  In particular we can choose $\hat A$ closed and so we see that
every attractor is closed. Furthermore, $|\CC f| \cap (\CC f)^{\infty}(\hat A) = |\CC f| \cap A^{\circ}$
and so the trace of the attractor is clopen in $|\CC f|$. It is $\CC f \cap(|\CC f| \times |\CC f|)$ invariant
by (a).

Conversely, if $A_0$ is a clopen $\CC f \cap (|\CC f| \times |\CC f|)$ invariant subset of $|\CC f|$
    then $K_0 = \CC f(A_0)$ is a $\CC f$ invariant subset of $X$  by (b), and it is contained in
    the open set $G = X \setminus (|\CC f| \setminus A_0)$ by $\CC f \cap (|\CC f| \times |\CC f|)$ invariance.
    By (a) there exists an inward set $A$ such that $K_0 \subset A \subset G$. Hence, $A \cap |\CC f| = A_0$
    and so $A_{\infty} = (\CC f)^{\infty}(A) = \CC f(A_0) = K_0$. That is, $K_0$ is the attractor associated
    with $A$ and the trace is $A_0$.

$\Box$ \vspace{.5cm}

{\bfseries Remark:} Notice that while an attractor is necessarily closed, a $\CC f$ invariant set need not be.
For example, if $X$ is the Cantor set and $f = 1_X$ then $\CC f = 1_X$ and every subset of $X$ is
$\CC f$ invariant.

\vspace{1cm}

\section{ Recurrence and Transitivity}
\vspace{.5cm}

We first consider recurrence.

 \begin{prop}\label{eqrel1}  Let $f$ be a relation on a uniform space $(X,\U)$ and let $d \in \Gamma(\U)$. Let
$F =  \G f, \A_d f, \A_{\U} f$, $\CC_d f$ or $\CC_{\U} f$.
  \begin{itemize}
 \item[(a)] The relation $F$ is an equivalence relation iff $f^{-1} \subset F$ and $Dom(F) = X$.

 \item[(b)] If $f$ is a continuous map on $X$ then $F$ is an equivalence relation iff $1_X \subset F$.

 \item[(c)] If $\G f$ is an equivalence relation then $\A_d f$, $\A_{\U} f$, $\CC_d f$ and $\CC_{\U} f$ are equivalence
 relations.
 \end{itemize}
 \end{prop}

 {\bfseries Proof:} (a) Clearly, if $F$ is an equivalence relation on $X$ which contains $f$
 then $1_X \cup f^{-1} \subset F$ and
 so $Dom(F) = X$.

 Conversely, if $f^{-1} \subset \G f$ then $\G f^{-1} \subset \G \G f = \G f$ and so, inverting,
 $\G f \subset \G f^{-1}$.  That is, $\G f$ is symmetric. Similarly, if $f^{-1} \subset F$
 for $F =  \A_d f, \A_{\U} f, \CC_d f$ or $\CC_{\U} f$
 then $F$ is symmetric.  If $F$ is symmetric and $Dom(F) = X$, then for any $x \in X$
 there exists $y \in X$ such that $(x,y) \in F$. By symmetry and transitivity, $(y,x), (x,x) \in F$. So $F$ is
 reflexive.

 (b)  If $f$ is a continuous map then  it is a cusc relation and so $F = f \cup  F\circ  f$
 by  Proposition \ref{prop16aj} (a).
  For any $x \in X$ assume $f(y) = x$.
 Since $(y,y) \in F$, either $(y,y) \in f$, i.e.\ $y = f(y) = x$ and so $(x,y) = (y,y) \in F$, or  $(y,x) \in f$
 and $(x,y) \in F$.  As $y$ was an arbitrary element of $f^{-1}(x)$ it follows that
 $f^{-1} \subset F$. Since $f$ is a map,
 $X = Dom(f) \subset Dom(F)$.

(c) If $\G f$ is an equivalence relation then, since it is contained in $F$ it follows
that $1_X \cup f^{-1} \subset F$ and so
$F$ is an equivalence relation by (a).

  $\Box$ \vspace{.5cm}

  \begin{df}\label{eqre11aaa}  Let $f$ be a relation on a uniform space $(X,\U)$ and let $d \in \Gamma(\U)$. For
$F =  \G f, \A_d f, \A_{\U} f$, $\CC_d f$ or $\CC_{\U} f$ we will say that $f$ is
totally $F$ recurrent when $F$ is an equivalence relation.\end{df}
\vspace{.5cm}

 \begin{df}\label{eqrel2} A topological space $X$ is \emph{completely Hausdorff}
 \index{space!completely Hausdorff}\index{completely Hausdorff space} if the Banach
 algebra $\B(X)$\index{b@$\B(X)$} of bounded,
 real-valued continuous functions distinguish the
 points of $X$. \end{df}

 Thus, if $X$ is completely Hausdorff and  $(x,y) \in X \times X \setminus 1_X$,
 there exists a continuous $L_{xy} : X \to [0,1]$ with $L_{xy}(x) = 0$ and
 $L_{xy}(y) = 1$. These maps define a continuous injection into a product
 of copies of $[0,1]$ indexed by the points of
 $X \times X \setminus 1_X$. Conversely, if there is a continuous
 injection from $X$ to a Tychonoff space, then $X$ is
 completely Hausdorff.

 In \cite{Bi} Bing constructs a simple example of a countable, connected Hausdorff space.
 On such a space the only continuous real-valued
 functions are constants and so the space is not completely Hausdorff.

 A subset $A$ of a topological space $X$ is called a \emph{zero-set}\index{zero-set}
 if there exists $u \in \B(X)$ such that $A = u^{-1}(0)$. Clearly, a zero-set in $X$
 is a closed, $G_{\d}$ subset of $X$. The constant functions $0$ and $1$ show that $X$ and
 $\emptyset$ are zero-sets.  If $u,v \in \B(X)$ and
 $w_1 = u \cdot v, w_2 = u^2 + v^2$ then $w_1^{-1}(0) = u^{-1}(0) \cup v^{-1}(0)$
 and $w_2^{-1} = u^{-1}(0) \cap v^{-1}(0)$. Thus, the collection of
 zero-sets is closed under finite unions and finite intersections. If $(X,d)$ is
 a pseudo-metric space and $A$ is a closed subset of $X$ then
 $u(t) = d(t,A)$ is an element of $\B(X)$ such that $A = u^{-1}(0)$.  That is,
 every closed subset of a pseudo-metric space is a zero-set.
 If $X$ is normal and $A$ is a closed,  then $A$ is a zero-set iff it is a
 $G_{\delta}$ set. If $h : X \to Y$ is a continuous function and
 $u \in \B(Y)$ then $h^*u = u \circ h \in \B(X)$ and $(h^*u)^{-1}(0) = h^{-1}(u^{-1}(0))$.
 That is, the continuous pre-image of a zero-set is a zero-set.
 It follows that if $u \in \B(X)$ and $K \subset \R$ is closed then since $\R$ is a
 metric space $K$ is a zero-set and so $u^{-1}(K)$ is a zero-set.

For a topological space $X$, we denote by $\tau X$
\index{t@$\tau X$}
the set $X$  equipped with the weak topology generated by the elements of $\B(X)$.
That is, it is the coarsest topology with respect to which every
 element of $\B(X)$ is continuous. Equivalently, if $h$ is a map to $X$ from a topological space $Y$, then
 $h : Y \to \tau X$ is continuous iff $h^*u \in \B(Y)$ for all $u \in \B(X)$.
 The set of complements of the zero-sets of $X$ forms a basis for the topology of $\tau X$.
 Thus, the closed sets are exactly those which are intersections of the zero-sets of $X$.
  Thus, the ``identity map'' from $ X$ to $\tau X$ is continuous and $\B(\tau X) = \B(X)$.

  \begin{prop}\label{tauprop} Let $X$ be a topological space.
  \begin{enumerate}
  \item[(a)] The following are equivalent.
    \begin{itemize}
  \item[(i)] $X$ is completely regular.
  \item[(ii)] Every closed subset of $X$ is an intersection of zero-sets.
  \item[(iii)] $X = \tau X$.
  \end{itemize}

   \item[(b)] The following are equivalent.
    \begin{itemize}
  \item[(i)] $X$ is completely Hausdorff.
  \item[(ii)] Every point of $X$ is an intersection of zero-sets.
  \item[(iii)] $X$ is a $T_1$ space and every compact subset of $X$ is an intersection of zero-sets.
  \item[(iv)] $X$ is a $T_1$ space and disjoint compact subsets can be distinguished by $\B(X)$.
  \item[(v)] $\tau X$ is a $T_1$ space.
  \item[(vi)] $\tau X$ is a Tychonoff space.
  \end{itemize}

  \item[(c)] The space $\tau X$ is completely regular and if $h : X \to Y$ is a
  continuous function with $Y$ completely regular, then
  $h : \tau X \to Y$ is continuous.

  \item[(d)] If $d$ is a  pseudo-metric on $X$ then $d$ is continuous on $X \times X$
  iff it is continuous on $\tau X \times \tau X$.
  The set of all continuous pseudo-metrics on $X$ is the gage of the maximum uniformity with topology that of $\tau X$.
  \end{enumerate}
  \end{prop}

  {\bfseries Proof:} (a) (i) $\Leftrightarrow$ (ii):  If $\{ u_i \} \subset \B(X)$ and $A = \bigcap_i u_i^{-1}(0)$ then
  for every $x \not\in A$ there exists $u_i$ with $u_i(x) \not= 0$,
 while for all $i$ $u_i = 0$ on $A$. On the other hand, if for every $x \not\in A$
 there exists a $v_x \in \B(X)$ with $v_x(x) \not\in \ol{v_x(A)}$
 then $u_x(y) = d(v_x(y),\ol{v_x(A)})$ then $A \subset u_x^{-1}(0)$ and $u_x(x) \not= 0$.
 Thus, $A$ is the intersection of the $u_x^{-1}(0)$'s as
 $x$ varies over $X \setminus A$. Thus, $A$ is an intersection of zero-sets
 iff $\B(X)$ distinguishes $A$ from the points of $X \setminus A$.

(ii) $\Leftrightarrow$ (iii): The closed sets of $\tau X$ are exactly the intersections of the zero-sets of $X$.

(c) Since $\B(X) = \B(\tau X)$ it is clear that $\tau (\tau X) = \tau X$ and
so $\tau X$ is completely regular by (a). If $A$ is a closed subset of $Y$
then because $Y$ is completely regular, $A$ is an intersection of zero-sets by
(a). Since $h : X \to Y$ is continuous, $h^{-1}(A)$ is an intersection
of zero-sets in $X$ and so is closed in $\tau X$. Thus, $h : \tau X \to Y$ is continuous.

(b) (i) $\Leftrightarrow$ (ii): Just as in (a).

(ii) $\Leftrightarrow$ (iii): If $A$ is compact and $x \not\in A$ then for each
$a \in A$ there exists $u_a \in \B(X)$ such that
$u_a(a) = 0$ and $u_a(x) = 1$. Let $v_a = 2 \max(u_a - \frac{1}{2}, 0)$. That is,
$v_a(x) = 1$ and $v_a = 0$ on a neighborhood of $a$.
By compactness there exists a finite subset $A_0$ of $A$ such that
$u = \Pi_{a \in A_0} v_a$ is $1$ at $x$ and $0$ on $A$. The converse is obvious.

 (iii) $\Leftrightarrow$ (iv): If $A$ and $B$ are disjoint compact sets then for
 every $x \in B$ there exists $u_x = 1$ on $A$ and has $v_x(x) = 0$.
 Use $u_x = 1 - u$ from the above proof. Again, let $v_x = 2 \max(u_x - \frac{1}{2}, 0)$.
 As above, there is a finite subset $B_0$ of $B$ so that
$u = \Pi_{x \in B_0} v_x$ is $1$ on $A$ and $0$ on $B$. Again, the converse is obvious.

(ii) $\Rightarrow$ (v): From (ii), every point is closed in $\tau X$.

(v) $\Rightarrow$ (vi): A $T_1$ completely regular space is Tychonoff.

(vi) $\Rightarrow$ (i): $X$ injects into the Tychonoff space $\tau X$.

(d) If $d$ is a continuous pseudo-metric on $\tau X$ then it is a continuous pseudo-metric on $X$ since $\tau X$ is coarser than $X$.
If $d$ is a continuous pseudo-metric on $X$ then $(X,d)$ is a pseudo-metric space with $X \to (X,d)$ continuous. Since a pseudo-metric space is
completely regular, (c) implies that $\tau X \to (X,d)$ is continuous. Since $d$ is a continuous function on $(X,d) \times (X,d)$, it is
continuous on $\tau X \times \tau X$.  For a completely regular space, like $\tau X$, the collection of all continuous pseudo-metrics is the gage of
the maximum uniformity.

$\Box$ \vspace{.5cm}

A clopen set is clearly a zero-set. Recall that the \emph{quasi-component} of a point $x \in X$ is the intersection of all the clopen sets
which contain $x$. In a compact space the quasi-components are the components, but even in a locally compact space this need not be true.
If $X_0 = [0,1] \times \{0, 1/n : n \in \N \}$ and $X = X_0 \setminus \{(\frac{1}{2},0) \}$ then the quasi-component of $(0,0)$ is
$([0,1] \setminus \{\frac{1}{2} \}) \times \{ 0 \}$.

\begin{df}\label{eqrel2zero} A topological space $X$ is
\begin{itemize}
\item \emph{totally disconnected} when the quasi-components are singletons.
\item \emph{zero-dimensional} when the clopen sets form a basis for the topology.
\item \emph{strongly zero-dimensional} when the clopen sets contain a neighborhood basis for every closed subset.
\end{itemize}
\end{df}
\index{space!totally disconnected}\index{totally disconnected space}\index{space!zero-dimensional}\index{zero-dimensional space}%
\index{space!strongly zero-dimensional}\index{strongly zero-dimensional space}

Recall from Appendix B, that we call a uniformity $\U$ zero-dimensional when it is generated by equivalence relations.

For a space $X$ let $\B_0(X)$ \index{b@$\B_0(X)$}
consist of those $u \in \B(X)$ with $u(X) \subset \{ 0, 1 \}$,
i.e.\ $\B_0(X)$ is the set of characteristic functions
of the clopen subsets. For a topological space $X$, we denote by $\tau_0 X$ \index{t@$\tau_0 X$}
the set $X$
equipped with the weak topology generated by the elements of $\B_0(X)$.
that is, it is the coarsest topology with respect to which every
 element of $\B_0(X)$ is continuous. Equivalently, if $h$ is a map to $X$ from a topological space $Y$, then
 $h : Y \to \tau X$ is continuous iff $h^{-1}(A)$ is clopen in $Y$ whenever $A$ is a clopen subset of $X$.

  \begin{prop}\label{taupropzero} Let $X$ be a topological space.
  \begin{enumerate}
  \item[(a)] The following are equivalent.
    \begin{itemize}
  \item[(i)] $X$ is zero-dimensional.
  \item[(ii)] Every closed subset of $X$ is an intersection of clopen sets.
  \item[(iii)] $X = \tau_0 X$.
  \end{itemize}

  If $X$ is zero-dimensional, then it is completely regular.

   \item[(b)] The following are equivalent.
    \begin{itemize}
  \item[(i)] $X$ is totally disconnected.
  \item[(ii)] Every point of $X$ is an intersection of clopen sets.
  \item[(iii)] $X$ is a $T_1$ space and every compact subset of $X$ is an intersection of clopen sets.
  \item[(iv)] $X$ is a $T_1$ space and if $A, B$ are disjoint compact subsets of $X$ then there exists a clopen set $U$
  with $A \subset U$ and $B \cap U = \emptyset$.
  \item[(v)] $\tau_0 X$ is a $T_1$ space.
    \end{itemize}

If $X$ is totally disconnected, then it is completely Hausdorff.

  \item[(c)] The space $\tau_0 X$ is zero-dimensional and if $h : X \to Y$ is a continuous function with $Y$ zero-dimensional, then
  $h : \tau_0 X \to Y$ is continuous.

  \item[(d)] If $d$ is a  pseudo-ultrametric on $X$ then $d$ is continuous on $X \times X$ iff it is continuous on $\tau_0 X \times \tau_0 X$.
  The set of all continuous pseudo-ultrametrics on $X$ is the gage of the maximum zero-dimensional uniformity with topology that of $\tau_0 X$.
  \end{enumerate}
  \end{prop}

  {\bfseries Proof:}  The proofs are completely analogous to those of Proposition \ref{tauprop}.  The details are left to the reader.

$\Box$ \vspace{.5cm}

\begin{prop}\label{taupropcompact}
\begin{enumerate}
\item[(a)] If a space is compact and totally
disconnected then it is strongly zero-dimensional.

\item[(b)] If a space is locally compact and totally
disconnected then it is zero-dimensional.

\item[(c)] A $T_1$ space is zero-dimensional iff it admits
an embedding into a compact, totally disconnected space.

\item[(d)] A space is totally disconnected iff it admits a continuous injection into a
compact, totally disconnected space.
\end{enumerate}
\end{prop}

{\bfseries Proof:} (a) If $X$ is a compact Hausdorff space then disjoint closed sets are disjoint compact sets.
So (i) $\Rightarrow$ (iv) of Proposition \ref{taupropzero} (b) implies that a compact, totally disconnected space is strongly zero-dimensional.

(b) If $x \in X$ and $x$ is contained in an open set $U$ with closure $\ol{U}$ compact, then there exists a clopen set $A_0$ containing $x$ and
disjoint from the compact set $\ol{U} \setminus U$. Hence, $A = A_0 \cap U = A_0 \cap \ol{U}$ is a clopen set containing $x$ and contained in $U$.

(c), (d) Using the elements of $\B_0(X)$ we can inject totally disconnected space $X$, or embed a $T_1$ zero-dimensional space into a product of
copies of $\{ 0,1 \}$, which is compact and totally disconnected.

Conversely, a subspace of a zero-dimensional space is zero-dimensional and if $X$ injects into a totally disconnected space then it is totally disconnected.

$\Box$ \vspace{.5cm}

\begin{queses} Does there exist a space which is completely Hausdorff and regular, but not completely regular?

Does there exist a completely regular, totally disconnected space which is not zero-dimensional?  In particular, for a totally disconnected
space $X$ is $\tau X = \tau_0 X$?
\end{queses}
\vspace{.5cm}

Call a Hausdorff space $X$ \emph{strongly $\s$-compact}\index{space!strongly $\s$-compact}\index{strongly $\s$-compact space} if
there is a sequence $\{ K_n \}$ of compacta covering $X$ such that $A \cap K_n$ closed for
all $n$ implies $A$ is closed. Equivalently, by taking complements, we have that $A \cap K_n$ is open in $K_n$ for all $n$ implies $A$ is open.
Consequently, if $A \cap K_n$ is clopen in $K_n$ for all $n$ then $A$ is clopen. Observe that the condition is a strengthening of the condition that
$X$ be a $k$-space.

\begin{prop}\label{propstronglysigma}
\begin{enumerate}
\item[(a)] If $X$ is a locally compact, $\s$-compact Hausdorff space then $X$ is strongly $\s$-compact.
\item[(b)] If $q : X \to Y$ is a quotient map with  $Y$ Hausdorff and $X$ Hausdorff and strongly $\s$-compact, then $Y$ is strongly $\s$-compact.
\item[(c)] $X$ is strongly $\s$-compact iff it is a Hausdorff quotient space of a locally compact, $\s$-compact Hausdorff space.
\item[(d)] If $X$ is a strongly $\s$-compact, Hausdorff space, then $X$ is normal.
\item[(e)] If $X$ is strongly $\s$-compact and totally disconnected, then $X$ is strongly zero-dimensional.
\end{enumerate}
\end{prop}

{\bfseries Proof:} (a) If $\{ K_n \}$ is an increasing sequence of compacta with $K_n \subset K_{n+1}^{\circ}$ and $\bigcup_n K_n = X$ then
$A \cap K_n$ open in $K_n$ implies $A \cap K_n^{\circ}$ is open in $X$ and so $A = \bigcup_n A \cap K_n^{\circ}$ is open.

(b) Assume that $\{ K_n \}$ is a sequence of compacta in $X$ which determine the topology. Assume that $B \subset Y$ is such that
$B \cap q(K_n)$ is closed for every $n$. Then $q^{-1}( B \cap q(K_n)) \cap K_n = q^{-1}(B) \cap K_n$ is closed for every $n$. Hence, $q^{-1}(B)$ is
closed since the sequence determines the topology of $X$.  Since $q$ is a quotient map, $B$ is closed. Thus, $\{ q(K_n) \}$ determines the topology of
$Y$.

(c) If the sequence $\{ K_n \}$ determines the topology of $X$ then $X$ is a quotient of the disjoint union of the $K_n$'s. The converse follows from
(a) and (b).

(d) Let $Y$ be a locally compact, $\s$-compact, Hausdorff space and $q: Y \to X$ be a quotient map.  Let $F$ be the closed equivalence relation
$(q \times q)^{-1}(1_X)$ on $Y$. Let $B, \bar B$ be disjoint closed subsets of $X$. Let $Y_0 = q^{-1}(B \cup \bar B)$. Define
$L_0 : Y_0 \to [0,1]$ by $L_0(x) = 0$ for $x \in q^{-1}(B)$ and $= 1$ for  $x \in q^{-1}(\bar B)$. Thus, $L_0$ is a Lyapunov function for
$F \cap (Y_0 \times Y_0)$.  By Theorem  \ref{theocom2}, there exists $L : Y \to [0,1]$ an $F$ Lyapunov function which extends $L_0$. Since
$F$ is an equivalence relation, $L$ is constant on the $F$ equivalence classes and so factors to define a continuous map on $X$ which is
$0$ on $B$ and $1$ on $\bar B$.

(e)  Replacing $K_n$ by $\bigcup_{i \leq n} K_i$, if necessary, we can assume that the determining sequence of compacta is non-decreasing.
  We may also assume $\emptyset = K_0$.
Let $B, \bar B$ be disjoint closed subsets of $X$. Let $A_0 = \emptyset$. Assume inductively, that $A_n$ is a subset of $K_n \setminus \bar B$ clopen
with respect to $K_n$
 with $B \cap K_n  \subset A_n$  and with
with $A_{n-1} = A_n \cap K_{n-1}$. Observe that $ A_n \cup (K_{n+1} \cap B)$ and
$(K_n \setminus A_n) \cup (K_{n+1} \cap B)$ are disjoint compact sets in $X$. Since
$X$ is totally disconnected,  Proposition \ref{taupropzero} implies there is a clopen subset $U$ of $X$ which contains $A_n$ and is
disjoint from $(K_n \setminus A_n) \cup (K_{n+1} \cap B)$. Hence, $A_{n+1} = U \cap K_{n+1}$ is the required subset clopen in $K_{n+1}$.
The set $A = \bigcup A_n$ is disjoint from $B$ and since $A \cap K_n = A_n$ for all $n$, $A$ is clopen in $X$.

$\Box$ \vspace{.5cm}

\begin{lem}\label{metriclem} (a) For a pseudo-metric space $(X,d)$ the relation $Z_d = \{ (x,y) : d(x,y) = 0 \}$ is a closed equivalence
relation and $d$ induces on the quotient space $X/Z_d$ a metric $\tilde d$, so that $d = \tilde q^* \tilde d = \tilde d \circ ( \tilde q \times \tilde  q)$, with  $\tilde q : (X,d) \to (X/Z_d, \tilde d)$ the induced ``isometry''. The map $\tilde q$ is an open
map and a closed map and so is a quotient map.

(b) Let $E$ be a closed equivalence relation on a topological space $X$ and let $q : X \to X/E$ be the quotient map. A continuous pseudo-metric $d$
on $X$ with $E \subset Z_d$ induces a continuous pseudo-metric $\bar d$ on $X/E$ so that $d = q^* \bar d$. Conversely,
if $\bar d$ is a continuous pseudo-metric on $X/E$, then $d = q^* \bar d$ is a continuous pseudo-metric
on $X$ with $E \subset Z_d$.
\end{lem}

{\bfseries Proof:} (a) A subset $A$ is closed in $(X,d)$ iff $d(x,A) = 0$ implies $x \in A$.  Hence, a closed set is $Z_d$ saturated and $\tilde q(A)$ is
a closed set in $(X/Z_d, \tilde d)$. Taking complements we see that an $\tilde q$ is an open map as well.

(b) If $d$ is a continuous pseudo-metric on $X$ with $Z_d \subset E$ then $\tilde q : X \to X/Z_d$ factors through the projection $q$ to
define a map $h : X/E \to X/Z_d$ so that $\tilde q = h \circ q$. Since $q$ is a quotient map, $h$ is continuous. Hence, $\bar d = h^* \tilde d$ is
a continuous pseudo-metric on $X/E$ with $d = q^* \bar d$. The converse is obvious.

$\Box$ \vspace{.5cm}

 \begin{theo}\label{eqrel3} Let $f$ be a relation on a Tychonoff space $X$.
\begin{enumerate}
 \item[(a)] If $\G f$ is an equivalence relation, then $\G f$ is the smallest closed equivalence relation which contains $f$.

 \item[(b)] Assume that $\A_{\U_M} f$ is an equivalence relation. The relation $\A_{\U_M} f$ is the
 smallest closed equivalence relation $E$ containing $f$ such that the quotient space $X/E$ is completely Hausdorff.  In particular,
 $\G f = \A_{\U_M} f$ iff $\G f$ is an equivalence relation with the quotient space $X/\G f$  completely Hausdorff.

 The set $\{ \ell^f_d = s\ell^f_d : d \in \Gamma(\U_M) \}$ projects to the gage of the maximum uniformity with topology $\tau (X/\A_{\U_M} f)$.

 If $X$ is a locally compact, paracompact Hausdorff space and either $X$ is $\s$-compact or $f$ is a proper relation, then $1_X \cup \G f = \A_{\U_M} f$.
 The space $X/\A_{\U_M} f$ is a Hausdorff and normal and so $X/\A_{\U_M} f = \tau (X/\A_{\U_M} f)$.

 \item[(c)]Assume that $\CC_{\U_M} f$ is an equivalence relation. The relation $\CC_{\U_M} f$ is the
 smallest closed equivalence relation $E$ containing $f$ such that the quotient space $X/E$ is totally disconnected.

 The set $\{ m^f_d = sm^f_d : d \in \Gamma(\U_M) \}$ projects to  the gage of the maximum zero-dimensional uniformity with topology $\tau_0 (X/\A_{\U_M} f)$.

 If $X$ is a locally compact, paracompact Hausdorff space and either $X$ is $\s$-compact or $f$ is a proper relation, then
  $X/\CC_{\U_M} f$ is a Hausdorff, strongly zero-dimensional space and so $X/\CC_{\U_M} f = \tau (X/\CC_{\U_M} f) = \tau_0 (X/\CC_{\U_M} f)$.
\end{enumerate}
\end{theo}

 {\bfseries Proof:} (a) If $E$ is a closed equivalence relation which contains $f$ then, because it is transitive,
 $\G f \subset E$.  Because $\G f$ is a closed equivalence relation which contains $f$, it is the smallest such.

 (b) If $d \in \Gamma(\U_M)$, i.e.\ $d$ is a continuous pseudo-metric on $X$ then since $\A_{\U_M} f$ is
 reflexive and symmetric, Proposition \ref{prop00a} together with Proposition \ref{prop1a} implies that $\ell^f_d = s \ell^f_d$ is a pseudo-metric
  on $X$ with $\A_{\U_M} f \subset Z_{\ell^f_d} $. On the other hand, if $d$ is a continuous pseudo-metric on $X$ with $\A_{\U_M} f \subset Z_d$ then
 by Lemma \ref{domlem} $\ell^f_d = d$. By  Lemma \ref{metriclem} these are exactly the pullbacks via $q : X \to X/\A_{\U_M} f$ of continuous
 pseudo-metrics on the quotient space, i.e.\ the gage of the maximum uniformity with topology $\tau (X/\A_{\U_M} f)$.

 If $E$ is an equivalence relation then a Lyapunov function $L$ for $E$ is exactly a continuous real-valued function
 which is constant on each equivalence class, i.e.\
$L$ factors through the projection $q : X \to X/E$ to define a continuous real-valued function on $X/E$. Hence, $-L$ is a Lyapunov
function for $E$ as well. Hence, $X/E$ is completely Hausdorff iff $E = \bigcap \leq_L$ with $L$ varying over the
Lyapunov functions for $E$. So  Corollary \ref{cor15u} implies that $X/\A_{\U_M} f$ is completely Hausdorff.

On the other hand, if $E$ is a closed equivalence relation which contains $f$ and which has a completely Hausdorff
quotient, then $E = \bigcap \leq_L$ with $L$ varying over the Lyapunov functions for $E$.  Each such $L$ is a Lyapunov function
for $f$ and so is an $\A_{\U_M} f$ Lyapunov function by Corollary \ref{cor15u} again. Hence,
$\A_{\U_M} f \subset \ \leq_L$ for each such $L$. Hence, $\A_{\U_M} f \subset E$.

If $X$ is a locally compact,  $\s$-compact, Hausdorff space, then by Proposition \ref{propstronglysigma} the quotient $X/\A_{\U_M} f$ is
a strongly $\s$-compact Hausdorff space and so it normal.  As it is completely regular, it follows that $X/\A_{\U_M} f = \tau(X/\A_{\U_M} f)$.

If $X$ is a locally compact, paracompact Hausdorff space and $f$ is proper, then by Lemma \ref{lemcom1} $X/\A_{\U_M} f$ is a disjoint union of
clopen strongly $\s$-compact Hausdorff subspaces and so it is normal. Again, $X/\A_{\U_M} f = \tau(X/\A_{\U_M} f)$.

Finally, $1_X \cup \G f = 1_X \cup \A_{\U_M} f = \A_{\U_M} f$ by Corollary \ref{corcom2a}.

(c) If $d \in \Gamma(\U_M)$, then since $\CC_{\U_M} f$ is
 reflexive and symmetric, Proposition \ref{prop00a} together with Proposition \ref{prop1a} implies that $m^f_d = s m^f_d$ is a pseudo-ultrametric
  on $X$ with $\CC_{\U_M} f \subset Z_{m^f_d} $. On the other hand, if $d$ is a continuous pseudo-ultrametric on $X$ with $\CC_{\U_M} f \subset Z_d$ then
 by Lemma \ref{domlem} $m^f_d = d$. By  Lemma \ref{metriclem} these are exactly the pullbacks via $q : X \to X/\CC_{\U_M} f$ of continuous
 pseudo-ultrametrics on the quotient space, i.e.\ the gage of the maximum zero-dimensional uniformity with topology $\tau_0 (X/\CC_{\U_M} f)$.

Assume that $(x,y) \not\in \CC_{\U_M} f$. There exists a continuous pseudo-metric $d$ on $X$ such that $m_d^f(x,y) = \ep > 0$.
Since  $m_d^f $ is a pseudo-ultrametric,
 $V^d_{\ep}(x)$ is a clopen set which contains $x$ but not $y$. Furthermore,  $V^d_{\ep}(x)$ is $\CC_{\U_M} f$
saturated.  Hence, if $q : X \to X/\CC_{\U_M} f$ is the projection, $q(\CC_{\U_M} f)$ is a clopen subset of $X/\CC_{\U_M} f$ which contains
$q(x)$ but not $q(y)$.  It follows that $X/\CC_{\U_M} f$ is totally disconnected.

On the other hand, let $E$ be a closed equivalence relation which contains $f$ and which has a totally disconnected
quotient with quotient map $q : X \to X/E$. If $(x,y) \not\in E$ then there exists a clopen set $A_1 \subset X/E$ with $q(x) \in A_1$ and
$q(y) \in B_1 = (X/E) \setminus A_1$. So $A = q^{-1}(A_1)$ and $B = q^{-1}(B_1)$ form a clopen partition of $X$.  Let $U = (A \times A) \cup (B \times B)$.
This is a clopen equivalence relation on $X$ with $f \subset E \subset U$.  It follows that if $[a,b] \in f^{\times n}$ is an
$(x,z), U$ chain, then with $b_0 = x, a_{n+1} = z$, $(a_i,b_i) \in f \subset U$ for $i = 1, \dots, n$ and $(b_i,a_{i+1}) \in U$ for $i = 0, \dots, n$.
Since $U$ is an equivalence relation $z \in U(x) = A$ and so $z \not= y$. Hence, $(x,y) \not\in \CC_{\U_M} f$. Contrapositively, $ \CC_{\U_M} f \subset E$.

If $X$ is a locally compact,  $\s$-compact, Hausdorff space, then by Proposition \ref{propstronglysigma} the quotient $X/\CC_{\U_M} f$ is
a strongly $\s$-compact, totally disconnected  space and so it strongly zero-dimensional.
As it is completely regular, it follows that $X/\CC_{\U_M} f = \tau(X/\CC_{\U_M} f)$.
As it is zero-dimensional, it follows that $X/\CC_{\U_M} f = \tau_0(X/\CC_{\U_M} f)$.

If $X$ is a locally compact, paracompact Hausdorff space and $f$ is proper, then by Lemma \ref{lemcom1} $X/\CC_{\U_M} f$ is a disjoint union of
clopen strongly $\s$-compact totally disconnected subspaces and so it is strongly zero-dimensional. Again,
$X/\CC_{\U_M} f = \tau (X/\CC_{\U_M} f) = \tau_0 (X/\CC_{\U_M} f)$.

 $\Box$ \vspace{.5cm}

 \begin{cor}\label{corquasi} For a Tychonoff space $1_X$, $\CC_{\U_M} 1_X$ is a closed equivalence relation with  equivalence classes the
 quasi-components of $X$. \end{cor}

 {\bfseries Proof:} Since $1_X$ is symmetric, $\CC_{\U_M} 1_X$ is a closed equivalence relation with a totally disconnected quotient via
 the quotient map $q : X \to X/ \CC_{\U_M} 1_X$ by Theorem \ref{eqrel3}. So if $q(x) \not= q(y)$ there is a clopen set $A \subset X/ \CC_{\U_M} 1_X$
 with $q(x) \in A$ and $q(y) \not\in A$. Since $U = q^{-1}(A)$ is clopen with $x \in U$ and $y \not\in U$, $x$ and $y$ lie in separate quasi-components.
 On the other hand, if $U_1 $ is a clopen subset of $X$ with $x \in X$ and $y \in U_2 = X \setminus U_1$ then $ E = (U_1 \times U_1) \cup (U_2 \times U_2)$
 is a clopen equivalence relation on $X$ and so $E \in \U_M$. If $[a,b] \in 1_X^{\times n}$ defining an $xz, E$ chain then $z \in U_1$ and so $z \not= y$.
 Hence, $(x,y) \not\in \CC_{\U_M} 1_X$.

 $\Box$ \vspace{.5cm}

 \begin{lem}\label{lemid} If $f$ is a relation on a Hausdorff uniform space $(X,\U)$, then $\A_{\U} (1_X \cup f) = 1_X \cup A_{\U} f$.\end{lem}

 {\bfseries Proof:} Clearly $1_X \cup A_{\U} f \subset \A_{\U} (1_X \cup f)$. If $(x,y) \not\in 1_X \cup A_{\U} f $ then because $(X,\U)$ is
 Hausdorff there exists $d_1 \in \Gamma(\U)$ such that $d_1(x,y) > 0$. Also, there exists  $d_2 \in \Gamma(\U)$ such that
 $\ell^{f}_{d_2}(x,y) > 0$. Hence, $d = d_1 + d_2 \in \Gamma(\U)$ with $(x,y) \not\in Z_d \cup \A_d f$. By (\ref{eq38})
 $(x,y) \not\in \A_d (1_X \cup f)$ and so is not in $\A_{\U} (1_X \cup f)$.

 $\Box$ \vspace{.5cm}

  \begin{cor}\label{eqrel3cor} Let $f$ be a relation on a Hausdorff uniform  space $(X,\U)$.

  The closed equivalence relations $1_X \cup (\A_{\U} f \cap \A_{\U} f^{-1})$ and $1_X \cup (\CC_{\U} f \cap \CC_{\U} f^{-1})$
  have completely Hausdorff quotients. On $|\CC_{\U} f|$ the equivalence relation $\CC_{\U} f \cap \CC_{\U} f^{-1}$ has a totally disconnected quotient.

   If $X$ is a locally compact, $\s$-compact Hausdorff space, then the quotients are
   Hausdorff and normal and $|\CC_{\U}|/[\CC_{\U} f \cap \CC_{\U} f^{-1}]$ is Hausdorff and strongly zero-dimensional.

   \end{cor}

   {\bfseries Proof:} $X$ is Tychonoff and so we can apply Lemma \ref{lemid},  (\ref{eq20ac}) together with
   monotonicity and idempotence of the operator $\A_{\U}$ to get
   \begin{align}\label{eqre13coreq1}
   \begin{split}
   1_X \cup (\A_{\U} f \cap \A_{\U} f^{-1}) \ \subset \quad &\A_{\U_M}[1_X \cup (\A_{\U} f \cap \A_{\U} f^{-1}]\\
   \subset  \ \A_{\U}[1_X \cup (\A_{\U} f \cap \A_{\U} f^{-1}] \
   =  \ &1_X \cup (\A_{\U}\A_{\U} f \cap \A_{\U}\A_{\U} f^{-1})\\ =  \ 1_X \cup (\A_{\U} f &\cap \A_{\U} f^{-1}).
   \end{split}
   \end{align}

    Thus, $E =  1_X \cup (\A_{\U} f \cap \A_{\U} f^{-1}) $ is a closed equivalence relation with $\A_{\U} E = E$. Similarly,
     $E =  1_X \cup (\CC_{\U} f \cap \CC_{\U} f^{-1})$ is a closed equivalence relation with $\A_{\U} E = E$.   By Theorem \ref{eqrel3} (b)
     each has a completely Hausdorff quotient and a normal Hausdorff quotient when $X$ is locally compact and $\s$-compact.

     Similarly, if $E = \CC_{\U} f \cap \CC_{\U} f^{-1}$ then $\CC_{\U_M} E = E$. Since $E \subset |\CC_{\U} f | \times |\CC_{\U} f |$,
     we can apply Theorem \ref{eqrel3} (c), replacing $X$ by $|\CC_{\U} f |$ on which $E$ is a closed equivalence relation. We obtain that
     the quotient is totally disconnected and is Hausdorff and strongly zero-dimensional when $X$ is locally compact and $\s$-compact.

 $\Box$ \vspace{.5cm}

 \begin{prop}\label{eqre14prop} Let $E$ be a closed equivalence relation on a Tychonoff space $X$.

 \begin{enumerate}
 \item[(a)] The relation $E$ is usc iff the quotient map $q : X \to X/E$ is a closed map.

 \item[(b)]  If $E$ is usc and $X$ is normal, then $X/E$ is a Hausdorff normal space.

\item[(c)] If $E$ is cusc, and $X$ is locally compact, then $X/E$ is locally compact.

\item[(d)]  If $E$ is cusc, and $X$ is second countable, then $X/E$ is second countable.
 \end{enumerate}
 \end{prop}

 {\bfseries Proof:} (a) If $A \subset X$ then $q^{-1}(q(A)) = E(A)$. To say that $E$ is usc is to say that $E(A)$ is closed whenever
 $A$ is.  To say that $q$ is closed is to say that $q^{-1}(q(A))$ is closed whenever $A$ is. So the equivalence is clear.

 (b) If $A_0, A_1 \subset X$ are disjoint closed sets with $E(A_0) = A_0$ and $E(A_1) = A_1$ then let $X_0 = A_0 \cup A_1$, $L_0(x) = 0$
 for $x \in A_0$ and $ = 1$ for $x \in A_1$. Thus, $L_0$ is a Lyapunov function for $E| (X_0 \times X_0)$ and so by Theorem \ref{theo16ae}
 extends to a Lyapunov function $L$ for $E$. This implies normality of $X/E$.

 (c), (d) We choose a basis $\B$ for $X$ which is closed under finite unions. Let $\tilde \B = \B = \{ E^*U : U \in \B \}$. Since
 $E$ is usc, each member of $\tilde \B$ is an $E$ saturated open set. If $x \in X$ and
 $V$ is open with $E(x) \subset V$ then since $E(x)$ is compact, there exists $U \in \B$ such that $E(x) \subset U \subset V$ and so
 $E(x) \subset E^*(U) \subset U \subset V$. Thus, $\B_E = \{ q(V) : V \in \tilde \B \}$ is a basis for $X/E$.

 For (c) we can choose $\B$ so that every member has compact closure. Since $q(E^*(U))\subset q(\ol{U})$ it follows that each $q(V)$ for
 $V \in \B_E$ has compact closure in $X/E$ and so $X/E$ is locally compact.

 For (d) choose $\B$ countable.  Then $\B_E$ is a countable basis for $X/E$.

 $\Box$ \vspace{.5cm}

 A second countable space which admits a complete metric is called a \emph{Polish space}\index{Polish space}\index{space!Polish}. Any
 $G_{\d}$ subset of a Polish space is a Polish space. A locally compact, second countable space is $\s$-compact and Polish.

 \begin{exes}\label{eqrel4} \begin{itemize}
 \item[(a)] There exists a homeomorphism $f$ on
 a separable metric space $X$ such that $\G f$ is an equivalence relation such that the quotient space
 $X/\G f$ is not Hausdorff and so $\G f$ is a proper subset of $\A_{\U_M} f$.

 \item[(b)] There exists a homeomorphism $f$ on
 a locally compact space $X$ such that $\G f$ is an equivalence relation such that the quotient space
 $X/\G f$ is not Hausdorff and so $\G
 f$ is a proper subset of $\A_{\U_M} f$.

 \item[(c)] There exists a homeomorphism $f$ on a Polish space $X$ with metric $d$, such that $\G f = \CC_{\U_M} f = \CC_d f$
  is an equivalence relation with a totally disconnected quotient which is not regular.

   \item[(d)] There exists a homeomorphism $f$ on a locally compact space, such that $\G f = \CC_{\U_M} f$
  is an equivalence relation with a totally disconnected quotient which is not regular and so is not zero-dimensional.

  \item[(e)] There exists a homeomorphism $f$ on a locally compact, $\s$-compact, metrizable space, such that $\G f = \CC_{\U_M} f$
  is an equivalence relation with a Hausdorff, strongly zero-dimensional quotient which is not first countable and so is not metrizable.
  \end{itemize}

 \end{exes}

 {\bfseries Proof:} (a) The following is a variation of the example in Problem 3J of \cite{GJ}.

 Let $g$ be a topologically transitive homeomorphism on a compact metric space $Y$ with a Cantor set
 $C \subset Y$ of fixed points.  Such maps can be constructed with $Y$ the torus or the Cantor set itself.

 Let $D$ be a countable dense subset of $C$ and $J = C \setminus D$ so that $J$ is a dense $G_{\delta}$ subset of $C$.
 Choose $e \in C$. For the homeomorphism $g \times g$ on $Y \times Y$, the compact set $Y \times \{ e \}$ and
 the $G_{\delta}$ set $J \times Y$ are $g \times g$ invariant. The restriction of $g \times g$ to $Y \times \{ e \}$
 is topologically transitive with $C\times \{ e \}$ a set of fixed points. For each $j \in J$, the restriction of $g \times g$ to
 $\{ j \} \times Y$ is topologically transitive with $(j,e)$ a fixed point. Let $X_0 =  Y \times \{ e \} \cup J \times Y$
 and $f_0$ be the restriction of $g \times g$ to this invariant set.

 Mapping $Y$ to $e$ we obtain a
 retraction $\pi : J \times Y \to J \times \{ e \} $.  By extending the definition of $\pi$ to
 be the identity on $Y \times \{ e \}$, we define the continuous retraction $\pi : X_0 \to Y \times \{ e \}$.

 Let $E_1$ denote the closed equivalence relation
 $$\pi^{-1} \circ \pi = (\pi \times \pi)^{-1}(1_{(Y \times \{ e \}) \times (Y \times \{ e \})}), $$
 Let $E_2 = 1_{J \times Y} \cup (Y \times \{ e \}) \times (Y \times \{ e \})$ which is also a closed equivalence
 relation.  Hence, $E_0 = E_1 \cup E_2$ is a closed, reflexive, symmetric relation on $X_0$. It is not, however, transitive.

 Let $X = X_0 \setminus (J \times \{ e \})$. Because we are removing a set of fixed points, $f_0$ restricts to a homeomorphism
 $f$ on $X$. Let $E$ denote the restriction $E_0 \cap (X \times X)$,
  a closed, reflexive, symmetric relation on $X$. We show that it is also transitive.

  Let $x,y \in X$.

 \begin{itemize}
 \item   $(x,y) \in E_1 \cap (X \times X) \setminus 1_X$ iff $x_1 = y_1 \in J$ and $x_2, y_2 \in Y \setminus \{ e \}$
 with $x_2 \not= y_2$.

 \item  $(x,y) \in E_2 \cap (X \times X) \setminus 1_X$ iff $x_2 = y_2 = e$ and $x_1,y_1 \in Y \setminus J$
 with $x_1 \not= y_1$.
 \end{itemize}

Assume $(x,y),(y,z) \in E$ if  $x = y$ (or $y = z$) then  $(x,z) = (y,z)$ (resp.\   $(x,z) = (x,y)$) and so $(x,z) \in E$.
So we may assume  $(x,y),(y,z) \in E \setminus 1_X$.

If $(x,y) \in E_1 \cap (X \times X) \setminus 1_X$ then $y_2 \not= e$ and so $ (y,z) \not\in E_2 \cap (X \times X) \setminus 1_X$.
Hence, $ (y,z) \in E_1\cap (X \times X) \setminus 1_X$ and so $(x,z) \in E_1\cap (X \times X) \subset E$.
If $(x,y) \in E_2 \cap (X \times X) \setminus 1_X$ then $y_2 = e$ and so, as before,
$ (y,z) \in E_2 \cap (X \times X) \setminus 1_X$ and  $(x,z) \in E_2\cap (X \times X) \subset E$. Thus, $E$ is transitive.

 From the invariance and transitivity results, it is clear that $f \subset E$ and $E \subset \G f$.
 Since $E$ is a closed, transitive relation which contains $f$, it contains $\G f$. Thus, $\G f = E$.

 Now consider the quotient space of $X$ by the equivalence relation $E$, with quotient map $q : X \to X/E$.
 We will see that $X/E$ is not Hausdorff even though $E$ is a closed relation.  In particular, this
 implies that $q \times q : (X \times X) \to (X/E) \times (X/E)$ is not a quotient map since $1_{X/E}$ not closed,
 because $X/E$ is not Hausdorff, but its pre-image is the closed set $E$.

 The set $(Y \setminus J) \times \{ e \}$ is mapped by $q$ to a single point which we will call $e^*$.  Let $G$ be a
 nonempty open subset of $X/E$. Since $(Y \setminus J) \times \{ e \}$ is not open in $X$, it follows that
 the $E$ saturated open set $U = q^{-1}(G) \cap (J \times (Y \setminus \{ e \})$ is nonempty.  The
 projection $\pi_1 : J \times Y \to J $ is an open map and so
 the image $\pi_1(U)$ is a nonempty open subset of $J$. Since $D = C \setminus J$ is dense in
 $C$, it follows that the closure in $C $ of $\pi_1(U)$ meets $D$.
 That is, there exists a sequence $\{ (j_n,y_n) \in U  \}$ such that
 $j_n \to d$ with $d \in D$. Since $U$ is $E$ saturated, we can vary $y_n$ arbitrarily in
 $Y \setminus \{ e \}$.  Because $g$ was topologically transitive, $e$ is not an isolated point in $Y$ and so
 we can choose $y_n \in Y \setminus \{ e \}$ converging to $e$. It follows that $\overline{U}$ contains
 the point $(d,e) \in (Y \setminus J) \times \{ e \}$. Hence, $e^{*} \in q(\overline{U}) \subset \overline{G}$.
 It follows that every neighborhood of $e^*$ is dense in $X/E$.

 Any Lyapunov function $L$ for $f$ is a Lyapunov
 function for $\G f = E$ and so
 factors through $q$ to yield a continuous real-valued function $\tilde L : X/E \to \R$. If
 $t \not= \tilde L(e^*)$ then we can choose disjoint open sets $U_1, U_2 \subset \R$ with $\tilde L(e^*) \in U_1, t \in U_2$.
 Thus, $e^*$ is in the open set $(\tilde L)^{-1}(U_1)$ which is disjoint from the open set $(\tilde L)^{-1}(U_2)$. Since
 $(\tilde L)^{-1}(U_1)$ is dense,  $(\tilde L)^{-1}(U_2)$ is empty. So $t$ is not in the image of $\tilde L$.
 Thus, $L = \tilde L \circ q$ is constant at the value $\tilde L(e^*)$.

 Thus, the only Lyapunov functions for $f$ are constant functions. It follows from Corollary \ref{cor15u} that
 $1_X \cup \A_{\U_M} f = X \times X$.  Since there are no isolated points in $X$, $X \times X \setminus 1_X$ is
 dense in $X \times X$. Since $\A_{\U_M} f$ is a closed relation, it follows that $\A_{\U_M} f = X \times X$. On the
 other hand, $E = \G f$ is a proper subset of $X \times X$.

 While $X_0$ is a $G_{\delta}$ subset of the compact metric space $Y \times Y$, $X$ is not.
 We do not know of examples like this with $X$ a Polish space. In particular, we do not know of an example of a
 closed equivalence relation on a Polish space with a non-Hausdorff quotient.

 (b), (c), (d), (e): Let $\om$ and $\Omega$ denote the first countable and first uncountable ordinal respectively.
 In particular, $\om$ is the set of non-negative integers. The ordered set $\R_+ = \om \times [0,1)$ with the lexicographical
 ordering is order-isomorphic with the half-open interval $[0,\infty)$ by $(n,t) \mapsto n+t$. With the order topology this bijection is
 a homeomorphism. The ordered set $L = \Omega \times [0,1)$ with the lexicographical ordering can be similarly equipped with the order topology
 to obtain the \emph{Long Line}.\index{Long Line} It is a non-paracompact, locally compact space and for every $\a \in \Omega$ the interval $[(0,0), (\a,0)]$ is
 order-isomorphic and thus homeomorphic with the unit interval. We double each example.  Let $\tilde \R = \R_+ \times \{ +,- \}$ with
 each $(n,0,+)$ identified with $(n,0,-)$. We identify $\om \subset \tilde \R$ by $n \mapsto (n,0,\pm)$.
 Let $\tilde L = L \times \{ +, - \}$ with each $(\a,0,+)$ identified with $(\a,0,-)$. We identify $\Omega \subset \tilde L$ by
 $\a \mapsto (\a,0,\pm)$.

 Let $\om^* = \om + 1 = \om \cup \{ \om \}$ and $\Omega^* = \Omega + 1 = \Omega \cup \{ \Omega \}$. These are the one-point compactifications of
 $\om$ and $\Omega$, respectively.  Similarly, let $\tilde \R^*$ and $\tilde L^*$ denote the one-point compactifications with points $\om$ and $\Omega$
 the respective points at infinity. The product $\Omega^* \times \om^*$ is compact and removing the point $(\Omega,\om)$ we obtain the locally
 compact  \emph{Tychonoff Plank} \index{Tychonoff Plank}$T$, see \cite{K} Example 4F.
 As described there, the Tychonoff Plank is not normal as the closed subsets
 $\Omega \times \{ \om \}$ and $\{ \Omega \} \times \om$ cannot be separated by open sets.

 On the unit interval $[0,1]$ let $u_+(t) = \sqrt{t}$ and $u_-(t) = t^2$. Each is a homeomorphism with fixed points $0$ and $1$. Observe that
 $u_+(t) > t$ and $u_-(t) < t$ for all $t \in (0,1)$. Thus, for every $t \in (0,1)$ the bi-infinite orbit sequence $\{ (u_-)^n(t) \}$ converges to
 $0$ and $n \to \infty$ and to $1$ as $n \to -\infty$. Since $u_+ = (u_-)^{-1}$ the reverse is true for the $u_+$ orbit sequences.
 On $\tilde \R$ define the homeomorphism $g$ by $g(n,t,\pm) = (n,u_{\pm}(t),\pm)$ and on $\tilde L$
 define the homeomorphism $G$ by $G(\a,t,\pm) = (\a,u_{\pm}(t),\pm)$. Observe that $\om \subset \tilde \R$ is the set of fixed points of $g$ and
 $\Omega \subset \tilde L$ is the set of fixed points of $G$. Notice that $\G g = \tilde \R \times \tilde \R$ and $\G G = \tilde L \times \tilde L$.

 We use these to construct our remaining examples.

 (b) Let $X$ equal $T \cup \tilde L \cup \tilde \R$ with $(\a, \om) \in T$ identified with $\a \in \tilde L$ for all $\a \in \Omega$ and
 with $(\Omega,n) \in T$ identified with $n \in \tilde \R$ for all $n \in \om$. Thus, $X$ is a locally compact, non-paracompact, Hausdorff space.
 The homeomorphism $f$ is
 the homeomorphism induced from $1_T \cup G \cup g$ via these identifications. Thus, $T$ is the set of fixed points of $f$.
 Clearly, $\G f$ is the equivalence relation $1_T \cup (\tilde L \times \tilde L) \cup (\tilde R \times \tilde R)$.
 The quotient space $X/\G f$ is the quotient space of
 the Tychonoff plank $T$ with the two closed subsets  $\Omega \times \{ \om \}$ and $\{ \Omega \} \times \om$  each smashed to a point.
 Since the closed sets cannot be separated in $T$, the quotient space is not Hausdorff.

 (c) Let $C \subset [0,1]$ be the Cantor Set and let $A = \{ a_1, a_2, \dots \}$ with $\{ a_k\}$
  a decreasing sequence in $C$ which converges to $0$. Let $\hat C$ be $C$ with the topology obtained by including
  $C \setminus A$ as an open set.  The new topology is $\{ U_1 \cup (U_2 \setminus A) : U_1, U_2 $ open in $C \}$.
  Thus, if $x \in C$ with $x \not= 0$ then a set is a neighborhood of $x$ iff it contains a $C$ open set $U$ with
  $x \in U$. A set is a neighborhood of $0$ iff it contains $U \setminus A$ with $U$ a $C$ open set such that $0 \in U$.
  Since the topology is finer than the original topology of $C$, the space $\hat C$ is completely Hausdorff. Note that it
  has a countable base. However, it is not regular. The closure of any neighborhood of $0$ meets $A$ and so there is
  no closed neighborhood of $0$ contained in the $\hat C$ open set $C \setminus A$.

   Observe that if $E$ is a closed equivalence relation on a Tychonoff space $X$ then the quotient $X/E$ is $T_1$ and
  so is Hausdorff if it is regular. If $X$ is a separable metric space, or, more generally, any Lindel\"{o}f space then
  the quotient is Lindel\"{o}f.  Since a regular, Lindel\"{o}f space is normal (see \cite{K} Lemma 4.1), it follows
  that if $E$ is a closed equivalence relation on a separable metric space $X$, then the quotient is Hausdorff and
  normal, and so completely regular, if it is regular.

  Let $X_0 = C \times \tilde \R$ with $f_0 = 1_C \times g$ and let $p_0 : X_0 \to C$ be the first coordinate projection.
  Clearly, $\G f_0 = p_0^{-1} \circ p_0$.
  That is, $\G f_0$ is a closed equivalence relation with equivalence classes the fibers of $p_0$. $X_0$ is a locally compact,
  metrizable space.

  Now let $Z_k =  \{ (n,t,\pm) \in \tilde \R : n < k \}$.

Let $X$ be the $G_{\delta}$ invariant subset $X_0 \setminus (\bigcup_{k = 1}^{\infty} \ \{ a_k \} \times Z_k)$ and
  let $f$ be the restriction of $f_0$ to $X$. Again $\G f  = p^{-1} \circ p$ where $p$ is the restriction of $p_0$.
  Notice that $p^{-1}(A)$ is a closed subset of $X$.  It easily follows that $p$ induces a homeomorphism of the quotient space
  $X/\G f$ onto $\hat C$. Thus, the quotient is not regular although it is completely Hausdorff.

   Notice that since $C$ is totally disconnected, it follows that for any metric $d$ on $X$, $\CC_d f = \G f$. Hence, $\CC_{\U_M} f = \G f$.
 Hence, for any uniformity $\U$ compatible with the topology on $X$, the inclusions
 $\G f \subset \A_{\U} f \subset \A_d f \subset \CC_d f$ and  $\G f \subset \CC_{\U} f \subset \CC_d f$
  imply that they are all equal. By Theorem \ref{eqrel3} the quotient space is totally disconnected.

  (d) We return to the Tychonoff Plank. Let $X$ equal $T \cup \tilde L $ with $(\a, \om) \in T$ identified with $\a \in \tilde L$ for all $\a \in \Omega$.
  Again $X$ is a locally compact, non-paracompact, Hausdorff space. In addition, it is zero-dimensional
  but not strongly zero-dimensional since it is not normal.

 The homeomorphism $f$ is
 the homeomorphism induced from $1_T \cup G $ via these identifications. Again $T$ is the set of fixed points of $f$.
 Clearly, $\G f$ is the equivalence relation $1_T \cup (\tilde L \times \tilde L)$.
 The quotient space $X/\G f$ is the quotient space of
 the Tychonoff plank $T$ with the closed subset  $\Omega \times \{ \om \}$  smashed to a point $e$.
Because $X$ is locally compact, it is completely regular. It follows that the quotient space $X/\G f$ is completely Hausdorff. However,
the point $e$ cannot be separated from the closed set $\{ \Omega \} \times \om$ and so the quotient is not regular.

Since $T$ is zero-dimensional we have that $\G f = \CC_{\U_M} f$. The quotient is totally disconnected but not zero-dimensional since it is not regular.

Notice that if we extend $f$ to the one-point compactification $X^*$ of $X$, by adjoining the point $(\Omega,\om)$ we obtain a homeomorphism $f^*$
$\G f^* = 1_T \cup (\tilde L^* \times \tilde L^*)$. The quotient space $X^*/\G f^*$ is a compact, Hausdorff space and the inclusion
$X \to X^*$ induces a continuous bijection $X/\G f \to X^*/\G f^*$ which is not a homeomorphism because $\{ \Omega \} \times \om$ is not closed
in $X^*/\G f^*$.

(e) Let $X = (\tilde \R \times \{ 0 \}) \cup \N \times \{ 1/k : k \in \N \}$. Let $f = g \times 1_{\{0 \}} \cup 1_{\N \times \{ 1/k : k \in \N \} }$.
Clearly, $\G f = \CC_{\U_M} f$ with quotient obtained by smashing $\tilde \R \times \{ 0 \}$ to a point $e$. The point $e$ does not have a countable
neighborhood base.  If $\{ U_n : n \in \N \}$ is a sequence of neighborhoods of $ \tilde \R \times \{ 0 \}$ in $X$ then for every $n \in \N$ there
exists $k_n \in \N$ such that $(n,1/k_n) \in U_n$. The set $\{ (n,1/k_n) : n \in \N \}$ is closed and disjoint from $\tilde \R \times \{ 0 \}$, but
meets every $U_n$.

For cases (b),(d) and (e) the relations $\G f$ are usc.  In general, if $A, B$ are disjoint closed subsets of $X$ then $(A \times A) \cup 1_X$ and
$(A \times A) \cup (B \times B) \cup 1_X$ are closed, usc equivalence relations.

 $\Box$ \vspace{.5cm}

Recall that a relation $f$ on $X$ is \emph{surjective}\index{relation!surjective}\index{surjective relation}
 if $Dom(f) = Dom(f^{-1}) = X$, i.e.\ $f(X) = f^{-1}(X) = X$.

 \begin{df}\label{trans1} A relation $f$ on a uniform space $(X,\U)$ is called $\U$ chain transitive when it is a surjective relation
 such that $\CC_{\U} f = X \times X$. \end{df}\index{relation!$\U$ chain transitive}\index{ch@$\U$ chain transitive}
\vspace{.5cm}

\begin{prop}\label{trans2} Let $f$ be a relation on a uniform space $(X,\U)$.
\begin{enumerate}
\item[(a)] If $f$ is $\U$ chain transitive then $f^{-1}$ is $\U$ chain transitive.
\item[(b)] If $f$ is a proper relation with $\CC_{\U} f = X \times X$ then $f$ is a surjective relation.
\item[(c)] If $f$ is a surjective relation then $f$ is $\U$ chain transitive iff for every $d \in \Gamma(\U)$
$M^f_d(x,y) = 0$ for all $x,y \in X$.
\item[(d)] If $g$ is a surjective relation on a uniform space $(Y,\V)$  and $h : X \to Y$ is a uniformly continuous surjective map
which maps $f$ to $g$, then  $g$ is $\V$ chain transitive if $f$ is $\U$ chain transitive.
\end{enumerate}
\end{prop}

{\bfseries Proof:} (a) The inverse of a surjective relation is clearly surjective and $\CC_{\U}  (f^{-1}) = (\CC_{\U}  f)^{-1}$.

(b) By Proposition \ref{prop16aj} and Proposition \ref{prop16ak} $Dom(f) = Dom(\CC_{\U} f) = X$ and  $Dom(f^{-1}) = Dom(\CC_{\U} f^{-1}) = X$.

(c) Since $m^f_d \leq M^f_d$ it is clear that $M^f_d(x,y) = 0$ implies $m^f_d(x,y) = 0$. So if for every $d \in \Gamma(\U)$
$M^f_d(x,y) = 0$ for all $x,y \in X$, then $\CC_{\U} f = X \times X$.

For the converse we cannot apply Proposition \ref{prop00uu} because we are not assuming that $f$ is usc. Given $d \in \Gamma(\U)$, $\ep > 0$ and
$x,y \in X$ there exists $z \in f(x)$ since $f$ is surjective. Because $(z,y) \in \CC_{\U} f$ there exists $[a,b] \in f^{\times n}$ with
the $zy$ chain-bound of $[a,b]$ less than $\ep$. Now define $[a,b]' \in f^{\times n+1}$ with $(x,z) = (a'_1,b'_1)$
and $(a'_i,b'_i) = (a_{i-1},b_{i-1})$ for
$i = 2, \dots, n+1$. Since the $xy$ chain-bound of $[a,b]'$ equals the $zy$ chain-bound of $[a,b]$ and $x = a'_1$ it follows that
$M^f_d(a,y) < \ep$.

(d) $Y \times Y = (h \times h)(X \times X) = (h \times h)(\CC_{\U} f) \subset \CC_{\V} g$ by Proposition \ref{prop6ub}.

$\Box$ \vspace{.5cm}

 \begin{df}\label{trans3} A relation $f$ on a uniform space $(X,\U)$ is called $\U$ chain mixing when it is a surjective relation and
 for every $d \in \Gamma(\U), \ep > 0, x,y \in X$ there exists a positive integer $N$ so that for all $n \geq N$
 there exists $[a,b] \in f^{\times n}$ with $a_1 = x$ and with the $xy$ chain-bound of $[a,b]$ with respect to $d$ less than $\ep$.\end{df}
 \index{relation!$\U$ chain mixing}\index{ch@$\U$ chain mixing}

 That is, for any  $d, \ep$ and $x,y$ for sufficiently large $n$ there is a chain of length $n$ from $x$ to $y$ with initial position $x$.

 Thus, $f$ is a $\U$ chain transitive relation iff $X \times X = \bigcup_{n=1}^{\infty} (V_{\ep}^d \circ f)^n$ for all $d \in \Gamma(\U)$ and
 $\ep > 0$. The relation $f$ is chain mixing iff $X \times X = \bigcup_{n=1}^{\infty} \bigcap_{i=n}^{\infty} (V_{\ep}^d \circ f)^i$
 for all $d \in \Gamma(\U)$ and  $\ep > 0$.

 For a positive integer $k$ the \emph{$k$-cycle} is the translation bijection $s(n) = n+1$ on the cyclic group $\Z_k = \Z/ k \Z$.

 \begin{theo}\label{trans4} Let $f$ be a $\U$ chain transitive relation on a uniform space $(X,\U)$.
 \begin{enumerate}
 \item[(a)]
 The following conditions are equivalent
 \begin{itemize}
  \item[(i)] The relation $f$ is $\U$ chain mixing.
 \item[(ii)] The relation $f \times f$ on $(X \times X, \U \times \U)$ is $\U$ chain mixing.
 \item[(iii)] The relation $f \times f$ on $(X \times X , \U \times \U)$ is $\U$ chain transitive.
  \item[(iv)] There does not exist for any integer $k > 1$ a uniformly continuous surjection from $X $ to $\Z_k$ which maps $f$ to $s$.
 \end{itemize}

 \item[(b)] If $f$ is $\U$ chain mixing then $f^{-1}$ is $\U$ chain mixing.

 \item[(c)] If for every positive integer $k$, the relation $f^k$ is $\U$ chain transitive, then $f$ is $\U$ chain mixing.
 Conversely, if $f$ is a uniformly continuous mapping which is $\U$ chain mixing, then
 for every positive integer $k$, the mapping $f^k$ is $\U$ chain mixing.
 \end{enumerate}
\end{theo}

{\bfseries Proof:} (a) (i) $\Leftrightarrow$ (ii): Easy to check.

(ii) $\Rightarrow$ (iii): A  chain mixing relation is chain transitive.

 If $h$ is uniformly continuous mapping $f $ onto a
surjective relation $g$ then $h \times h$ maps $f \times f$ to  $g \times g$ and $h$ maps $f^n$ to $g^n$.
Observe that with $k > 1$ $s \times s$ on $\Z_k \times \Z_k$ is not chain transitive since it is the disjoint union of $k$ separate periodic orbits.
Furthermore, $s^k = 1_{\Z_k}$ and so $s^k$ is not chain transitive. So Proposition \ref{trans2} (d) implies (iii) $\Rightarrow$ (iv) and
and if $f^k$ is chain transitive for all positive $k$ then (iv) holds.

We prove the contrapositive of  (iv) $\Rightarrow$ (i) following Exercise 8.22 of  \cite{A93}.   See also \cite{RW}. Assume $f$ is $\U$ chain
transitive but not $\U$ chain mixing. With $d \in \Gamma(\U)$ and $\ep > 0$ fixed we define for $x,y \in X$ the set
of positive integers $N(x,y)$ by $n \in N(x,y)$ iff there exists $[a,b] \in f^{\times n}$ with $a_1 = x$ and with
the $xy$ chain-bound of $[a,b]$ with respect to $d$ less than $\ep$. Since $f$ is assumed to be $\U$ chain transitive, Proposition
\ref{trans2} (c) implies that $N(x,y)$ is non-empty for every pair $x,y$.  With $A,B$ nonempty subsets of $\N$ we let $A + B$ denote
$\{ a + b : a \in A, b \in B \}$. By concatenating chains we observe that for $x,y,z \in X$
\begin{equation}\label{transeq1}
N(x,y) + N(y,z) \ \subset \ N(x,z). \hspace{4cm}
\end{equation}
In particular, $N(x,x)$ is an additive sub-semigroup of $\N$.  Let $k(x)$ be the greatest common divisor of the elements of $N(x,x)$.
We will need the following classic result.

\begin{lem}\label{trans5} If $A$ is a nonempty additive sub-semigroup of $\N$ then there exists $N$ such that $nk \in A$ for all $n \geq N$ where
$k$ is the greatest common divisor of $A$. \end{lem}

{\bfseries Proof:} $A - A$ is a non-trivial additive subgroup of $\Z$ and so equals $k \Z$ where $k$ is the smallest positive element of $A - A$.
Dividing through by $k$ we may assume that that greatest common divisor is $1$. So there exists $m \in \N$ such that $m, m+1 \in A$.
If $n \geq m^2 $ then with $0 \leq r < m$ and $q \geq m - 1$, $n = q m + r = (q - r)m + r(m + 1) \in A$.

$\Box$ \vspace{.5cm}

By assumption, we can choose $d, \ep, x_0$ and $y_0$ so that infinitely often $i \not\in N(x_0,y_0)$. Since
$N(x_0,x_0) + N(x_0,y_0) \subset N(x_0,y_0)$ it cannot happen that eventually $i \in N(x_0,x_0)$. That is, $k(x_0) > 1$.
Observe that $k(x)$ divides every element of $  N(x,y) + N(y,x) \subset N(x,x)$ and every element of $N(x,y) + N(y,y) + N(y,x) \subset N(x,x)$.
Consequently, $k(x)$ divides every element of $N(y,y)$ and so $k(x)|k(y)$. Interchanging $x$ and $y$ we see that there is an integer $k > 1$
such that $k(x) = k$ for all $x \in X$. It then follows that all of the elements of $N(x,y)$ are congruent mod $k$ with congruence class
inverse to to congruence class of the elements of $N(y,x)$. If $(x,y) \in f^p$ then $p \in N(x,y)$ and so the elements of $N(x,y)$ are
congruent to $p$ mod $k$. Fix a base point $x_0 \in X$. Map $X$ to $\Z_k$ by letting $h(x)$ be the mod $k$ congruence class of
the elements of $N(x_0,x)$. Observe that if $(x,y) \in f$ then $h(y) = h(x) + 1 = s(h(x))$. Since $f$ is surjective, $h$ maps $X$ onto
$\Z_k$ and maps $f$ onto $s$.

For uniform continuity, we prove that $h$ is constant on $V^d_{\ep}(x)$ for all $x$. Let $y \in X$ with $d(x,y) = \ep_1 < \ep$ and let
$\ep_2 = \ep - \ep_1$.
Since $f$ is $\U$ chain transitive, there exists $[a,b] \in f^{\times n}$ with $a_1 = x_0$ and $x_0 x$ chain-bound with respect to $d$ less than
$\ep_2$. Hence, $n \in N(x_0,x)$. Furthermore, the $x_0 y$ chain-bound with respect to $d$ is less than $\ep$.
Hence, $n \in N(x_0,y)$. Thus, $h(x) = h(y)$ is the congruence class of $n$ mod $k$.

(b) If $h : X \to \Z_k$ is a uniformly continuous surjection mapping $f^{-1}$ to $s$ then it maps $f$ to $s^{-1}$. The bijection $inv : t \mapsto -t$
maps $s^{-1}$ to $s$ and so $inv \circ h : X \to Z_k$ is a uniformly continuous surjection mapping $f$ to $s$. It follows from (a) that if
$f^{-1}$ is not $\U$ chain mixing then $f$ is not $\U$ chain mixing.

(c) We saw in the proof of (a) that if $f$ is not $\U$ chain mixing then, by (iv), there exists a positive integer such that
$f^k$ is not chain transitive.
Now assume that $f$ is a uniformly continuous map which is $\U$ chain mixing and that $k$ is a positive integer.

\begin{lem}\label{trans6} If $f$ is a uniformly continuous map, then for every $d \in \Gamma(\U), \ep > 0$ and positive integer $k$, there exists
$\bar d \in \Gamma(\U), \d > 0$ such that $ (V^{\bar d}_{\d} \circ f)^k \subset V^d_{\ep} \circ f^k$. \end{lem}

{\bfseries Proof:} By induction on $k$. For $k = 1$ let $d_1 = d$ and $\d = \ep$.

Assume $d_1 \in \Gamma(\U), \d_1 > 0$ such that $ (V^{d_1}_{\d_1} \circ f)^n \subset V^d_{\ep/2} \circ f^n$.  By uniform continuity of $f^n$
there exists $d_2 \in \Gamma(\U), \d_2 > 0$ such that $f^n \circ V^{d_2}_{\d_2} \subset V^d_{\ep/2} \circ f^n$. If $\bar d = d_1 + d_2 $ and
$\d = \min(\d_1,\d_2)$, then
$$ (V^{\bar d}_{\d} \circ f)^{n+1} \subset V^d_{\ep/2} \circ f^n \circ V^{d_2}_{\d_2} \circ f \subset V^d_{\ep} \circ f^{n+1}.$$

$\Box$ \vspace{.5cm}

Given $d \in \Gamma(\U), \ep > 0$ and a positive integer $k$ choose $\bar d$ and $\d$ as in Lemma \ref{trans6}. For $x,y \in X$
there exists $N$ so that $y \in (V^{\bar d}_{\d} \circ f)^n(x)$ for all $n \geq N$. Since $nk \geq N$,
$$ y \in (V^{\bar d}_{\d} \circ f)^{nk}(x) \subset (V^d_{\ep} \circ f^k)^n(x). $$
Thus, $f^k$ is $\U$ chain mixing.

$\Box$ \vspace{.5cm}

Assume that $T$ is a set of positive integers directed by divisibility\index{directed by divisibility},
i.e.\ if $k_1, k_2 \in T$ then there exists $k_3 \in T$ with
$k_1|k_3$ and $k_2|k_3$.  If $k_1|k_2$ we let $\pi : \Z_{k_2} \to \Z_{k_1}$ be the cyclic group surjection induced by the inclusion
$k_2 \Z \subset k_1 \Z$. For the directed set $T$ we let $\Z_T = \{ t \in \Pi_{k \in T} \Z_k :  k_1|k_2 \Rightarrow \pi(t_{k_2}) = t_{k_1} \}$.
If $T$  is finite then $\Z_T$  is isomorphic to $\Z_{k}$ where $k$ is the maximum element of $T$. If $T$  is infinite, then $\Z_T$ is a compact
monothetic group, i.e.\ if $1 \in \Z_T$ the unit element which projects to $1 \in \Z_k$ for all $k \in T$, then the cyclic group
generated by $1$ is dense in $\Z_T$. We let $s_T$  be the translation by $1$ in $\Z_T$ which projects to $s_k$ on $\Z_k$ for all $k \in T$.
When $T$ is infinite, the dynamical system
consisting of the homeomorphism $s_T$ on the compact space $\Z_T$ is called the \emph{odometer} \index{odometer} associated
with $T$.

 \begin{theo}\label{trans7} Assume that $f$ is a $\U$ chain transitive relation on a uniform space $(X,\U)$. Let $T$ be the set of
 positive integers $k$ such that there is a $\U$ uniformly continuous map $h_k : X \to \Z_k$ which maps $f$ to $s_k$.

 \begin{enumerate}
 \item[(a)] The set $T$  is directed by divisibility.

 \item[(b)] If $T$  is infinite, then there exists a uniformly continuous map $h : X \to \Z_T$ with a dense image which maps $f$ to $s_T$.

 \item[(c)] If $T$  is finite with maximum element $k$ and the uniformly continuous $h_k : X \to \Z_k$ maps $f$ to $s_k$ then for each
 $i \in \Z_k$, $X_i = (h_k)^{-1}(i)$ is an $f^k$ invariant subset. If, in addition, $f$ is a $\U$ uniformly continuous map then the
 restriction $f^k|X_i$ is $\U$ chain mixing for each $i \in \Z_k$.
 \end{enumerate}
 \end{theo}

 {\bfseries Proof:} Fix a base point $e \in X$. If $h_k(e) = p$ then by replacing $h_k$ by the composition $(s_k)^{-p} \circ h_k$ we can
 assume that $h_k(e) = 0$. For each $k \in T$ we will assume that $h_k(e) = 0$. Let $E_k = (h_k \times h_k)^{-1}(1_{\Z_k})$. Since
 $h_k$ is $\U$ uniformly continuous, $E_k \in \U$ and it is a clopen equivalence relation on $X$. If $(x,y) \in f$ and $(y,y_1) \in E_k$
 then $h(y_1) = h(y) = s_k(h(x)) = h(x) + 1$. Thus, $h_k$ maps $E_k \circ f$ to $s_k$ and, since $h_k(e) = 0$, we see that
 \begin{equation}\label{transeq2}
 x \in (E_k \circ f)^n(e) \quad \Longrightarrow \quad h(x) = n \in \Z_k.
 \end{equation}
 Because $f$ is assumed to be $\U$ chain transitive, every $x \in X$ lies in  $(E_k \circ f)^n(e)$ for some $n \in \Z$.

 (a) For $k_1, k_2 \in T$ let $E = E_{k_1} \cap E_{k_2}$, a clopen equivalence relation in $\U$. From (\ref{transeq2}) it follows
 that $(h_{k_1},h_{k_2}) : X \to \Z_{k_1} \times \Z_{k_2}$ maps $E \circ f$ to the restriction of $s_{k_1} \times s_{k_2}$ on the cyclic
 subgroup generated by $(1,1)$, which has order the least common multiple $k$ of $k_1$ and $k_2$. This restriction can be identified with
 $s_k$ on $\Z_k$.  Thus, $k \in T$.

 (b) If $k_1|k_2$ in $T$ and $E = E_{k_1} \cap E_{k_2}$ then for $x \in (E \circ f)^n(e)$, $h_{k_1}(x) = n \in \Z_{k_1}$ and $h_{k_2}(x) = n \in \Z_{k_2}$.
 Hence, with $\pi : \Z_{k_2} \to \Z_{k_1}$ the projection we see that $\pi(h_{k_2}(x)) = h_{k_1}(x)$.  It follows that
 the product $h_T = \Pi_{k \in T} h_k$ maps $X$ to $\Z_T$ taking $f$ to $s_T$. Since each fact is uniformly continuous, the map $h_T$
 is uniformly continuous. Since each $h_k$ is surjective, it follows that the image is dense in $\Z_T$.

 Notice that from (\ref{transeq2}) it follows that the $h_k$'s and $h_T$ are uniquely determined by the condition that $e$ is mapped to $0$.

 (c) Let $k \in T$. If $(x,y) \in f^k$ then $h_k(x) = h_k(y)$.  Since $f^k$ is a surjective relation,
 it follows that  $X_i$ is $f^k$ invariant for each $i \in \Z_k$.

 Now assume that $f$ is a uniformly continuous map and that some $f^k|X_i$ is not $\U$ chain mixing. By changing the choice of base point
 and translating, we may assume that $i = 0$. We will show that $k$ is not the maximum element of $T$.

 Since $f^k|X_0$ is not $\U$ chain mixing, there is an integer $p > 1$ and a uniformly continuous map $g_p : X_0 \to \Z_p$ taking  $f^k|X_0$
 to $s_p$. Label the congruence classes of $\Z_k$ by $i = 0,\dots k-1$, of $\Z_p$ by $j = 0,\dots p-1$ and of $\Z_{kp}$
 by $kj + i$. Observe that if $x \in X_i$ then $f^{k-i}(x) \in X_0$. Define the map $H : X \to \Z_{kp}$ by
 \begin{equation}\label{transeq3}
 H(x) \ = \ k g_p(f^{k-i}(x)) + i \quad \text{if} \ \ x \in X_i,
 \end{equation}
 We see that if $i < k-1$ then $f(x) \in X_{i+1}$ and so $f^{k - (i+1)}(f(x)) = f^{k - i}(x)$. If $i = k - 1$ then $f(x) \in X_0$
 and so $H(f(x)) = g_p(f^k(f(x))) = g_p(f(x)) + 1$ provided $g_p(f(x)) < p - 1$ and $ = 0$ if $g_p(f(x)) = p - 1$. Hence,

 \begin{equation}\label{transeq4}
 H(f(x)) \ = \ \begin{cases} k g_p(f^{k-i}(x)) + i + 1 \quad \text{if} \ i < k -1, \\
 k (g_p(f(x)) + 1) + 0 \quad \ \text{if} \ i = k -1, g_p(f(x)) < p - 1, \\
 \qquad \quad 0 \qquad \qquad \qquad \ \text{if} \ i = k -1, g_p(f(x)) = p - 1. \end{cases}
 \end{equation}
 It is clear that $H$ is $\U$ uniformly continuous since $h_k, g_p$ and $f$ are. From (\ref{transeq4}) we see that $H$ maps $f$ to $s_{pk}$.
 Hence, $pk \in T$ and so $k$ is not the maximum element.

$\Box$ \vspace{.5cm}

 {\bfseries Remark:} Without compactness of $X$ the map $h_T$ in (b) need not be surjective.
 For example, let $X $ be the dense cyclic subgroup generated by $1_T$ in $\Z_T$, or, more generally, any
 proper, $s_T$ invariant subset of an odometer $\Z_T$ which includes $0$. With the  uniformity induced from $\Z_T$ the homeomorphism $s_T$
 is a uniform isomorphism of $X$. Choose $e = 0$. Since every orbit of $s_T$ is dense, $s_T$ is $\U$ chain transitive
 on $X$. For every  $k \in T$, the projection map $\Z_T \to \Z_k$ maps $s_T$ to $s_k$ and is surjective on $X$. But
 $h_T : X \to \Z_T$ is just the inclusion.
 \vspace{.5cm}

 \begin{cor}\label{trans8} Let $f$ be a surjective relation on a connected uniform space $(X,\U)$.

 The following conditions are equivalent.
 \begin{itemize}
 \item[(i)] The relation $f$ is $\U$ chain mixing.

 \item[(ii)] The relation $f$ is $\U$ chain transitive.

 \item[(iii)] The relation $f$ is $\U$ chain recurrent, i.e.\ $\CC_{\U} f$ is an equivalence relation.

 \item[(iv)] The relation $\CC_{\U} f$ is reflexive, i.e.\ $1_X \subset \CC_{\U} f$.
 \end{itemize}
 \end{cor}

 {\bfseries Proof:} It is obvious that (i) $\Rightarrow$ (ii) $\Rightarrow$ (iii) $\Rightarrow$ (iv).

 Since a connected space does not admit a continuous surjection onto a nontrivial finite set, (ii) $\Rightarrow$ (i) by Theorem \ref{trans4} (a).

As in the proof of Corollary  \ref{corquasi}
 $X/\CC_{\U} 1_X$ is totally disconnected, but as the continuous image of the connected space $X$
 it is connected and so the quotient is a singleton. Hence, $\CC_{\U} 1_X = X \times X$, i.e.\ the identity map is $\U$ chain transitive.
 So if $1_X \subset \CC_{\U} f$ then $X \times X = \CC_{\U} 1_X \subset \CC_{\U} \CC_{\U} f = \CC_{\U} f$. Thus, (iv) $\Rightarrow$ (ii).

$\Box$

\vspace{1cm}

\section{ The Ma\~{n}\'{e} Set in the Compact, Metrizable Case}
\vspace{.5cm}

Throughout this section $X$ is a compact metrizable space.  A compact space is metrizable iff it is Hausdorff and second countable. In that case,
every continuous metric $d$ on $X$ is an element of $\Gamma(\U)$ where $\U$ is the unique uniformity which consists of all neighborhoods of the
diagonal. In particular, $\U = \U(d)$ for each such metric. Thus, for a compact metrizable space with unique uniformity $\U$,
$\Gamma_m(\U) = \Gamma_m(X)$.

If $E$ is a closed equivalence relation on $X$ then the quotient $X/E$ is a compact metrizable space by Proposition \ref{eqre14prop}. If the quotient is
totally disconnected then it is strongly zero-dimensional by Proposition \ref{propstronglysigma}.

We let $\CC f $ denote
$\CC_{\U} f$ where $\U$ is the unique uniformity. By Theorem \ref{theo11u}
  $\CC f = \CC_d f$ for every
$d \in \Gamma_m(X)$ and
 $\G f = \bigcap_{d \in \Gamma_m(X)} \ \A_d f$.  On the other hand, the union
is not obviously closed or transitive. We prove that it is both using an idea from \cite{W16}.

For $V$ a neighborhood of the diagonal $1_X \subset X \times X  $ and a pair
$x,y \in X$, $[a,b] \in f^{\times n}$
defines an \emph{$xy, V$ chain } if $(x,a_1), (b_n,y)$ and $(b_i,a_{i+1})$ are in $V$ for  $i = 1, \dots, n-1$.
We will call $n$ the \emph{length} of the chain.\index{chain!length}

\begin{df} Let $\W f$ denote the set of pairs $(x,y) \in  X \times X$ such that
for every neighborhood $W$ of $1_X$ there exists a closed, symmetric neighborhood $V$ of $1_X$ and $n \in \N$
such that there is an $xy, V$ chain of length $n$ and $V^{3^{n}} \subset W$. \end{df}
\vspace{.5cm}\index{w@$\W f$}

\begin{theo}\label{Wise1} For a relation $f$ on a compact, metrizable space $X$, the
relation $\W f$ is a closed, transitive relation and
$\W f  = \bigcup_{d \in \Gamma_m(X)}  \ \A_d f $. \end{theo}

{\bfseries Proof:}  We will prove that
$\bigcup_d \A_d f  \subset \W f $ and $ \ol{\W f \ \cup \ (\W f)^2} \subset \bigcup_d \A_d f$.

Let $(x,y) \in \A_d f$ for some metric $d$ on $X$ and let $W$ be a neighborhood of the diagonal.
Choose $\ep > 0$ so that $\bar V^d_{3 \ep} \subset W$. Since $(x,y) \in \A_d f$ there exists $[a,b] \in f^{\times n}$ for
some $n \in \N$ so that with respect to $d$ the $xy$ chain-length of $[a,b]$ is less than $ \ep$. Write $b_0 = x$ and
$a_{n+1} = y$. Define
$ \ep_i = d(b_i,a_{i+1}) $ for $i = 0, \dots , n$. Thus, $\sum_i \ep_i < \ep$.
Define
$$ V \quad = \quad \bar V^d_{\ep /3^{n+1}} \ \cup \ ( \bigcup_{i= 0}^n \
\bar V^d_{\ep_i}(b_i) \times \bar V^d_{\ep_i}(b_i) ).$$
Clearly, $[a,b]$ defines an $xy, V$.

We show that if $ (w,z) \in V^{3^{n+1}}$ then $(w,z) \in \bar V^d_{3 \ep}$.  There exists a sequence
$w = u_0, u_1, \dots , u_N = z$ with $(u_i, u_{i+1} ) \in V$ for $i = 0, \dots , N-1$ and with $N \leq 3^{n+1}$.
Choose the sequence so that $N$ is minimal. If $u_j, u_{j + k} \in \bar V^d_{\ep_i}(b_i) $ with $k > 0$ then
$k = 1$ for otherwise we could eliminate the terms $u_{j+1}, \dots, u_{j + k -1}$ and obtain a sequence with $N$ smaller.
Thus, for each $i$ there is at most one $j$ such that $u_j, u_{j+1} \in \bar V^d_{\ep_i}(b_i) $. For the remaining
$j$'s, $(u_j,u_{j+1}) \in  \bar V^d_{\ep /3^{n+1}}$. It follows that
$$ \Sigma_{i=0}^{N-1} d(u_j,u_{j+1}) \quad \leq \quad N \cdot
(\ep /3^{n+1}) \ + \ 2 \Sigma_i \ep_i \quad \leq \quad 3 \ep.$$
By the triangle inequality $d(w,z) = d(u_0,u_N) \leq 3 \ep$.

It follows that $\bigcup_d \ \A_d f \subset \W f$.

Now assume that  $(x,y) \in \ol{\W f}$. We will use the Metrization Lemma for uniform spaces, \cite{K} Lemma 6.12,
to construct a metric $d$ such that
$(x,y) \in \A_d f$.  We will then indicate how to adjust the proof to obtain the required metric when
$(x,y) \in \ol{(\W f)^2}$.

Fix some metric $d_0 $ on $X$ which is bounded by $1$.

Let $U_0 = A_0 = X \times X = V^{d_0}_{1}$ and $M_0 = 0$.  Assume that, inductively, the closed
symmetric neighborhood of the diagonal (= csn)\index{csn} $U_{M_k} = A_k \subset V^{d_0}_{2^{-k}}$ has been
constructed. There exists $(x_k,y_k) \in \W f$ such that
$(x,x_k), (y_k,y) \in A_k$. Hence, there exists $n_k \in \N$ and a csn $B_{k}$ such
that there is a $x_k y_k, B_k$ chain length $n_k $ and $(B_k)^{3^{n_k}} \subset A_k \cap V^{d_0}_{2^{-k-1}}$.

We now interpolate powers of $B_k$ between $B_k = A_{k+1}$ and $A_k$.

For $i = 1, \dots , n_k + 1$ let $U_{M_k + i} = (B_k)^{3^{n_k + 1 - i}}$. Let
$M_{k+1} =  M_k + n_k + 1$  and $A_{k+1} = U_{M_{k+1}}$.

Thus, $\{ U_j \}$ is a sequence of csn's with $(U_{j+1}^3 \subset U_j$  and $B_j = A_{j+1} = U_{M_{j+1}}$ for $j \geq 0$.

From the Metrization Lemma we obtain a metric $d$ such that
$ U_j \subset V^d_{2^{-j}} \subset U_{j - 1} $ for all $j \in \N$.

It follows that with respect to $d$ the $x_k y_k$ length of the
$B_k$ chain is bounded by
$$(n_k + 1) 2^{-M_{k + 1}}
\ = \ (n_k + 1) 2^{-(M_k + n_k + 1)} \ \leq \ 2^{-M_k}$$
and
since $(x,x_k), (y_k,y) \in A_k$ we have that the $xy$ length is bounded by
$3 \cdot 2^{-M_k}$.  Since $M_k \geq k$ it follows
that $(x,y) \in \A_d f$.

If $(x,y) \in \ol{(\W f)^2}$ then there exist $(x_k,z_k), (z_k,y_k) \in \W f$ with $(x,x_k),$ $ (y_k,y) \in A_k$. We begin
with a $n_k \in \N$ and a csn $B_{k}$ such that there is a $x_k z_k, B_k$ chain
of size $n_k $ and $(B_k)^{3^{n_k}} \subset A_k \cap V^{d_0}_{2^{-k-1}}$.
Then choose an $m_k \in \N$ and a csn $C_{k}$ such
that there is a $z_k y_k, C_k$ chain
of size $m_k $ and $(C_k)^{3^{m_k}} \subset B_k$.

This time for $i = 1, \dots , n_k + 1$ let $U_{M_k + i} = (B_k)^{3^{n_k + 1 - i}}$ and
$j = 1, \dots , m_k + 1$ let $U_{M_k + n_k + 1 + j} = (C_k)^{3^{m_k + 1 - j}}$. Let $M_{k+1} =  M_k + n_k + m_k + 2$
and $A_{k+1} = U_{M_{k+1}} = C_k$. Estimate as before to get that the $xy$ length of the $B_k$ chain followed by the
$C_k$ chain (with $z_k$ omitted between them) is at most $ 4 \cdot 2^{-M_k}$. Again $(x,y) \in \A_d f$.

$\Box$ \vspace{.5cm}

Following Fathi and Pageault \cite {FP}, we call $|\W f|$
 the \emph{Ma\~{n}\'{e} set}.\index{Ma\~{n}\'{e} set}

For every  $d \in \Gamma_m(X)$ on the compact metrizable space $X$ we have
\begin{equation}\label{mane01}
\G f \ \subset \ \A_d f \ \subset \ \W f \ \subset \ \CC f.
\end{equation}

Using Theorem \ref{Wise1} we follow \cite{W17} to prove the following extension of a theorem of Fathi and Pageault,
see \cite{FP}.

\begin{theo}\label{manetheo} Let $f$ be a continuous map on a compact, metrizable space $X$ such that $f^{-1}(|f|) = |f|$ and
let $K = X \setminus |f|^{\circ}$. $\W f = 1_{|f|} \cup \CC (f|K)$. Hence,
$|\W f| = |f| \cup |\CC (f|K)|$. \end{theo}

{\bfseries Proof:} From $f^{-1}(|f|) = |f|$ it follows that $K$ is $f$ $^+$invariant.
Because $f$ is a map, $f = 1_{|f|} \cup (f|K)$.

For any metric $d$, equation (\ref{eq38}) implies that
$$\A_d f = \A_d (1_{|f|} \cup (f|K)) = 1_{|f|} \cup \A_d (f|K)$$
and so
$$\W f = \W (1_{|f|} \cup (f|K)) = 1_{|f|} \cup \W (f|K) \subset 1_{|f|} \cup \CC (f|K).$$

To complete the proof we assume that $(x,y) \in \CC (f|K)$ and show that $(x,y) \in \W f$. Fix a metric $d$ on $X$ and
let $W$ be an arbitrary neighborhood of $1_X$. Choose $\ep > 0$ so that $V^d_{4 \ep} \subset W$. Let $\d > 0$ be
such that $\d < \ep/2$ and $d(x,y) < \d $ implies $d(f(x),f(y)) < \ep/2$. Choose $[a,b] \in (f|K)^{\times n}$ of
minimum size $n$ such that the $xy$ chain-bound is less than $\d$. We may perturb so that $a_i \not\in |f|$ for
$i = 1, \dots, n$ and so $b_i = f(a_i) \not\in |f|$ since $|f| = f^{-1}(|f|)$ by assumption.
Let $b_0 = x$ and $a_{n+1} = y$.
If $1 \leq i < j \leq n+1$ then
$a_i \not= a_j$, $b_{i-1} \not= a_j$ and $b_{i-1} \not= b_{j-1}$ for if not we could shorten the chain by removing
the pairs $(a_k,b_k)$ for $k = i, \dots, j-1$ contradicting the minimality of $n$. Now if $a_i = b_{j-1}$ then
$j > i + 1$ since $a_j \not\in |f|$. Let $i'$ be the smallest index such that $a_{i'} = b_{j-1}$ for some $j > i' + 1$ and
let $j'$ be the largest such $j$ for $i'$. Eliminate the pairs $(a_k,b_k)$ for $k = i'+1, \dots, j'$. Observe that
\begin{equation}
\begin{split}
d(b_{i'},a_{j'+1}) \leq  d(b_{i'},b_{j'}) + d(b_{j'},a_{j'+1}) = \hspace{2cm}\\
d(f(b_{j'-1}),f(a_{j'})) + d(b_{j'},a_{j'+1}) \leq \ep. \hspace{1cm}
\end{split}
\end{equation}
Moving right we may have to do several of these truncations, which do not overlap, and so eventually,
we obtain $[a,b] \in f^{\times n'}$ with $\{ a_i,b_{i-1} \} \cap \{ a_j,b_{j-1} \} = \emptyset$ if $i \not= j$.

Choose $0 < \d_C < \ep$ small enough that the sets $C_i =  \bar V^d_{\d_C}(\{ a_i,b_{i-1} \})$ are pairwise disjoint
for $i = 1, \dots n'+1$. Let $\ep > \ep_0 > 0$ be smaller than the distance between $C_i$ and $C_j$ if $i \not= j$.
Let $V = \bar V^d_{\ep_0/3^{n'}} \cup \bigcup_{i = 1}^{n+1} C_i \times C_i$.
Clearly, $[a,b]$ defines an $xy, V$ chain.

If $z_1,\dots , z_M$ satisfies $(z_i,z_{i+1}) \in V$ and $M \leq 3^{n'}$
then since $\ep_0$ is smaller than the distance between
the $C_i$'s at most one pair $\{ z_k,z_{k+1} \}$ lies in some $C_i$. Hence,
\begin{equation}
\begin{split}
d(z_1,z_M) \leq  \Sigma_{k = 1}^{M-1} d(z_k,z_{k+1}) \leq \hspace{3cm} \\
3^{n'} \cdot (\ep_0/3^{n'}) + 2 \max \ diam C_i \leq \ep_0 + \ep + 2 \d_C \leq 4 \ep.\hspace{1cm}
\end{split}
\end{equation}
Hence, $(z_1,z_M) \in W$.

$\Box$ \vspace{.5cm}

The following extension of Corollary \ref{cor6ubb} is easy to check.

\begin{prop} \label{prop19}  If $f$ is a Lipschitz map on $(X,d)$, then
 \begin{equation}\label{eq63}
 \A_d (f^n)  \quad \subset \quad \A_d f \quad = \quad f^{[1,n]} \cup [(\A_d (f^n)) \circ f^{[0,n]}],
 \end{equation}
 and so $|\A_d(f^n)| = |\A_d f|$.
  \end{prop}

$\Box$ \vspace{.5cm}

 For a continuous map $f$ on $(X,d)$ let $Per(f)$ denote the set of periodic points,
 so that $Per(f) = \bigcup_{n = 1}^{\infty} \ |f^n|$.
Let $Per(f)^{\circ \circ} = \bigcup_{n = 1}^{\infty} \ |f^n|^{\circ}$.

\begin{lem}\label{lem20a} The open set $Per(f)^{\circ \circ}$ is dense in $Per(f)^{\circ}$, the interior of the
set of periodic points. \end{lem}

{\bfseries Proof:} Each $|f^n|$ is closed in $X$. Let $U$ be a nonempty open subset of $Per(f)$. It is the
countable union of the relatively closed sets $|f^n|\cap U$ and so by the Baire Category Theorem at least one of these
has a nonempty interior.

 $\Box$ \vspace{.5cm}

While, $Per(f)^{\circ \circ}$ is contained in the interior of $Per(f)$,
but might be a proper subset of it. By periodicity each $|f^n|$ is $f$ invariant, i.e.\ $f(|f^n|) = |f^n|$. So if $f$ is a
homeomorphism  each $|f^n|^{\circ}$ is invariant as well. Thus, if $f$ is a homeomorphism, $Per(f)^{\circ \circ}$ is an open
invariant set and its complement in $X$ is a closed invariant set. Notice also that if $A$ is any closed subset of $X$
which is $f$ $^+$invariant then it is $f^n$ $^+$invariant and $(f|A)^n = (f^n)|A$.

In the Lipschitz case we can extend the above results.

\begin{cor}\label{cor21} Let $f$ be a homeomorphism on $(X,d)$.  If $f$ is a Lipschitz map then
\begin{equation}\label{eq66a}
|\A_d f| \quad \subset \quad Per(f) \cup |\CC (f|(X \setminus Per(f)^{\circ \circ}))|.
\end{equation}
\end{cor}

{\bfseries Proof:} Let $X_n = X \setminus |f^{n!}|^{\circ}$.  By Proposition \ref{prop3} and Proposition \ref{prop19} we have
\begin{equation}\label{eq66b}
|\A_d f | \ = \ |\A_d (f^{n!})| \ \subset \ |f^{n!}| \cup |\CC ((f|X_n)^{n!})|
 \ \subset \ Per(f) \cup |\CC (f|X_{n})|.
\end{equation}
Now $\{ X_n \}$ is a decreasing sequence of closed invariant sets with intersection
$X_{\infty} = X \setminus Per(f)^{\circ \circ}$. Hence,  $\{ f|X_{n} = f \cap (X_n \times X_n) \}$
is a decreasing sequence of closed relations with intersection $f|X_{\infty}$. By \cite{A93} Theorem 7.23, the map
$R \mapsto |\CC R|$ is a monotone, usc function on closed sets and so
\begin{equation}\label{eq66c}
\bigcap_n \ |\CC (f|X_{n})| \quad = \quad | \CC (\bigcap_n \ (f|X_n))| \quad = \quad |\CC (f|X_{\infty})|.
\end{equation}
Together with (\ref{eq66b}) this implies (\ref{eq66a})

$\Box$ \vspace{.5cm}

\begin{ex}\label{counterex1} Without the Lipschitz assumption the result is not true. \end{ex}

{\bfseries Proof:}  On $I = [0,1]$ let $\mu$ be a full, nonatomic
probability measure concentrated on a dense countable union
of Cantor sets of Lebesgue measure zero. Let $\pi : I \to I$ be the distribution function so that $\pi(t) = \mu([0,t])$.
Then $\pi$ is a homeomorphism on $I$ fixing the end-points. Let $X_0 = I \times \{ -1, 0, 1 \} $ with the
metric $d_0((s,a),(t,b)) = |s - t| + |a - b|$. Let $\tilde{\pi} : X_0 \to X_0$ be the homeomorphism defined by
$\tilde{\pi}(t,-1) = (\pi(t),-1)$ and  $\tilde \pi(t,a) = (t,a)$ for
$a = 0,1$. Let $d$ be the metric $d_0$ pulled back by
$\tilde \pi$.  Thus,  if $s < t \in I$ then $d((s,a),(t,a)) = t - s$ if $a = 0,1$ and $\mu([s,t])$ if $a = -1$.
Let $L = \{ 0 \} \times \{ -1, 0, 1 \}, R = \{ 1 \} \times \{ -1, 0, 1 \}$. Let
$E = 1_{X_0} \cup (L \times L) \cup (R \times R)$ and
use the metric $d =\ell_d^E = s\ell_d^E$ on the quotient space $X = X_0/E$ with quotient map $q : X_0 \to X$. For $a = -1, 0, 1$ let $I_a$ denote
$q(I \times \{ a \})$. It is easy to check that  each restriction $q : I \times \{ a \} \to I_a$
is an isometry.

Now define the homeomorphism $f$ on $X$ by
\begin{equation}\label{eq67}
f(t,a) \quad  = \quad \begin{cases} \quad  (t,-a)\qquad \mbox{for} \ a = \pm 1 \\
 \quad (t^2,a) \qquad \mbox{for} \ \ a = 0. \end{cases}
\end{equation}

Let $Y$ be the subspace of $X$ which is the quotient of $I \times \{ 0,1 \}  \subset X_0$, i.e.\
$Y = I_1 \cup I_0 \subset X$, and define
$g$ on $Y$ and $h : X \to Y$ by
\begin{equation}\label{eq68}
\begin{split}
g(t,a) \quad  = \quad \begin{cases} \quad  (t,a)\qquad \mbox{for} \ a =  1 \\
 \quad (t^2,a) \qquad \mbox{for} \ \ a = 0 \end{cases} \\
h(t,a) \quad  = \quad \begin{cases} \quad  (t,1)\qquad \mbox{for} \ a = \pm 1 \\
 \quad (t,a) \qquad \mbox{for} \ \ a = 0 \end{cases}
\end{split}
\end{equation}

Neither $f$ nor $h$ is Lipschitz.  Because $g$ is the identity on $I_1$ and $f^2$ is the
identity on $I_1 \cup I_{-1}$ it follows from Proposition \ref{prop3} that
\begin{equation}\label{eq69}
\begin{split}
\A_d (f^2) \quad = \quad 1_{(I_1 \cup I_{-1})} \cup \A_d ((f|I_0)^2) \\
\mbox{and} \qquad \A_d g \quad = \quad 1_{I_1} \cup \A_d (f|I_0).
\end{split}
\end{equation}
For $f|I_0 = g|I_0$, $L(t,0) = 1 - t$ is a Lipschitz Lyapunov function
which is increasing on all orbits except the -fixed- endpoints.
It follows that $|\CC (f|I_0)| = |\CC ((f|I_0)^2)|  = \{ q(0,0), q(1,0) \}$.

We will show that
\begin{equation}\label{eq70}
\A_d f \quad = \quad X \times X
\end{equation}

Thus, for $0 < t < 1$ the point $(t,0) \in |\A_d f|$ but is not in $|\A_d f^2|$ and $(t,0) = h(t,0)$ is not
in $|\A_d g|$.

Let $s < t$ in $I$. Because $\mu$ and Lebesgue measure $\lambda$ are mutually singular
we can choose for any $\ep > 0$ an increasing sequence $s = u_1,....,u_{2n+1} = t$ so that
$\mu(\bigcup_{i=1}^n \ [u_{2i -1},u_{2i}]) <  \ep$ but $\lambda(\bigcup_{i=1}^n \ [u_{2i -1},u_{2i}]) \geq 1 - \ep$
and so $\lambda(\bigcup_{i=1}^n \ [u_{2i},u_{2i + 1}] < \ep$. On $I_1$ the length of an interval is its
Lebesgue measure while on $I_{-1}$ the length is its $\mu$ measure. Thus, if $x = (s,1)$ and $(y = (t,1)$ then
\begin{equation}\label{eq71}
(u_1,1),(u_2,-1),(u_3,1),(u_4,-1),....,(u_{2n},-1)
\end{equation}
each paired with its image under $f$, defines a sequence in $f^{\times 2n}$ whose $xy$ chain-length is less than
$2 \ep$. Since  $f$ is symmetric on $I_1 \cup I_{-1}$ we can reverse the sequence to get one whose $yx$ chain-length
is the same.  Thus, any two elements of $I_1 \cup I_{-1}$ are $\A_d f$ equivalent.  On the other hand, it is easy
to check that for any $t \in (0,1)$, $q(t,0) \in \G (f|I_0)(q(1,0))$ and $q(0,0) \in \G (f|I_0)(t,0)$.  It follows that
any two elements of $X$ are $\A_d f$ equivalent.

On the invariant set $X_{\pm 1} \ =_{def} \ I_1 \cup I_{-1}$ the restriction of $f$ has order $2$ and so
$E = 1_{X_{\pm 1}} \cup f|X_{\pm 1}$ is a closed equivalence relation with each equivalence class having one or two points.
However, the pseudo-metric $s \ell^K_d =\ell^K_d$ is identically zero and so does not induce a metric on the space of
equivalence classes.

$\Box$ \vspace{.5cm}

\begin{ex}\label{counterex2} In (\ref{eq66a}) the inclusion may fail if $Per(f)^{\circ \circ}$ is replaced by $Per(f)^{\circ}$
and it may fail
if $|\CC (f|(X \setminus Per(f)^{\circ \circ}))|$ is replaced by $|\A_d (f|(X \setminus Per(f)^{\circ \circ}))|$.\end{ex}

{\bfseries Proof:} Let $S$ be the unit circle in the complex plane.  Let
\begin{equation}\label{eq72}
X = ([-1,1] \cup S) \times \{ 0 \} \ \cup \ (\bigcup_{n=1}^{\infty} S \times \{ 1/n \}).
\end{equation}
equipped the restriction of the Eucidean metric from $\R^3$.

On $X$ define the Lipschitz homeomorphism $f$ by
\begin{equation}\label{eq73}
f(x,t) \quad = \quad \begin{cases} \ \qquad (x \cdot e^{2 \pi i t},t) \quad \mbox{for} \ x \in S, \hspace{.5cm} \\
(\frac{1}{2}(x^2 + 2x - 1),0) \quad \mbox{for} \ x \in [-1,1], t = 0.\end{cases}
\end{equation}
That is, on $[-1,1] \times \{ 0 \}$ the map is conjugate to $x \mapsto x^2$ via the homeomorphism
$(x,0) \mapsto (x + 1)/2$ from $[-1,1] \times \{ 0 \}$ to $[0,1]$.

$Per(f) = S \times (\{ 0, 1, 1/2,... \})$ and so $Y =_{def}  X \setminus Per(f)^{\circ}$ is
$[-1,1] \times \{ 0 \}$.  For the restriction
of $f$ to this set, the only chain recurrent points are the endpoints, i.e.\ $|\CC_d (f|Y)| = \{ (-1,0), (1,0) \}$.

$X_{\infty} = X \setminus Per(f)^{\circ \circ} = (S \cup [-1,1]) \times \{ 0 \} $. For the restriction of $f$ to
$X_{\infty}$ every point is chain recurrent, i.e.\ $|\CC_d (f|X_{\infty})| = X_{\infty}$, but from Proposition \ref{prop3},
\begin{equation}\label{eq74}
S \times \{ 0 \} \quad = \quad |(f|X_{\infty})| \quad = \quad | \G (f|X_{\infty})| \quad = \quad |\A_d (f|X_{\infty})|
\end{equation}

Finally, it is easy to check that for $f$ itself
\begin{equation}
X \quad = \quad |\G f| \quad = \quad |\A_d f|.
\end{equation}

Thus, in  (\ref{eq66a}) the equation fails if $Per(f)^{\circ \circ}$ is replaced by $Per(f)^{\circ}$ and it fails
if $|\CC (f|(X \setminus Per(f)^{\circ \circ}))|$ is replaced by $|\A_d (f|(X \setminus Per(f)^{\circ \circ}))|$.

$\Box$

\newpage

\section{ Appendix A: Directed Sets and Nets}\label{appendix-nets}
\vspace{.5cm}

We review the theory of nets, following \cite[Chapter 2]{K}.

A set $I$ is directed by a reflexive, transitive relation $\prec$
 if for every $i_1, i_2 \in I$  there exists $j \in I$ such that $i_1, i_2 \prec j$. We call
 $I$ a  \emph{directed set}.\index{directed set} If $I_1, I_2$ are directed sets then the product $I_1 \times I_2$ is directed
 by the  product ordering $(i_1,j_1) \prec (i_2,j_2)$ when $i_1 \prec i_2$ and $j_1 \prec j_2$.

For $i \in I$ let $\prec_i = \{ j : i \prec j \}$. A set $F \subset I$ is  called \emph{terminal}\index{subset!terminal}\index{terminal subset}
if $F \supset \prec_i $ for some $i \in I$.   $F$ is called \emph{cofinal} if $F \cap \prec_i \not= \emptyset $ for all
$i \in I$.\index{subset!cofinal}\index{cofinal subset} In the family language of \cite{A97} these are dual families of subsets of $I$. Because the
set $I$ is directed by $\prec$ it follows that the family of terminal sets is a filter.  That is, a finite intersection of
terminal sets is terminal. The cofinal sets satisfy the dual, \emph{Ramsey Property}:\index{Ramsey Property} If a finite union of subsets
of $I$ is cofinal then at least one of them is cofinal.

For example, if $A \subset X$ then the set $\NN_A$ of neighborhoods of $A$ is
directed by $\supset$ and a subset of $\NN_A$ is cofinal iff it is a neighborhood base. If $A$ is the singleton $ \{ x \}$,
then we write $\NN_x$ for $\NN_A$.
 The sets $\Z_+$ and $\N$ are directed by $\leq$ and a subset is terminal iff it is cofinite. A subset is
 cofinal iff it is infinite.

A \emph{net}\index{net} in a set $Q$ is a function from a directed set $I$ to $Q$, denoted $\{ x_i : i \in I \}$. If $A \subset Q$ we
say that the net is \emph{eventually} (or \emph{frequently}) in $A$ if $\{ i : x_i \in A \}$ is terminal (resp.\  is cofinal).
\index{eventually in $A$}%
\index{frequently in $A$}%

A map $k : I' \to I$ between directed sets is a \emph{directed set morphism}\index{directed set morphism} if
$k^{-1}(F)$ is terminal in $I'$ whenever
$F$ is terminal in $I$. If $k$ is order-preserving, i.e.\ $i_1' \prec i_2' $
implies $k(i_1') \prec k(i_2') $, and, in addition,
 the image, $k(I'),$ is cofinal in
$I$ then $k$ is a morphism.

A map $k : I' \to I$ is a morphism iff whenever $F $ is cofinal in $I'$,
then $k(F)$ is cofinal in $I$. This follows because
$$k(F) \cap A \not= \emptyset \quad \Longleftrightarrow \quad F \cap k^{-1}(A) \not= \emptyset$$
and a set is  cofinal iff it meets every terminal set and vice-versa.

With this definition of morphism, the class of directed sets becomes a category.

If $i \mapsto x_i$ is a net, then the composite $i' \mapsto x_{k(i')}$ is the \emph{subnet}\index{subnet}
induced by the morphism $k$. We will usually suppress the mention of $k$ and just write $\{ x_{i'} : i' \in I' \}$ for the subnet.

If $x $ is a point of a topological space $X$ then a net in $X$ converges to $x$
(or has $x$ is a cluster point) if for every $U \in \NN_x$
the net is eventually in $U$ (resp.\   is frequently in $U$). Thus, if a net in $A$ has $x$ as a cluster point then
$x$ is in the closure of $A$. Conversely, if $x \in \ol{A}$ then we can use
$I = \NN_x$ and choose $x_U \in A \cap U$. We thus
obtain a net in $A$ converging to $x$. For a net $\{ x_i : i \in I \}$ in $X$ the set
of cluster points is $\bigcap_{i \in I} \ol{ \{ x_j : i \prec j \}}$. Equivalently, this is the set of limit points of
convergent subnets of $\{ x_i \}$.

\begin{lem}\label{applem02}  If $\{ x_i : i \in I \}$ is a net in $X$ and $A$
is a compact subset of $X$, then $A$ contains a cluster point of the
net iff $x_i \in U$ frequently for every open set containing $A$. \end{lem}

{\bfseries Proof:} Clearly, if $x \in A$ is a cluster point of the net, then it frequently enters every neighborhood of $x$ and
a fortiori it frequently enters every neighborhood of $A$. If for some $i$ the set
$K_i = \ol{ \{ x_j : i \prec j \}}$ is disjoint from $A$ then its complement
is a open set containing $A$ which the net does not enter frequently. So if the net frequently enters every neighborhood of $A$ then
$\{ K_i \cap A : i \in I \}$ is a collection of closed subsets of $A$ with the finite intersection property.
Hence, the intersection is nonempty by compactness.

$\Box$ \vspace{1cm}

\section{Appendix B: Uniform Spaces}\label{appendix-uniform}
\vspace{.5cm}

We review from \cite{K} Chapter 6 the facts we will need about uniform spaces.\index{space!uniform}\index{uniform space}

A uniformity $\U$ \index{uniformity} on a set $X$ is a filter of reflexive relations on $X$ which satisfies

\begin{itemize}
\item $U \in \U$ implies $U^{-1} \in \U$.

\item If $U \in \U$, then there exists $W \in \U$ such that $W \circ W \subset U$.
\end{itemize}

We say that a collection $\U_0$ of reflexive relations generates a uniformity when
$\U = \{ U : U \supset V$ for some $V \in \U_0\}$
is a uniformity.  This requires that if $V_1, V_2 \in \U_0 $,  there exists $V_3 \in \U_0$ so that
$V_3 \circ V_3 \subset V_1 \cap V_2^{-1}$.

For example, if $d$ is a pseudo-metric on $X$ then $V^d_{\ep} = \{ (x,y) : d(x,y) \leq \ep \}$ with $\ep > 0$
generates a uniformity   $\U(d)$ which we call
the uniformity associated with $d$.\index{u@$\U(d)$}

The \emph{gage} $\Gamma$\index{uniformity!gage}\index{gage} of a uniformity $\U$ (or $\Gamma(\U)$\index{g@$\Gamma(\U)$} when
we need to keep track of the uniformity)
is the set of all bounded pseudo-metrics $d$ such that $V^d_{\ep} \in \U$ for all
$\ep > 0$, or, equivalently, $\U(d) \subset \U$. From the
Metrization Lemma for uniformities, Lemma 6.12 of \cite{K}, it follows that if $U \in \U$ then
there exists $ d \in \Gamma$ such that $V^d_1 \subset U$.

A collection $\Gamma_0$ of pseudo-metrics generates a uniformity when \\ $\bigcup_{d \in \Gamma_0} \U(d)$ is a uniformity.
It suffices that if $d_1, d_2 \in \Gamma_0$, there exists $d_3 \in \Gamma_0$ such that $d_3 \leq K(d_1 + d_2)$ for some
positive $K$.

Since $\U$ is a filter, it is directed by $\supset$. If $d_1, d_2 \in \Gamma$ then $d_1 + d_2 \in \Gamma$, and so $\Gamma$ is
directed by $\leq$.

\begin{lem}\label{applem01a} Let $\{ d_1, d_2, \dots \}$ be a sequence in $\Gamma$ with  $d_i$ bounded by $K_i \geq 1$.
If $\{ a_1, a_2, \dots \}$ is a summable sequence of positive reals then $d = \Sigma_{k=1}^{\infty} \ (a_i/K_i) d_i $
is a pseudo-metric in
$\Gamma$.\end{lem}

{\bfseries Proof:} Dividing by $\Sigma_{k = 1}^{\infty} \ (a_i/K_i)$ we can assume the sum is $1$. Given $\ep > 0$ choose
$N$ so that $\Sigma_{k = N+1}^{\infty} \ (a_i/K_i) < \ep/2 $. Then $\bigcap_{k = 1}^N \ V^{d_k}_{\ep/2} \subset V^d_{\ep}$.

$\Box$ \vspace{.5cm}

Associated to a uniformity $\U$\index{uniformity!associated topology} is the $\U$ topology with $G$ open
iff $x \in G$ implies $U(x) \subset G$ for some $U \in \U$. The topology is
Hausdorff iff $1_X = \bigcap \ \{ U : U \in \U \}$, in which case we call $\U$ a Hausdorff uniformity.
If $X$ is a topological space then $\U$ is called \index{uniformity!compatible with the topology of $X$} compatible with the topology on $X$ if
$X$ has the $\U$ topology.

If $(X,\U), (Y,\V)$ are uniform spaces then the product uniformity\index{uniformity!product} $\U \times \V$ on $X \times Y$ is generated by the
product relations $U \times V$ for $U \in \U, V \in \V$. Given pseudo-metrics
in $\Gamma(\U)$ and $\Gamma(\V)$ the product pseudo-metrics on $X \times Y $
generate the gage $\Gamma(\U \times \V)$. The associated topology is the product of the $\U$ topology on $X$ with
the $\V$ topology on $Y$.

If $A \subset X$ then $\U|A$, the set of restrictions to $A$ of the relations $U \in \U$, is the
induced uniformity on $A$ with associated topology the subspace topology.  The restrictions to $A \times A$ of the
pseudo-metrics in $\Gamma(\U)$ generate the gage $\Gamma(\U|A)$.

Observe that if $E$ is an equivalence relation which contains the diagonal $1_X$ in its interior then every equivalence class is a neighborhood
of each of its points and so is open. It follows that $E = \bigcup_{x \in X} \{ E(x) \times E(x) \}$ is open in $X \times X$  and its complement
 $\bigcup_{(x,y) \not\in E} \{ E(x) \times E(y) \}$ is open as well. Thus, $E$ is a clopen equivalence relation. For a clopen
 equivalence relation $E$ on $X$, the characteristic function of $X \times X \setminus E$ is a continuous pseudo-ultrametric on $X$.

 We call a uniformity $\U$ \emph{zero-dimensional}\index{uniformity!zero-dimensional}
 when it is generated by equivalence relations. Equivalently, the gage is generated by
 pseudo-ultrametrics. In that case, the associated topology is zero-dimensional, i.e.\ the clopen subsets form a basis for the topology.
 Conversely, if $X$ is a zero-dimensional space then the set of all clopen equivalence relations on $X$ generates the maximum zero-dimensional
 uniformity compatible with the topology on $X$.  We denote it by $\U_{M_0}$. \index{u@$\U_{M_0}$} The gage $\Gamma(\U_{M_0})$ is generated by
 the pseudo-ultrametrics which are continuous on $X$. The class of zero-dimensional uniform spaces is closed under the
 operations of products and taking subspaces.

\begin{prop}\label{appprop02} Let $X$ be a  topological space. The following conditions are equivalent.

\begin{itemize}
\item[(a)] There exists a uniformity compatible with the topology on $X$.

\item[(b)] The topology on $X$ is completely regular.
That is, the continuous real-valued functions distinguish points and closed sets.
\end{itemize}

If $X$ is Hausdorff, then these are equivalent to
\begin{itemize}

\item[(c)] There exists a homeomorphism onto a subset of a compact Hausdorff space.
\end{itemize}
 \end{prop}

 {\bfseries Proof:} (a) $\Leftrightarrow$ (b) If $x$ is not in a closed set $A$ then
 there is a $d \in \Gamma$ such that
 $V^d_{\ep}(x) \cap A = \emptyset$ for some $\ep > 0$. The continuous function
 $y \mapsto min(d(x,y),1)$ is $0$ at $x$ and
 $1$ on $A$. If $X$ is completely regular then the uniformity generated by the pseudo-metrics
 $d_u(x,y) = |u(x) - u(y)|$, with $u$ varying over continuous
 real-valued functions, is compatible with the topology.

 (b)  $\Leftrightarrow$ (c)  Using bounded real-valued continuous functions we can
 embed a Hausdorff, completely regular space
 into a product of intervals.  On the other hand, by the Urysohn Lemma a compact
 Hausdorff space is completely regular and
 so any subspace is completely regular as well.

$\Box$ \vspace{.5cm}

A completely regular, Hausdorff space is called a \emph{Tychonoff} space. \index{space!Tychonoff}\index{Tychonoff space}
Clearly, a completely regular space $X$ is Tychonoff iff the points are closed,
i.e.\  iff $X$ is $T_1$.

If there is a metric in the gage then the $\U$ topology is Hausdorff, but the gage of a
Hausdorff uniformity  need not contain
a metric.

A map $h : X_1 \to X_2$ between uniform spaces is uniformly continuous \index{map!uniformly continuous}\index{uniformly continuous map}if $U \in \U_2$
implies $(h \times h)^{-1}(U) \in \U_1$, or, equivalently, if $h^*d \in \Gamma(\U_1)$ for
all $d \in \Gamma(\U_2)$ where $h^*d(x,y) = d(h(x),h(y))$.\index{h@$h^*d$}\index{d@$h^*d$}
A pseudo-metric $d$ on $X$ is in the gage of $\U$ iff $1_X : (X,\U) \to (X,\U(d))$ is uniformly continuous.
With the uniformity induced by the usual metric on $\R$, a pseudo-metric $d$ on $X$ is in the gage of $\U$ iff the map
$d :(X \times X, \U \times \U) \to \R$ is uniformly continuous.

In general, there may be many uniformities with the same associated topology. Given a completely regular space there is a maximum
uniformity $\U_M$ \index{u@$\U_{M}$}compatible with the topology. It is characterized by the condition
that any continuous map from $X$ to a uniform space is
uniformly continuous with respect to $\U_M$. If $X$ is paracompact then the
set of all neighborhoods of the diagonal is a uniformity
which is therefore $\U_M$.  If $X$ is compact, then this is the unique uniformity compatible with the topology on $X$.

A uniformity $\U$ on $X$ is \emph{totally bounded}\index{uniformity!totally bounded}\index{totally bounded uniformity}
 if for every $V \in \U$ the cover $\{ V(x) : x \in X \}$ has a finite
subcover, or, equivalently, if for every $d \in \Gamma(\U)$ the pseudo-metric space $(X,d)$ is totally bounded. Let
$\B(X,\U)$\index{b@$\B(X,\U)$} denote the Banach algebra of bounded, uniformly continuous, real-valued functions. If $u \in \B(X,\U)$
then the pseudo-metric $d_u$ defined by\index{d@$d_u$}
\begin{equation}\label{eqapp1}
d_u(x,y) \ = \  \ |u(x) - u(y)|. \hspace{2cm}
\end{equation}
is a totally bounded pseudo-metric in $\Gamma(\U)$. For $\B \subset \B(X,\U)$ a
closed subalgebra (assumed to contain the constant
functions) the pseudo-metrics $d_F = \Sigma_{u \in F} d_u$, with $F$ a finite subset of $ \B$, generate a totally bounded
uniformity $\T(\B) \subset \U$. \index{t@$\T(\B)$} If $\B$ is separable then the uniformity $\T(\B)$ is pseudo-metrizable. In fact,
if $\{ u_i \}$ is a dense sequence in the unit ball of $\B$ then
\begin{equation}\label{eqapp2}
d(x,y) \ = \ \Sigma_{i=1}^{\infty} \ 2^{-i} d_{u_i} \hspace{2cm}
\end{equation}
is a metric such that $\U(d) = \T(\B)$.

Recall that if $u \in \B$ then we can use the series expansion of the square root to show that
$|u| = \sqrt{u^2} \in \B$.  Hence, if $u_1, u_2 \in \B$ then $\max(u_1,u_2) = \frac{1}{2}(|u_1 - u_2| + u_1 + u_2)$ and
$\min(u_1,u_2) = - \max(-u_1,-u_2)$ are in $\B$.

The
subalgebra $\B$ \emph{distinguishes points and closed sets } when for every closed subset $A $ of $X$ and any
$x \in X \setminus A$ there exists $u \in \B$ such that $u(x) \not\in \ol{u(A)}$. Notice that if $|u(x) - t| < \ep$ implies
$t \not\in u(A)$, then $v(z) = \frac{1}{\ep} \min (|u(x) - u(z)|,\ep)$ is an element of $\B$ with
$v(x) = 0$ and $v = 1$ on $A$.  In that case, the topology associated
with $\T(\B)$ is that of $(X,\U)$, i.e.\ $\T(\B)$ is compatible with the topology of $X$. The uniformity
$\T(\B(X,\U))$ is the maximum totally bounded uniformity contained in $\U$ and we will denote it $\T(\U)$.
The gage of $\T(\U)$ consists of all the totally bounded pseudo-metrics in the gage of $\U$.

If $\{ x_i : i \in D \}$ and $\{ y_j : j \in I \}$ are nets in $X$, then they are \emph{$\U$-asymptotic}\index{a@$\U$-asymptotic nets}
for a uniformity
$\U$ on $X$ if the product net $\{ (x_i,y_j) : (i,j) \in D \times I \}$ is eventually in $U$ for all $U \in \U$.
The net $\{ x_i \}$ converges to $x \in X$ exactly when it is $\U$-asymptotic to a net constant at $x$. The
$\U$-asymptotic relation on nets on $X$ is symmetric and transitive, but not reflexive.
A net $\{ x_i \}$ is \emph{Cauchy}\index{net!Cauchy}\index{Cauchy net} when it
is $\U$-asymptotic to itself. The uniform space $(X,\U)$ is \emph{complete}\index{uniformity!complete}\index{complete uniformity}
 when every Cauchy net converges.
For a Hausdorff uniform space $(X,\U)$ there exists $j$ a uniform isomorphism from $(X,\U)$ onto
a dense subset of a complete, Hausdorff uniform space $(\bar X, \bar \U)$. Regarding $j$ as an inclusion, we call
$(\bar X, \bar \U)$ the \emph{completion}\index{uniformity!completion}\index{completion} of $(X,\U)$.
We can regard $\bar X$ as the space of the $\U$-asymptotic
equivalence classes of Cauchy nets in $X$. In general, if $(Y,\V)$ is a complete, Hausdorff uniform space
and $h : A \to Y$ is a uniformly continuous map on a subset $A$ of $X$
then $h$ extends uniquely to $\bar h : \ol{A} \to Y$
a uniformly continuous map on the closure. If $x \in \ol{A}$ and $\{ x_i \}$ is a net in $A$ converging to $x$
then $\{ h(x_i) \}$ is a Cauchy net in $Y$ and so converges to a unique point $\bar h(x)$.
It follows that the completion of a
Hausdorff uniform space is unique up to uniform isomorphism. For each
$d \in \Gamma(\U)$, the map $\bar d : \bar X \times \bar X \to
[0,M]$, where $M = \sup d$ is a pseudo-metric on $\bar X$ and these form the gage of $\bar \U$. If $d$ is a metric with
$\U = \U(d)$ then $ \bar d $ is a metric with $\bar \U = \U(\bar d)$.  That is, the completion
of a metric space is a metric
space.

A uniform space is compact iff it is totally bounded and complete. So the completion of a totally bounded, Hausdorff
uniform space is compact. In particular,  if $(X,\U)$ is Hausdorff and
$\B$ is a closed subalgebra of $\B(X,\U)$ which distinguishes points
and closed sets then the completion $(\bar X, \bar \T(\B))$ is a compact Hausdorff space. If $Y$ is a compact, Hausdorff
space (with its unique uniformity) and $h : (X,\U) \to Y$ is uniformly continuous then $h : (X,\T(\U)) \to Y$ is
uniformly continuous and so extends uniquely to $\bar h : (\bar X, \bar \T(\U)) \to Y$. If $\B$ is
closed subalgebra of $\B(X,\U)$ which distinguishes points and closed sets,
then with $\bar X$ the $\bar \T(\B)$ completion of $X$, the map $u \mapsto \bar u$ is a Banach algebra isomorphism
from $\B$ onto the Banach algebra of continuous, real-valued maps on $\bar X$. Thus, $\bar X$ is version the compactification
of $X$ obtained from $\B$ by the Gelfand space construction, see, e.g. \cite{A97} Chapter 5. In particular, if
$X$ is a Tychonoff space with $\U_M$ the maximum uniformity compatible with the topology then $(\bar X,\bar \T(\U_M))$ is
a version of the Stone-Cech compactification of $X$.

Finally, notice that  $(X,\U)$ has a second countable topology iff there exists a separable, closed subalgebra $\B$ of $\B(X,\U)$
which distinguishes points and closed sets. In that case, there is a metric $d$ such that $\T(\B) = \U(d)$ and
the associated compactification $\bar X$ is metrizable with metric $\bar d$.

\vspace{1cm}

\section{Appendix C: Proper Maps}\label{appendix-proper}
\vspace{.5cm}

A \emph{proper map} \index{map!proper}\index{proper map}$f : X \to Y$ is a
continuous map such that $f \times 1_Z : X \times Z \to Y \times Z$ is a closed map
for every topological space $Z$. Using a singleton for $Z$ we see that a proper map is closed.  We collect the elementary
properties of proper maps from \cite{B} Section 1.10.

\begin{prop}\label{appprop03} \begin{enumerate}
\item[(a)] If $f : X \to Y$ is injective, then it is proper iff it is a homeomorphism onto a closed subset of $Y$.

\item[(b)] Assume that  $f : X \to Y$ and $g : Y \to Z$ are continuous.
\begin{itemize}
\item[(i)] If $f$ and $g$ are proper, then $g \circ f$ is proper.
\item[(ii)]  If $g \circ f$ is proper and $f$ is surjective, then $g$ is proper.
\item[(iii)]  If $g \circ f$ is proper and $g$ is injective, then $f$ is proper.
 \end{itemize}

 \item[(c)] If $f_1 : X_1 \to Y_1$ and $f_2 : X_2 \to Y_2$ are continuous maps with $X_1$ and $X_2$ nonempty,
 then $f_1 \times f_2 : X_1 \times X_2 \to Y_1 \times Y_2$ is proper iff both $f_1$ and $f_2$ are proper.

 \item[(d)] Let $f : X \to Y$ be a proper map.  If $A$ is a closed subset of $X$ then the restriction
 $f|A : A \to Y$ is a proper map.

 \item[(e)] If $B$ is an arbitrary subset of $Y$ then the restriction
 $f|B : f^{-1}(B) \to B$ is a proper map.

 \item[(f)] If $f_1 : X \to Y_1$ and $f_2 : X \to Y_2$ are proper maps with $X$ Hausdorff
then the map $x \mapsto (f_1(x),f_2(x))$ is proper. In particular, its image is closed.
 \end{enumerate}
 \end{prop}

 {\bfseries Proof:} These results are Propositions 2-5 of \cite{B} Section 1.10.1.

 (a) An injective continuous map is a homeomorphism onto a closed subset iff it it is a closed map.

 (b) If $A \subset X \times Z$ then $[(g \circ f) \times 1_Z](A) = (g \times 1_Z)\circ (f \times 1_Z)(A)$ and
 if $g$ is injective, $(f \times 1_Z)(A) = (g  \times 1_Z)^{-1}((g \circ f) \times 1_Z](A)]$. If $f$ is surjective
 and $B \subset  Y \times Z$ then $(g \times 1_Z)(B) = [(g \circ f) \times 1_Z]((f  \times 1_Z)^{-1}(B))$.

 (c) $f_1 \times f_2$ is the composition  $(f_1 \times 1_{Y_2}) \circ (1_{X_1} \times f_2)$.

 (d) If $K$ is a closed subset of $A \times Z$ then $K$ is a closed subset of $X \times Z$.

(e) If $A \subset f^{-1}(B) \times Z$ is closed relative to $f^{-1}(B) \times Z$ then there exists $A_1$ closed in
 $X \times Z$ with $A = A_1 \cap (f^{-1}(B) \times Z)$ and $(f \times 1_Z)(A) = (f \times 1_Z)(A_1) \cap B \times Z$.

(f) Since $X$ is Hausdorff, the diagonal $1_X$ is closed in $X \times X$ and so the map $\Delta : X \to X \times X$
 $x \mapsto (x,x)$ is a proper map from $X $ to $X \times X$ by (a). Since $f_1 \times f_2$ is proper by (c),
 the composition $(f_1 \times f_2) \circ \Delta$ is proper.

$\Box$ \vspace{.5cm}

 The condition that $f$ be proper can be described in terms of compactness.  For convenience we restrict attention
 to Tychonoff spaces, i.e.\ completely regular Hausdorff spaces.

 \begin{prop}\label{appprop04} \begin{enumerate}
\item[(a)] Assume that  $f : X \to Y$ is continuous with X a Tychonoff space. The following are equivalent.
\begin{itemize}
\item[(i)] The map $f$ is proper.

\item[(ii)] $f \times 1_Z : X \times Z \to Y \times Z$ is a closed map
for every compact Hausdorff space $Z$.

\item[(iii)] The map $f$ is closed and $f^{-1}(y)$ is compact for every $y \in Y$.

\item[(iv)] Whenever $\{ x_i : i \in I \}$ is a net in $X$ such that $\{ f(x_i) \}$ converges to a point
$y \in Y$ then $\{ x_i \}$ has a cluster point in $f^{-1}(y)$.
\end{itemize}

\item[(b)] If $p$ is a singleton space and $X$ is a Tychonoff space, then the map
$p : X \to p$ is proper iff $X$ is compact.

\item[(c)] If  $f : X \to Y$ is proper with $X, Y$ Tychonoff spaces and $B \subset Y$
is compact, then $f^{-1}(B) \subset X$
is compact.
\end{enumerate}
\end{prop}

  {\bfseries Proof:} These results are essentially Theorem 1 and Lemma 1 of \cite{B} Section 1.10.2.

 (a) (i) $\Rightarrow$ (ii) Obvious.

 Let $Z$ be a compactification of $X$, i.e.\ there is a continuous embedding $k : X \to Z$ with
  $Z$ a compact Hausdorff space.
 Because $Z$ is Hausdorff, the map $k$ is a closed subset of $X \times Z$. The map
 $p \times 1_Z : X \times Z \to p \times Z$ is isomorphic to the projection
 $\pi_2: X \times Z \to Z$. If the map $p \times 1_Z$
 is closed, then $k(X) = \pi_2(k)$ is a closed subset of $Z$ and so is compact.
 Since $k$ is an embedding $X$ is compact. In particular,
 this proves one direction of (b).

 (ii) $\Rightarrow$ (iii) Using $Z$ as a singleton we see that $f$ is closed. As in Proposition \ref{appprop03}(e)
 we see that $f \times 1_Z : f^{-1}(y) \times Z \to y \times Z$ is closed for any compact Hausdorff space. From the
 above argument it follows that $f^{-1}(y)$ is compact.

 (iii) $\Rightarrow$ (iv) If for some $i \in I$ the set $A_i = \ol{\{ x_j : i \prec j  \}}$ is disjoint from
 $f^{-1}(y)$ then $f(A_i)$ is a closed set disjoint from $y$ and so $\{ f(x_i) \}$ does not converge to $y$.
 Hence, $\{ A_i \cap f^{-1}(y) \}$ is a collection of closed sets satisfying the finite intersection property.
 Since $f^{-1}(y)$ is compact, the intersection is nonempty and the intersection is the set of cluster points
 of $\{ x_i \}$ in $f^{-1}(y)$.

 (iv) $\Rightarrow$ (i) Let $A$ be a closed subset of $X \times Z$ and $(y,z)$ a point
 of the closure of $(f \times 1_Z)(A)$.
 There exists a net $\{ (x_i,z_i) \}$ in $A$ such that $\{ (f(x_i),z_i) \}$ converges to $(y,z)$. From (iv)
 it follows that there exists $x \in f^{-1}(y)$ and a subnet $\{ x_{i'} \}$ which converges to $x$.
 Hence, the subnet $\{ (x_{i'}, z_{i'}) \}$ converges to $(x,z)$ and since $A$ is closed $(x,z) \in A$.
 So $(y,z) = (f \times 1_Z)(x,z) \in (f \times 1_Z)(A)$. Thus $(f \times 1_Z)(A)$ is closed.

 (b) If $X$ is compact, then $X \to p$ satisfies condition (iii) of (a) and so is a proper  map.

 (c) Since $B$ is compact, $B \to p$ is proper. Since $f$ is proper, the restriction $f^{-1}(B) \to B$ is proper.
 Hence, the composition $f^{-1}(B) \to p$ is proper and so $f^{-1}(B)$ is compact.

$\Box$ \vspace{.5cm}

A Hausdorff space $X$ is called a \emph{k-space}\index{space!k-space}\index{k-space} when the topology is compactly generated.
That is, $A \cap K$ compact for
every compact subset $K$ of $X$ implies $A$ is closed. A locally compact space is clearly a k-space. Since a convergent
sequence together with its limit is compact, any Hausdorff sequential space is a k-space, where
$X$ is \emph{sequential} when $x \in \ol{A}$ implies $x$ is the limit of a sequence in $A$.
So any Hausdorff,  first countable
space is a k-space. In particular, a metrizable space is a k-space.

\begin{prop}\label{appprop05} Let $f : X \to Y$ be a continuous map with $X$ and $Y$ Tychonoff spaces.

(a) If $Y$ is a k-space and for every compact $B \subset Y$, the pre-image $f^{-1}(B)$ is
compact, then $f$ is a proper map.

(b) If $X$ is a k-space and $A \subset X$ such that the restriction $f|A : A \to Y$ is proper
then $A$ is a closed subset of $X$.
\end{prop}

{\bfseries Proof:} (a) From Proposition \ref{appprop04} (a)(iii) it suffices to
show that $f$ is closed. Let $A \subset X$ be
closed and let $K \subset Y$ be compact. By hypothesis, $f^{-1}(K)$ is compact and so $A \cap f^{-1}(K)$ is compact.
It follows that $f(A) \cap K = f(A \cap f^{-1}(K))$ is compact. As $K$ was arbitrary,
$f(A)$ is closed because $Y$ is a k-space.

(b) Let $K \subset X$ be compact so that $f(K) \subset Y$ is compact. By Proposition \ref{appprop04} (c) applied to
$f|A$, $(f|A)^{-1}(f(K))$ is compact. Hence, $K \cap A = K \cap (f|A)^{-1}(f(K))$ is compact.  Since $K$ was arbitrary and
$X$ is a k-space, $A$ is closed.

$\Box$

\vspace{1cm}

\bibliographystyle{amsplain}

\printindex

\end{document}